\documentclass[12pt,times]{elsarticle}

\usepackage{geometry}
\usepackage{graphicx}

\usepackage{epsfig}
\usepackage{epstopdf}
\usepackage{subfig}
\usepackage{enumitem}
\usepackage{multicol}
\usepackage[toc,page]{appendix}
\usepackage{amssymb}

\usepackage{mathtools}
\usepackage{amsmath}
\usepackage{algpseudocode}
\usepackage{algorithm}

\usepackage{siunitx}

\usepackage[pdfborder={0 0 0},colorlinks,allcolors=blue]{hyperref}
\usepackage{xcolor}

\usepackage{amsthm}
\theoremstyle{definition} 
\newtheorem{remark}{Remark} 
\newtheorem{example}{Example} 

\usepackage{cleveref}
\hypersetup{colorlinks=true, allcolors=blue,pdfauthor=author}
\crefname{equation}{Equation}{Equations}
\crefname{figure}{Figure}{Figures}
\crefname{section}{Section}{Sections}
\crefname{remark}{Remark}{Remarks}
\crefname{example}{Example}{Examples}
\crefname{algorithm}{Algorithm}{Algorithms}
\crefname{appendix}{}{}

\usepackage{lineno}

\usepackage{enumitem}
\newlist{steps}{enumerate}{1}
\setlist[steps, 1]{label = Step \arabic*:,leftmargin = 2cm}

\usepackage{xparse}

\NewDocumentCommand \an{ m }{%
    \langle {#1} \rangle%
}

\newcommand{\NAT}[1]{{\color{black}{#1}}}

\newcommand\notype[1]{\unskip}

\begin{document}

\begin{frontmatter}

\title{A General-Purpose, Inelastic, Rotation-Free Kirchhoff–Love Shell Formulation for Peridynamics}

\author[brownaddress]{Masoud Behzadinasab\corref{mycorrespondingauthor}}
\ead{masoud\_behzadinasab@brown.edu}

\author[brownaddress]{Mert Alaydin}

\author[sandiaaddress]{Nathaniel Trask}

\author[brownaddress]{Yuri Bazilevs\corref{mycorrespondingauthor}}
\ead{yuri\_bazilevs@brown.edu}

\cortext[mycorrespondingauthor]{Corresponding authors.}

\address[brownaddress]{School of Engineering, Brown University, 184 Hope St., Providence, RI 02912, USA}
\address[sandiaaddress]{Center for Computing Research, Sandia National Laboratories, Albuquerque, NM 87185, USA\textsuperscript{1}\footnote{\textsuperscript{1}Sandia National Laboratories is a multi-mission laboratory managed and operated by National Technology and Engineering Solutions of Sandia, LLC., a wholly owned subsidiary of Honeywell International, Inc., for the U.S. Department of Energy’s National Nuclear Security Administration under contract DE-NA0003525. This paper describes objective technical results and analysis. Any subjective views or opinions that might be expressed in the paper do not necessarily represent the views of the U.S. Department of Energy or the United States Government.}}

\begin{abstract}
We present a comprehensive rotation-free Kirchhoff–-Love (KL) shell formulation for peridynamics (PD) that is capable of modeling large elasto-plastic deformations and fracture in thin-walled structures. To remove the need for a predefined global parametric domain, Principal Component Analysis is employed in a meshfree setting to develop a local parameterization of the shell midsurface. The KL shell kinematics is utilized to develop a correspondence-based PD formulation. A bond-stabilization technique is employed to naturally achieve stability of the discrete solution. Only the mid-surface velocity degrees of freedom are used in the governing thin-shell equations. 3D rate-form material models are employed to enable simulating a wide range of material behavior. A bond-associative damage correspondence modeling approach is adopted to use classical failure criteria at the bond level, which readily enables the simulation of brittle and ductile fracture. \NAT{Discretizing the model with asymptotically compatible meshfree approximation provides a scheme which converges to the classical KL shell model while providing an accurate and flexible framework for treating fracture.} A wide range of numerical examples, ranging from elastostatics to problems involving plasticity, fracture, and fragmentation, are conducted to validate the accuracy, convergence, and robustness of the developed PD thin-shell formulation. It is also worth noting that the present method naturally enables the discretization of a shell theory requiring higher-order smoothness on a completely unstructured surface mesh.
\end{abstract}


\end{frontmatter}

\section{Introduction}
\label{sec:intro}

Low-weight components are widely used in the design of engineering structures, where it is preferable to keep the structural weight low by using high-strength material (e.g., aircraft). These structures - having through-thickness dimension much smaller than the other two - are referred to as {\em shell structures}~\cite{timoshenko1959theory}. Despite great progress made over the years to develop both structural shell theories and their discretization, challenges remain to efficiently and accurately predict the shell structural response under extreme loading conditions exhibiting large deformation and material failure using existing computational methods and tools. 

Formulations based upon the Reissner--Mindlin (RM) shell theory~\cite{reissner1945effect,mindlin1951influence} are widely implemented and used in commercial Finite Element Methods (FEM) software. The RM shell theory is suitable for thicker shells and includes through-thickness shearing deformations. The discrete equations of motion for the RM shell are written in terms of displacement and rotational degrees of freedom (DOFs) and are typically discretized with traditional $C^0$-continuous finite elements. While the KL shell formulation~\cite{love1888xvi} is a recent newcomer in the world of FEM, however, its popularity and adoption are soaring~\cite{kiendl2009isogeometric,kiendl2010bending} mainly due to the introduction of Isogeometric Analysis (IGA)~\cite{hughes2005isogeometric,cottrell2009isogeometric}. The KL shell formulation neglects transverse shear deformations and is thus suitable for describing thin structures. In the discrete setting, the KL equations of motion are formulated in terms of the midsurface displacement DOFs, without the need to introduce rotational DOFs. This efficiency gain comes with the need to use smooth basis functions (e.g., NURBS~\cite{cottrell2009isogeometric}), which is a distinguishing feature of IGA. It is also worth noting that the lack of rotational DOFs and transverse shear deformations in the KL shell formulation circumvents transverse-shear locking issues that need to be addressed for the RM shells.

Peridynamics (PD)~\cite{silling2000reformulation,silling2007peridynamic} has been in development in the past two decades. \NAT{It provides an integrodifferential model for solid mechanics which naturally accommodates both nonlocal physics and low-regularity solutions typical of fracture and interface problems. PD has been particularly effective in handling fracture mechanics problems involving:} 
concrete cracking~\cite{gerstle2007peridynamic,yaghoobi2017fracture,nikravesh2018improved,bobaru2018intraply}; crack branching phenomena~\cite{ha2010studies,bobaru2015cracks}; ductile fracture in metallic materials and structures~\cite{wu2015stabilized,behzadinasab2019third,hu2020numerical,boyce2014sandia,kramer2019third}; pitting corrosion damage~\cite{chen2015peridynamic,jafarzadeh2018peridynamic,rokkam2019nonlocal}; and fracture in porous materials~\cite{ouchi2015fully,behzadinasab2018peridynamics,chen2019peridynamic}, among others. An in-depth review of recent developments in PD may be found in~\cite{madenci2014peridynamic,bobaru2016handbook,javili2019peridynamics}. 

There are several PD formulations that were developed for thin structures, modeling membranes~\cite{silling2005peridynamic,sarego2016linearized,oterkus2020peridynamic} and plates and flat shells~\cite{o2014peridynamic,taylor2015two}. A number of PD theories, mostly based on RM kinematics, were proposed for thick-shell analysis~\cite{diyaroglu2015peridynamics,chowdhury2016peridynamic,nguyen2019peridynamics,vazic2020peridynamic,zhang2021peridynamic}, consequently involving computationally unfavorable rotational DOFs for thin-shell analysis. The existing KL-based PD models~\cite{o2014peridynamic,yang2020kirchhoff} were developed only for linearly elastic flat plates. In cases where simple curved geometries were considered (e.g., cylindrical or spherical surfaces), an approach based upon global, analytical parameterization of the shell was employed \cite{chowdhury2016peridynamic,nguyen2019peridynamics,zhang2021peridynamic}.

The correspondence principle~\cite{silling2007peridynamic} has been developed to enable the use of classical constitutive models within the PD framework. In correspondence-based PD, a deformation or strain measure (e.g., the deformation gradient) is computed using a PD estimator and used to evaluate the stress tensor using traditional constitutive laws. The stress state is then converted to the force state, which is more natural to PD. While the correspondence technique was initially found to suffer from instabilities~\cite{breitenfeld2014non,tupek2014extended,silling2017stability,chowdhury2019modified,behzadinasab2020stability}, several stabilization techniques have been proposed to circumvent these issues~\cite{littlewood2010simulation,silling2017stability,breitzman2018bond,chen2018bond,chowdhury2019modified,behzadinasab2020semi}. A higher-order, stable correspondence-based PD framework was recently presented in~\cite{behzadinasab2021unifiedI,behzadinasab2021unifiedII} and was shown to improve the performance of correspondence-based PD in both static and dynamic regimes. 

The objective of the present work is to develop a comprehensive correspondence-based PD formulation for the analysis of thin-shell structures undergoing large elasto-plastic deformation and damage in the static and dynamic regimes, \NAT{which asymptotically converges to the classical KL shell theory}. In addition, to enable complex-geometry computations, we no longer rely on the global parameterization of the mid-surface and utilize Principal Component Analysis (PCA)~\cite{wold1987principal} to define a local parameterization for the neighborhood of each PD node. The KL kinematic assumptions are employed to obtain the 3D spatial velocity gradient using only the mid-surface velocity unknowns. For this, we adopt the approach detailed in~\cite{benson2011large,alaydin2021updated} for IGA-based KL shells. As a result, the governing equations do not use rotational DOFs and is thus more computationally efficient. The mid-surface velocity gradient is calculated using a higher-order PD differential operator~\cite{madenci2016peridynamic,hillman2020generalized,behzadinasab2021unifiedI}. The rate-of-deformation, defined as the symmetrized spatial velocity gradient, is used to evolve the stress using 3D rate-form constitutive relations. The field variables (e.g., strain and stress) are used to incorporate classical fracture criteria at the bond level to naturally accommodate modeling of material failure, thus eliminating extra effort required to develop native PD shell constitutive models. 

\NAT{In the broader meshfree discretization literature, a number of works have developed discretizations of general surface PDE ~\cite{gross2020meshfree,trask2020compatible,lai2013local,shankar2018rbf,suchde2021meshfree,suchde2019fully,mohammadi2019generalized,torres2020approximation,shankar2018mesh,macdonald2013simple,fuselier2013high,piret2012orthogonal,leung2011grid}. Such schemes treat the manifold \textit{extrinsically}, working in ambient space and projecting to manifold, or \textit{intrinsically}, using compact reconstructions of the manifold to obtain local coordinates. The latter generally apply polynomial reconstruction to obtain estimates of the metric tensor necessary for computing surface differential operators ~\cite{lancaster1981surfaces,hoppe1992surface,amenta2004defining,gross2020meshfree,trask2020compatible,lai2013local}. In contrast, by working in the reference configuration, the linear PCA reconstruction used in the current approach provides high-order treatment of surface curvature, avoiding the need for either high-order polynomial reconstruction or use of greater than $C^0$ continuity shape functions. When the resulting model is discretized with an asymptotically compatible meshfree approximation ~\cite{tian2014asymptotically,hillman2020generalized,chi2013gradient,trask2019asymptotically,leng2021asymptotically,leng2020asymptotically}, we preserve high-order $\delta$-convergence \cite{bobaru2009convergence} to solutions of the local KL shell theory.}

The remainder of the paper is structured as follows. In \cref{sec:manifold}, we show how PCA is used to develop local parametric domains from a 3D point cloud. The basics of the KL shell kinematics, generalized from the IGA-based formulations in~\cite{benson2011large,alaydin2021updated}, are reviewed. In \cref{sec:PD_KL}, the rate form of the energy balance law is employed to derive the PD force state in the KL shell framework. Other important aspects of the formulation, i.e., the evaluation of the PD parametric gradients, the co-rotational formulation of standard J2 plasticity, the enforcement of the plane stress condition, the local thickness update, and continuum damage modeling, are presented in the same section. A battery of numerical tests, ranging from linear elastostatics to nonlinear dynamics with plasticity and fracture, are presented in \cref{sec:results} to demonstrate the accuracy, robustness, and the general-purpose nature of the proposed PD KL shell formulation. In \cref{sec:conclusions}, we provide concluding remarks and outline future research directions. 

In what follows, material points and PD bonds are indicated by parentheses and angle brackets by subscripts, respectively. For example, $\mathbf{X}_{({\rm P})}$ denotes the position of a material point P, and $\mathbf{T}_{\rm \an{P-Q}}$ is the force state for the bond ${\rm \an{P-Q}}$, which is a bond from P to its neighbor Q. The field variables varying along the thickness direction are denoted by the superscript ``3D'', and those that are defined on the mid-surface only are indicated by the superscript ``2D''.

\section{Shell Structures}
\label{sec:shell}

The main characteristic of shell structures is their small thickness relative to the other two dimensions. These structures are commonly described using a mid-surface, which is a 2D manifold in a 3D space, with a given distribution of thickness in the surface-normal direction. To achieve a good description of the mechanics of shell structures, the mid-surface geometry should be accurately represented and physically realistic through-thickness kinematic assumptions should be employed.

\subsection{Local Parameterization from a Point Cloud}
\label{sec:manifold}

\begin{figure}[!hbpt]
    \centering
    \includegraphics[width=0.45\columnwidth]{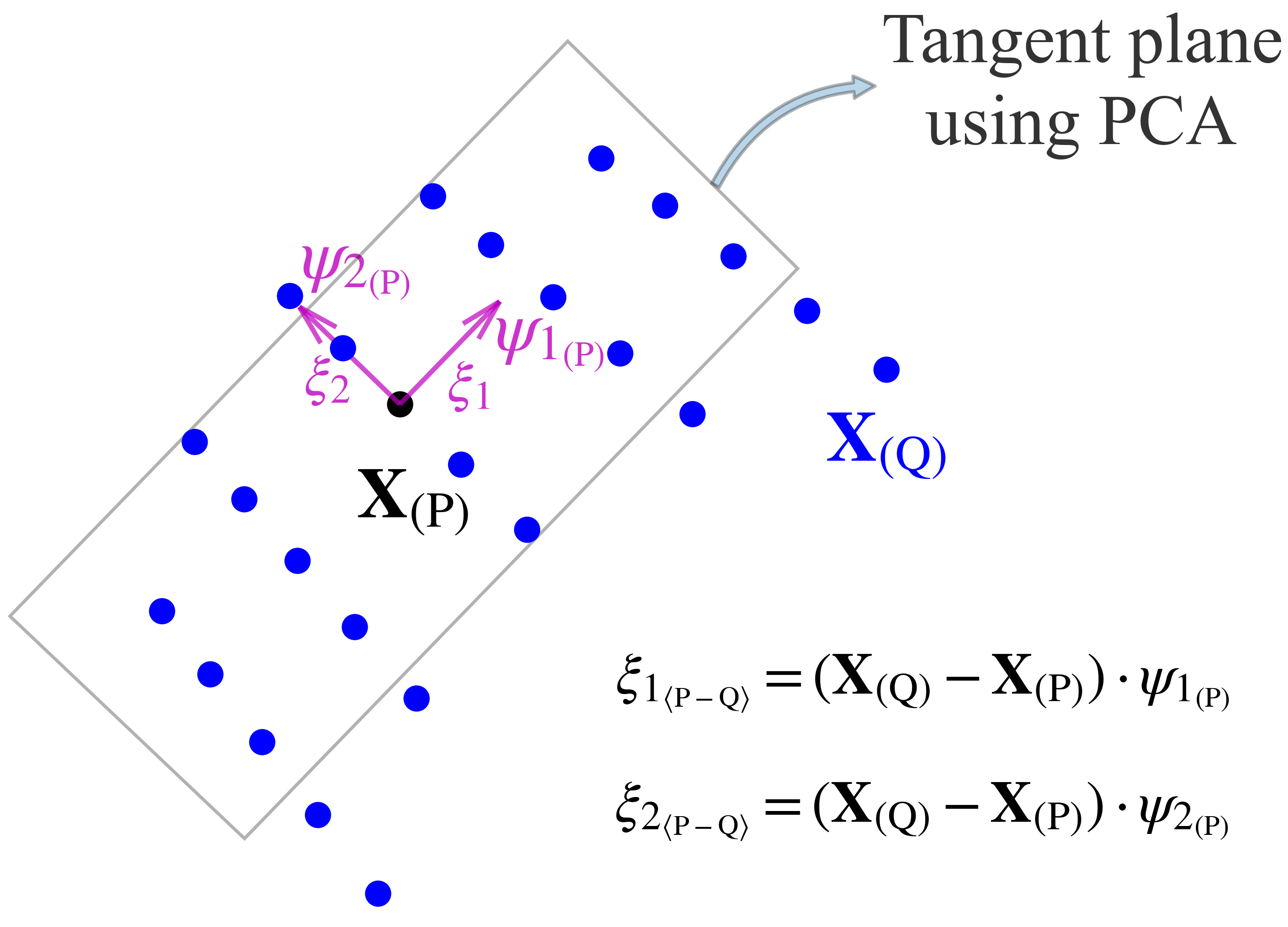}
    \caption{Schematic of the local parametric space from a point cloud using the PCA algorithm. First, PCA is used to obtain a plane approximately tangent to the point P using the reference positions of its neighborhood. Next, the neighbors Q are projected to the tangent space using the base vectors of the plane $\boldsymbol\psi_1$ and $\boldsymbol\psi_2$ to obtain the local parametrization.}
    \label{fig:PCA}
\end{figure}

We start with a meshfree, local parameterization of the shell mid-surface. We employ the approach detailed in~\cite{trask2020compatible} to define a local parameterization of a point cloud without resorting to using a global parametric map. Considering the undeformed configuration of a mid-surface $\mathcal{B}^{\rm 2D}$, we approximate a tangent plane at each surface mesh point ${\rm P} \in \mathcal{B}^{\rm 2D}$ using Principal Component Analysis (PCA)~\cite{wold1987principal}. In this approach, the positions of points in the local neighborhood of ${\rm P}$, i.e., ${\rm Q}\in\mathcal{H}_{({\rm P})}$, are used to estimate the local tangent plane. The centering point is defined as
\begin{equation} 
\bar{\mathbf{X}}_{({\rm P})} = \frac{1}{N_{({\rm P})}} \sum_{{\rm Q}\in \mathcal{H}_{({\rm P})}} \mathbf{X}_{({\rm Q})},
\end{equation}
where $\mathbf{X}$ is the reference-configuration position vector, and $N_{({\rm P})}$ is the number of points in the local neighborhood. The covariance matrix for the neighborhood of ${\rm P}$ is given by
\begin{equation} 
\pmb{\mathcal{C}}_{({\rm P})} = \text{Cov}(\{{\rm P}\}) = \frac{1}{N_{({\rm P})}} \sum_{{\rm Q}\in \mathcal{H}_{({\rm P})}} \left( \mathbf{X}_{({\rm Q})} - \bar{\mathbf{X}}_{({\rm P})} \right) \left( \mathbf{X}_{({\rm Q})} - \bar{\mathbf{X}}_{({\rm P})} \right)^{\intercal}.
\end{equation}
The eigenvectors corresponding to the two largest eigenvalues of the covariance matrix $\pmb{\mathcal{C}}_{({\rm P})}$ are orthogonal and define a good approximation to the tangent plane at ${\rm P}$. After normalizing, we denote the orthonormal bases of the tangent plane at ${\rm P}$ by $\boldsymbol\psi_{1_{({\rm P})}}$ and $\boldsymbol\psi_{2_{({\rm P})}}$. This, in turn, allows us to define the \textit{local} parametric coordinates for each neighbor ${\rm Q}$ as
\begin{equation}
\begin{aligned}
\xi_{1_{\rm \an{P-Q}}} &= \left( \mathbf{X}_{({\rm Q})} - \mathbf{X}_{({\rm P})} \right) \cdot \boldsymbol\psi_{1_{({\rm P})}} , \\
\xi_{2_{\rm \an{P-Q}}} &= \left( \mathbf{X}_{({\rm Q})} - \mathbf{X}_{({\rm P})} \right) \cdot \boldsymbol\psi_{2_{({\rm P})}} .
\end{aligned}
\label{eqn:xi12}
\end{equation}
A schematic of the process for obtaining the local parametric space from a given point set is provided in \cref{fig:PCA}.
\begin{remark}
It is important to note that the tangent plane at ${\rm P}$ defined above is not thought of as the local approximation of the mid-surface, but rather as the \textit{local parametric domain} for the mid-surface, much like in the isoparametric FEM or IGA approaches. This is in contrast to the local surface reconstruction approaches for the discretization of PDEs on surfaces that enhance the tangent plane with an estimate of the local curvature to improve solution accuracy (see, e.g.,~\cite{trask2020compatible,gross2020meshfree,lai2013local}.) 
\end{remark}
\begin{remark}
Note that the above approach constructs a local parametric domain for each ${\rm P}$ in the reference configuration, which remains unchanged throughout the deformation. This is akin to the Lagrangian formulation typically employed in FEM or IGA. A semi-Lagrangian approach, where the local parameterization is periodically recomputed based on the current coordinates of the neighbors as the structure deforms, could also be adopted~\cite{shankar2018mesh,guan2011semi,behzadinasab2020semi}); we leave this development for the future work.
\end{remark}

\subsection{Thin Shell Kinematics}
\label{sec:kinematics}

Equipped with the local parameterization of the mid-surface, in this section we adopt the KL shell kinematics from the formulation presented in~\cite{benson2011large,alaydin2021updated} in the context of IGA. 

\begin{figure}[!hbpt]
    \centering
    \includegraphics[width=0.95\columnwidth]{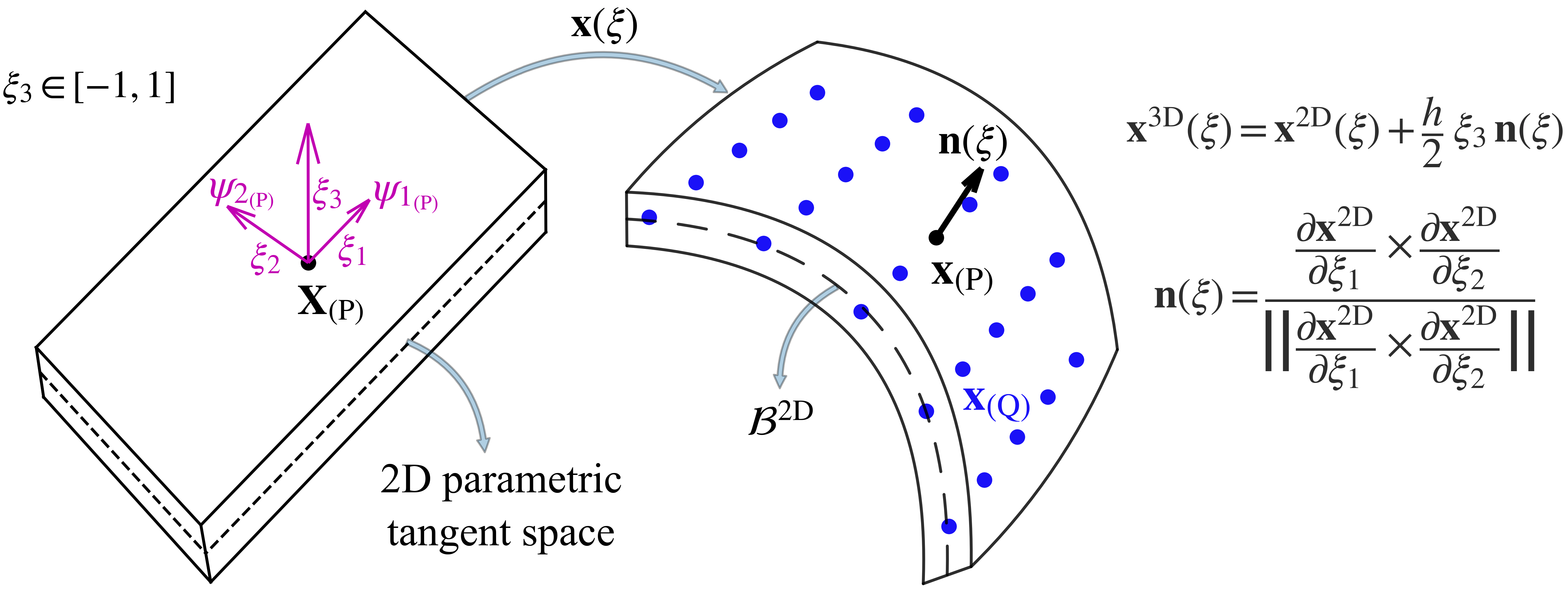}
    \caption{Schematic of the parametrization enhanced with the KL shell kinematics.}
    \label{fig:KL}
\end{figure}

In addition to the in-plane parametric coordinates $\xi_1, \xi_2$, we denote by $\xi_3 \in [-1, +1]$ the through-thickness parametric coordinate. We define the 3D position vector $\mathbf{x}^{\rm 3D}$ in the current configuration as
\begin{equation}
\mathbf{x}^{\rm 3D}(\xi_1,\xi_2,\xi_3) = \mathbf{x}^{\rm 2D}(\xi_1,\xi_2) + \frac{h}{2} \xi_3 \, \mathbf{n}(\xi_1,\xi_2) ,
\label{eqn:x3d}
\end{equation}
where $\mathbf{x}^{\rm 2D}(\xi_1,\xi_2)$ is the mid-surface position vector, $h$ is the local shell thickness, and $\mathbf{n}(\xi_1,\xi_2)$ is the unit surface normal vector, all defined in the current configuration. The normal vector may be expressed as
\begin{equation}
\mathbf{n}(\xi_1,\xi_2) = \frac{ \dfrac{\partial \mathbf{x}^{\rm 2D}}{\partial \xi_1} \times \dfrac{\partial \mathbf{x}^{\rm 2D}}{\partial \xi_2} }{\left|\left| \dfrac{\partial \mathbf{x}^{\rm 2D}}{\partial \xi_1} \times \dfrac{\partial \mathbf{x}^{\rm 2D}}{\partial \xi_2} \right|\right|},
\label{eqn:normal}
\end{equation}
using a suitable definition of the parametric derivatives, which we present in later sections. A schematic of the KL shell kinematics with respect to the parametric space is provided in \cref{fig:KL}.
Taking the material time derivative of the position vector and dropping functional dependence on $\xi_{i=1,2}$ for brevity, the velocity vector $\mathbf{v}^{\rm 3D}$ may be expressed in index notation as
\begin{equation}
v_i^{\rm 3D} = v_i^{\rm 2D} + \frac{h}{2} \xi_3 \, \dot{n}_i,
\label{eqn:v3d}
\end{equation}
where
\begin{equation}
\dot{\mathbf{n}} = \frac{ \mathbf{I}-\mathbf{n}\otimes\mathbf{n} }{\left|\left| \dfrac{\partial \mathbf{x}^{\rm 2D}}{\partial \xi_1} \times \dfrac{\partial \mathbf{x}^{\rm 2D}}{\partial \xi_2} \right|\right|} \left( \frac{\partial \mathbf{v}^{\rm 2D}}{\partial \xi_1} \times \frac{\partial \mathbf{x}^{\rm 2D}}{\partial \xi_2} + \frac{\partial \mathbf{x}^{\rm 2D}}{\partial \xi_1} \times \frac{\partial \mathbf{v}^{\rm 2D}}{\partial \xi_2} \right).
\label{eqn:ndotFull}
\end{equation}
Denoting the first fraction on the right-hand side of \cref{eqn:ndotFull} as $\mathbf{A}$, 
\begin{equation}
\mathbf{A} = \frac{ \mathbf{I}-\mathbf{n}\otimes\mathbf{n} }{\left|\left| \dfrac{\partial \mathbf{x}^{\rm 2D}}{\partial \xi_1} \times \dfrac{\partial \mathbf{x}^{\rm 2D}}{\partial \xi_2} \right|\right|} \, ,
\label{eqn:A}
\end{equation}
$\dot{\mathbf{n}}$ can be written in index notation as
\begin{equation}
\dot{n}_i = B^1_{ik} \, v^{\rm 2D}_{k,\xi_1} + B^2_{il} \, v^{\rm 2D}_{l,\xi_2},
\label{eqn:ndot}
\end{equation}
where the auxiliary second-order tensors $\mathbf{B}^1$ and $\mathbf{B}^2$ are defined as
\begin{equation}
\begin{aligned}
& B^1_{ik} = A_{ij} \, \epsilon_{jkl} \, x^{\rm 2D}_{l,\xi_2} , \\
& B^2_{il} = A_{ij} \, \epsilon_{jkl} \, x^{\rm 2D}_{k,\xi_1} ,
\end{aligned}
\label{eqn:B12}
\end{equation}
and $\epsilon_{ijk}$ is the alternator or Levi–Civita tensor.

Using \cref{eqn:v3d,eqn:ndot}, the full 3D velocity vector may be expressed as
\begin{equation}
v_i^{\rm 3D} = v_i^{\rm 2D} + \frac{h}{2} \xi_3 \left( B^1_{ik} \, v^{\rm 2D}_{k,\xi_1} + B^2_{il} \, v^{\rm 2D}_{l,\xi_2} \right).
\label{eqn:v3dind}
\end{equation}
Taking the partial derivatives of $\mathbf{v}^{\rm 3D}$ with respect to the in-plane and through-thickness parametric coordinates, we obtain
\begin{equation}
\begin{aligned}
v^{\rm 3D}_{i,\xi_j} &= v^{\rm 2D}_{i,\xi_j} + \frac{h}{2} \xi_3 \left( B^1_{ik,\xi_j} \, v^{\rm 2D}_{k,\xi_1} + B^1_{ik} \, v^{\rm 2D}_{k,\xi_1\xi_j} + B^2_{il,\xi_j} \, v^{\rm 2D}_{l,\xi_2} + B^2_{il} \, v^{\rm 2D}_{l,\xi_2\xi_j} \right), \quad j=1,2 \  (\text{in-plane}), \\
v^{\rm 3D}_{i,\xi_j} &= \frac{h}{2} \left( B_{ik}^1 \, v^{\rm 2D}_{k,\xi_1} + B_{il}^2 \, v^{\rm 2D}_{l,\xi_2} \right), \quad j=3 \  (\text{through-thickness}),
\end{aligned}
\label{eqn:v3dd}
\end{equation}
where, using \cref{eqn:B12}, we can express
\begin{equation}
\begin{aligned}
& B^1_{ik,\xi_j} = \epsilon_{pkl} \left( A_{ip,\xi_j} \, x^{\rm 2D}_{l,\xi_2} + A_{ip} \, x^{\rm 2D}_{l,\xi_2\xi_j} \right) , \\
& B^2_{il,\xi_j} = \epsilon_{pkl} \left( A_{ip,\xi_j} \, x^{\rm 2D}_{k,\xi_1} +  A_{ip} \, x^{\rm 2D}_{k,\xi_1\xi_j} \right) .
\end{aligned}
\label{eqn:B12d}
\end{equation}
Using \cref{eqn:A}, $A_{ip,\xi_j}$ is computed as
\begin{equation}
A_{ip,\xi_j} = - \, \epsilon_{klm} \frac{ A_{ik} \, n_p + A_{pk} \, n_i + A_{ip} \, n_k }{\left|\left| \dfrac{\partial \mathbf{x}^{\rm 2D}}{\partial \xi_1} \times \dfrac{\partial \mathbf{x}^{\rm 2D}}{\partial \xi_2} \right|\right|} \left( x^{\rm 2D}_{l,\xi_1\xi_j} \, x^{\rm 2D}_{m,\xi_2} + x^{\rm 2D}_{l,\xi_1} \, x^{\rm 2D}_{m,\xi_2\xi_j} \right) .
\label{eqn:dAdxi}
\end{equation}
Partial derivatives of the position vector $\mathbf{x}^{\rm 3D}$ with respect to the parametric coordinates are computed in a similar fashion,
\begin{equation}
\begin{aligned}
x^{\rm 3D}_{i,\xi_j} &= x^{\rm 2D}_{i,\xi_j} + \frac{h}{2} \xi_3 \left( B^1_{ik} \, x^{\rm 2D}_{k,\xi_1\xi_j} + B^2_{il} \, x^{\rm 2D}_{l,\xi_2\xi_j} \right), \quad j=1,2 \  (\text{in-plane}), \\
x^{\rm 3D}_{i,\xi_j} &= \frac{h}{2} n_i, \quad j=3 \  (\text{through-thickness}).
\end{aligned}
\label{eqn:x3dd}
\end{equation}
Finally, using the chain rule, the spatial velocity gradient tensor in index notation may be expressed as
\begin{equation}
v^{\rm 3D}_{i,j} = v^{\rm 3D}_{i,\xi_k} F^{\rm 3D^{-1}}_{kj}, 
\label{eqn:transferV}
\end{equation}
or in matrix notation as
\begin{equation}
\begin{aligned}
\nabla \mathbf{v}^{\rm 3D} & = 
\begin{bmatrix}
\uparrow & \uparrow & \uparrow \\
\dfrac{\partial \mathbf{v}^{\rm 3D}}{\partial \xi_1} & \dfrac{\partial \mathbf{v}^{\rm 3D}}{\partial \xi_2} & \dfrac{h}{2}\dot{\mathbf{n}} \\
\uparrow & \uparrow & \uparrow
\end{bmatrix} 
\cdot
\mathbf{F}^{\rm 3D^{-1}} ,
\end{aligned}
\label{eqn:defDradV}
\end{equation}
where $F^{\rm 3D}_{kj} = x^{\rm 3D}_{k,\xi_j}$ is the Jacobian of the mapping between the local parametric domain and the physical domain of the current configuration,
\begin{equation}
\begin{aligned}
\mathbf{F}^{\rm 3D} & = 
\begin{bmatrix}
\uparrow & \uparrow & \uparrow \\
\dfrac{\partial \mathbf{x}^{\rm 3D}}{\partial \xi_1} & \dfrac{\partial \mathbf{x}^{\rm 3D}}{\partial \xi_2} & \dfrac{h}{2}\mathbf{n} \\
\uparrow & \uparrow & \uparrow
\end{bmatrix}.
\end{aligned}
\label{eqn:defF}
\end{equation}

\section{Peridynamic KL Shell Formulation}
\label{sec:PD_KL}

In this section, the mathematical formulation for the proposed PD shell framework is detailed. 

\subsection{Derivation from the Energy Balance Law}
\label{sec:balance_laws}

Let $\mathcal{B}^{\rm 3D}$ denote a body in the undeformed configuration. The rate form of the energy balance law states that the rate of change of total energy of the body is balanced by the supplied power $\mathcal{W}_{\rm supp_{(\mathcal{B})}}$. Neglecting the thermal effects, the rate-form energy balance law may be expressed as:
\begin{equation} 
\dot{\mathcal{K}}_{(\mathcal{B})} + \dot{\mathcal{U}}_{(\mathcal{B})} = \mathcal{W}_{\rm supp_{(\mathcal{B})}} .
\label{eqn:energy_balace}
\end{equation}
Here, $\mathcal{K}_{(\mathcal{B})}$ is the total kinetic energy
\begin{equation} 
\begin{aligned}
\mathcal{K}_{(\mathcal{B})} = \int_{\mathcal{B}^{\rm 3D}} \frac{1}{2} \, \rho_{(\rm P)} \, \mathbf{v}_{(\rm P)} \cdot \mathbf{v}_{(\rm P)} \, {\rm dP} ,
\end{aligned}
\label{eqn:kinetic_E}
\end{equation}
where ${\rm dP}$ is the infinitesimal volume of the material point ${\rm P}$, $\rho_{(\rm P)}$ is the mass density of ${\rm P}$ in the undeformed configuration, and $\mathbf{v}_{(\rm P)}$ is the velocity of ${\rm P}$.

The total strain power $\dot{\mathcal{U}}_{(\mathcal{B})}$ may be expressed as
\begin{equation}   
\dot{\mathcal{U}}_{(\mathcal{B})} = \int_{\mathcal{B}^{\rm 3D}} \dot{\mathcal{U}}_{(\rm P)} \, {\rm dP},
\label{eqn:Udef}
\end{equation}
where $\dot{\mathcal{U}}_{(\rm P)}$ is the strain-power density at P. The PD theory regularizes $\dot{\mathcal{U}}_{(\rm P)}$ as~\cite{silling2010peridynamic,behzadinasab2020peridynamic}
\begin{equation}   
\dot{\mathcal{U}}_{(\rm P)} = \int_{\mathcal{B}^{\rm 3D}} \alpha_{\rm \an{P-Q}} \, {\dot{\mathcal{U}}}_{\rm \an{P-Q}} \, {\rm dQ},
\label{eqn:U_PD}
\end{equation}
where ${\dot{\mathcal{U}}}_{\rm \an{P-Q}}$ is the strain-power density associated with a bond ${\rm \an{P-Q}}$, and $\alpha_{\rm \an{P-Q}}$ is a normalized weighting function constrained to satisfy the identity
\begin{equation}
\int_{\mathcal{B}^{\rm 3D}} \alpha_{\rm \an{P-Q}} \, {\rm dQ} = 1.
\end{equation}
The bond-associated strain-power density is defined as:
\begin{equation} 
{\dot{\mathcal{U}}}_{\rm \an{P-Q}} = \boldsymbol{\tau}_{\rm \an{P-Q}} : \mathbf{L}_{\rm \an{P-Q}} ,
\label{eqn:U_classical}
\end{equation}
where $\boldsymbol{\tau}_{\rm \an{P-Q}}$ and $\mathbf{L}_{\rm \an{P-Q}}$ are the bond-associated, power-conjugate Kirchhoff stress and velocity gradient tensors, respectively. Combining \cref{eqn:Udef,eqn:U_PD,eqn:U_classical}, we obtain
\begin{equation}   
\dot{\mathcal{U}}_{(\mathcal{B})} = \int_{\mathcal{B}^{\rm 3D}} \int_{\mathcal{B}^{\rm 3D}} \alpha_{\rm \an{P-Q}} \ \boldsymbol{\tau}_{\rm \an{P-Q}} : \mathbf{L}_{\rm \an{P-Q}} \, {\rm dQ} \, {\rm dP}.
\label{eqn:U_stress_based}
\end{equation}

The above equations have general applicability to any 3D continuum body. In order to specialize the formulation to a thin shell, we need to appropriately restrict the kinematics to respect the Kirchhoff hypothesis, which states that the shell director remains orthogonal to the mid-surface throughout the deformation. For this, let $\mathcal{B}^{\rm 2D}$ denote the shell mid-surface in the reference configuration. Introducing the Kirchhoff--Love shell kinematics given by \cref{eqn:v3d} into \cref{eqn:kinetic_E} and integrating through the thickness at each material point ${\rm P}$, we obtain the rate of kinetic energy per unit area as
\begin{equation} 
\begin{aligned}
\dot{\mathcal{K}}_{(\mathcal{P})} &= \int_{-1}^{1} \rho \, \dot{\mathbf{v}}^{\rm 3D}_{(\rm P)} \cdot \mathbf{v}^{\rm 3D}_{(\rm P)} \frac{h_{\rm (P)}}{2} \, {\rm d}\xi_3 \\
&= \frac{\rho^{\rm 2D}_{(\rm P)}}{2} \int_{-1}^{1} \left( \dot{\mathbf{v}}^{\rm 2D}_{(\rm P)} + \frac{h_{\rm (P)}}{2} \xi_3 \, \ddot{\mathbf{n}}_{\rm (P)} \right) \cdot \left( \mathbf{v}^{\rm 2D}_{(\rm P)} + \frac{h_{\rm (P)}}{2} \xi_3 \, \dot{\mathbf{n}}_{\rm (P)} \right) \, {\rm d}\xi_3 \\
&= \rho^{\rm 2D}_{(\rm P)} \, \dot{\mathbf{v}}^{\rm 2D}_{(\rm P)} \cdot \mathbf{v}^{\rm 2D}_{(\rm P)} + \frac{\rho^{\rm 2D}_{(\rm P)} \, h^2_{\rm (P)}}{12} \, \ddot{\mathbf{n}}_{\rm (P)} \cdot \dot{\mathbf{n}}_{\rm (P)} ,
\end{aligned}
\label{eqn:K_P}
\end{equation}
where $\rho^{\rm 2D}_{(\rm P)}$ is the body mass per unit area at ${\rm P}$ given by
\begin{equation} 
\rho^{\rm 2D}_{(\rm P)} = \rho \, h_{\rm (P)} .
\label{eqn:rho_2D}
\end{equation}
Integrating over the shell undeformed surface results in the following expression for the kinetic energy rate in terms of the translational and rotational contributions:
\begin{equation} 
\begin{aligned}
\dot{\mathcal{K}}_{(\mathcal{B})} &= \dot{\mathcal{K}}^{\rm tran}_{(\mathcal{B})} + \dot{\mathcal{K}}^{\rm rot}_{(\mathcal{B})} \\
& = \int_{\mathcal{B}^{\rm 2D}} \rho^{\rm 2D}_{(\rm P)} \, \dot{\mathbf{v}}^{\rm 2D}_{(\rm P)} \cdot \mathbf{v}^{\rm 2D}_{(\rm P)}~{\rm dP} + \int_{\mathcal{B}^{\rm 2D}} \frac{\rho^{\rm 2D}_{(\rm P)} \, h^2_{\rm (P)}}{12} \, \ddot{\mathbf{n}}_{\rm (P)} \cdot \dot{\mathbf{n}}_{\rm (P)}~{\rm dP}.
\end{aligned}
\label{eqn:kinetic_E_trans_rot}
\end{equation}
Starting from \cref{eqn:U_stress_based}, the stress power terms may be re-expressed as follows:
\begin{equation}   
\dot{\mathcal{U}}_{(\mathcal{B})} = \int_{\mathcal{B}^{\rm 3D}} \int_{\mathcal{B}^{\rm 3D}} \alpha_{\rm \an{P-Q}} \ \boldsymbol{\tau}_{\rm \an{P-Q}} : \mathbf{L}_{\rm \an{P-Q}} \, {\rm dQ} \, {\rm dP} = \int_{\mathcal{B}^{\rm 2D}} \int_{\mathcal{B}^{\rm 2D}} \mathbf{T}^{\rm int}_{\rm \an{P-Q}} \cdot \mathbf{v}^{\rm 2D}_{\rm \an{P-Q}} \, {\rm dQ} \, {\rm dP},
\label{eqn:U_T_based_1}
\end{equation}
where $\mathbf{T}^{\rm int}$ is the PD force state (with units of force per unit area squared in the shell framework) and $\mathbf{v}^{\rm 2D}_{\rm \an{P-Q}}$ is the mid-surface velocity state defined as
\begin{equation}   
\mathbf{v}^{\rm 2D}_{\rm \an{P-Q}} = \mathbf{v}^{\rm 2D}_{(\rm Q)} - \mathbf{v}^{\rm 2D}_{(\rm P)} .
\label{eqn:v2D_state}
\end{equation}
The expression for the velocity gradient in \cref{eqn:defDradV,eqn:defF}, which depends on $\mathbf{v}^{\rm 2D}$, determines the exact relationship between the force state and stress tensor in \cref{eqn:U_T_based_1}, which is to be defined in \cref{sec:PDforce}. Using \cref{eqn:v2D_state} and changing the order of integration in \cref{eqn:U_T_based_1} results in
\begin{equation}   
\dot{\mathcal{U}}_{(\mathcal{B})} = \int_{\mathcal{B}^{\rm 2D}} \left( \int_{\mathcal{B}^{\rm 2D}} \left( \mathbf{T}^{\rm int}_{\rm \an{Q-P}} - \mathbf{T}^{\rm int}_{\rm \an{P-Q}} \right) \, {\rm dQ} \right) \cdot \mathbf{v}^{\rm 2D}_{(\rm P)} \, {\rm dP}.
\label{eqn:U_T_based_2}
\end{equation}

Noting that $\dot{\mathbf{n}}$ is a function of $\mathbf{v}^{\rm 2D}$ as per \cref{eqn:ndot}, the rate of the rotational kinetic energy may also be re-expressed as
\begin{equation}
\begin{aligned}
\dot{\mathcal{K}}^{\rm rot}_{(\mathcal{P})} = 
\int_{\mathcal{B}^{\rm 2D}} \frac{\rho^{\rm 2D}_{(\rm P)} \, h^2_{\rm (P)}}{12} \, \ddot{\mathbf{n}}_{\rm (P)} \cdot \dot{\mathbf{n}}_{\rm (P)} \, {\rm dP} = \int_{\mathcal{B}^{\rm 2D}} \left( \int_{\mathcal{B}^{\rm 2D}} \left( \mathbf{T}^{\rm rot}_{\rm \an{Q-P}} - \mathbf{T}^{\rm rot}_{\rm \an{P-Q}} \right) \, {\rm dQ} \right) \cdot \mathbf{v}^{\rm 2D}_{(\rm P)} \, {\rm dP} ,
\end{aligned}
\label{eqn:T_rot}
\end{equation}
where the detailed expression for the corresponding force state $\mathbf{T}^{\rm rot}$ is also provided in \cref{sec:PDforce}.

Combining \cref{eqn:energy_balace,eqn:kinetic_E_trans_rot,eqn:K_P,eqn:U_T_based_2,eqn:T_rot}, we obtain
\begin{equation}
\int_{\mathcal{B}^{\rm 2D}} \left[ \rho^{\rm 2D}_{(\rm P)} \, \dot{\mathbf{v}}^{\rm 2D}_{(\rm P)} + \int_{\mathcal{B}^{\rm 2D}}\left( \mathbf{T}_{\rm \an{Q-P}} - \mathbf{T}_{\rm \an{P-Q}} \right) {\rm dQ} - \rho^{\rm 2D}_{(\rm P)}~\mathbf{b}_{(\rm P)} \right] \cdot \mathbf{v}^{\rm 2D}_{(\rm P)} \, {\rm dP} = 0 ,
\label{eq:WeakForm}
\end{equation}
where
\begin{equation}
\mathbf{T}_{\rm \an{P-Q}} = \mathbf{T}^{\rm int}_{\rm \an{P-Q}} + \mathbf{T}^{\rm rot}_{\rm \an{P-Q}}
\label{eqn:T_tot}
\end{equation}
and $\rho^{\rm 2D}_{(\rm P)}~\mathbf{b}_{(\rm P)}$ is the body force per unit area. Using a localization argument that \cref{eq:WeakForm} must hold \textit{for all} choices of $\mathbf{v}^{\rm 2D}_{(\rm P)}$ and at each material point ${\rm P}$, we arrive at a PD version of the local balance of linear momentum:
\begin{equation} 
\rho^{\rm 2D}_{(\rm P)} \, \dot{\mathbf{v}}^{\rm 2D}_{(\rm P)} = \int_{\mathcal{B}^{\rm 2D}} \left( \mathbf{T}_{\rm \an{P-Q}} - \mathbf{T}_{\rm \an{Q-P}} \right) \, {\rm dQ} + \rho^{\rm 2D}_{(\rm P)}~\mathbf{b}_{(\rm P)}.
\label{eqn:PD_EOM}
\end{equation}
\cref{eqn:PD_EOM} represents the state-based PD formulation~\cite{silling2007peridynamic} of the equation of motion for the KL shell structure. In the discrete setting, only a finite number of material points ${\rm P}$ are employed, each corresponding to a mesh node, and the integral over ${\mathcal{B}^{\rm 2D}}$ is carried out using nodal quadrature equipped with appropriate corrections (i.e., asymptotically compatible schemes).

\subsection{Evaluation of the Parametric Derivatives}
\label{sec:PD_gradient}

The spatial derivatives in meshfree methods, including PD, are typically computed directly from the point cloud with the aid of higher-order kernel functions  and/or numerical integration~\cite{chi2013gradient,madenci2016peridynamic,hillman2020generalized,behzadinasab2021unifiedI,trask2019asymptotically}. In the proposed thin shell framework, however, the spatial derivatives are computed using a chain rule that involves a mapping from the parametric to the physical domain, as outlined in the previous sections. The chain rule involves computation of the parametric derivatives at each material point ${\rm P}$, which is done as follows. In the continuous form, an integral-based first-order partial derivative operator acting on a vector-valued function $\mathbf{f}$ is defined as
\begin{equation} 
\mathbf{f}_{{,\xi_j}_{({\rm P})}} = \int_{\mathcal{B}^{\rm 2D}} \mathbf{f}_{\rm \an{P-Q}} \, \Phi_{\xi_j} \left(\xi_{1_{\rm \an{P-Q}}}, \xi_{2_{\rm \an{P-Q}}}\right) \, {\rm dQ}, \qquad j=1,2 .
\label{eqn:grad1}
\end{equation}
In practice, the discrete form of the partial-derivative operator is obtained by applying nodal quadrature in the above expression, namely,
\begin{equation} 
\mathbf{f}_{{,\xi_j}_{({\rm P})}} = \sum_{{\rm Q}\in \mathcal{B}^{\rm 2D}} \mathbf{f}_{\rm \an{P-Q}} \, \Phi_{\xi_j} \left(\xi_{1_{\rm \an{P-Q}}}, \xi_{2_{\rm \an{P-Q}}}\right) A_{\rm (Q)}, \qquad j=1,2 .
\label{eqn:grad1disc}
\end{equation}
Here,
\begin{equation} 
\mathbf{f}_{\rm \an{P-Q}} = \mathbf{f}_{({\rm Q})} - \mathbf{f}_{({\rm P})} ,
\end{equation}
and $\Phi\left(\xi_1, \xi_2\right)$ is the kernel function of two in-plane parametric variables that can be defined explicitly as in the reproducing kernel (RK) methods~\cite{chi2013gradient,hillman2020generalized} or using alternative algorithms as in the generalized moving least squares (GLMS) methods~\cite{trask2019asymptotically,leng2021asymptotically}. As shown in~\cite{behzadinasab2021unifiedI}, similar performance can be generally expected from both approaches. Because the KL shell kinematics involves both first- and second-order in-plane parametric derivatives, we also calculate the second parametric derivatives of $\mathbf{f}$ as
\begin{equation} 
\mathbf{f}_{{,\xi_j\xi_k}_{({\rm P})}} = \sum_{{\rm Q}\in \mathcal{B}^{\rm 2D}} \mathbf{f}_{\rm \an{P-Q}} \Phi_{\xi_j\xi_k} \left(\xi_{1_{\rm \an{P-Q}}}, \xi_{2_{\rm \an{P-Q}}}\right) A_{\rm (Q)}, \qquad j=1,2, \quad k=1,2 .
\label{eqn:grad2}
\end{equation}
The RK-based kernel functions and their derivatives~\cite{hillman2020generalized} (which is closely related to the RKPM implicit gradient~\cite{chi2013gradient} and the PD differential operator~\cite{madenci2016peridynamic}) are constructed based on the Taylor series expansion by including the parametric-domain monomials up to order $n$ as
\begin{equation}
  \boldsymbol{\Phi}_{\rm \an{P-Q}}^{[n]} = \omega_{\rm \an{P-Q}} \, \tilde{\boldsymbol{\Xi}}^{[n]} \, \mathbf{M}^{[n]^{-1}}_{({\rm P})} \, \boldsymbol{\Xi}^{[n]}_{\rm \an{P-Q}} ,
\label{eqn:kernel_vector}  
\end{equation}
where the vector $\boldsymbol{\Phi}_{\rm \an{P-Q}}^{[n]}$ has the following structure:
\begin{equation}
  \boldsymbol{\Phi}_{\rm \an{P-Q}}^{[n]} = \left[\Phi_{\xi_1}, \, \Phi_{\xi_2}, \, \Phi_{\xi_1\xi_1}, \, \Phi_{\xi_1\xi_2}, \, \Phi_{\xi_2\xi_2}\right]_{\rm \an{P-Q}}^\intercal ,
\end{equation}
$\omega_{\rm \an{P-Q}}$ is the scalar weighting function that depends on the relative distance between the material points normalized by the PD support size or {\em horizon} $\delta$ (note, $\omega_{\rm \an{P-Q}}=0$ for the bonds outside of the horizon), and $\boldsymbol{\Xi}^{[n]}_{\rm \an{P-Q}}$ is a column vector of the parametric-domain monomials. The discrete form of the moment matrix $\mathbf{M}^{[n]}$ in \cref{eqn:kernel_vector} is given by
\begin{equation}
  \mathbf{M}^{[n]}_{({\rm P})} = \sum_{{\rm Q}\in \mathcal{B}^{\rm 2D}} \omega_{\rm \an{P-Q}} \, \boldsymbol{\Xi}^{[n]}_{\rm \an{P-Q}} \, \boldsymbol{\Xi}^{[n]^{\intercal}}_{\rm \an{P-Q}} \, A_{\rm (Q)} ,
\end{equation}
and $\tilde{\boldsymbol{\Xi}}^{[n]}$ is the following $5\times d_{\Xi}$ matrix, where $d_{\Xi}$ is the dimension of $\boldsymbol{\Xi}^{[n]}$:
\begin{equation}
  \tilde{\boldsymbol{\Xi}}^{[n]} = 
  \begin{bmatrix} 
    1, \ 0, \ 0, \ 0, \ 0 , \ 0, \ \dots \ , \ 0 \\
    0, \ 1, \ 0, \ 0, \ 0 , \ 0, \ \dots \ , \ 0 \\
    0, \ 0, \ 2, \ 0, \ 0 , \ 0, \ \dots \ , \ 0 \\
    0, \ 0, \ 0, \ 1, \ 0 , \ 0, \ \dots \ , \ 0 \\
    0, \ 0, \ 0, \ 0, \ 2 , \ 0, \ \dots \ , \ 0 
  \end{bmatrix}.
\end{equation}
Note that $n\geq2$ is required in order to have well-defined second-derivative kernel functions.

\begin{example}
  \normalfont
  For n=3 (cubic kernel functions) we have
  \begin{equation}
    \boldsymbol{\Xi}^{[3]} = [\xi_1, \, \xi_2, \, \xi_1^2, \, \xi_1\xi_2, \, \xi_2^2, \, \xi_1^3, \, \xi_1^2\xi_2, \, \xi_1\xi_2^2, \, \xi_2^3]^\intercal ,
  \end{equation}
    \begin{equation}
      \tilde{\boldsymbol{\Xi}}^{[3]} = 
      \begin{bmatrix} 
        1, \ 0, \ 0, \ 0, \ 0 , \ 0, \ 0 , \ 0 , \ 0 \\
        0, \ 1, \ 0, \ 0, \ 0 , \ 0, \ 0 , \ 0 , \ 0 \\
        0, \ 0, \ 2, \ 0, \ 0 , \ 0, \ 0 , \ 0 , \ 0 \\
        0, \ 0, \ 0, \ 1, \ 0 , \ 0, \ 0 , \ 0 , \ 0 \\
        0, \ 0, \ 0, \ 0, \ 2 , \ 0, \ 0 , \ 0 , \ 0 
      \end{bmatrix},
    \end{equation}
and the vector $\boldsymbol{\Phi}_{\rm \an{P-Q}}^{[3]}$ is computed using the expression given by \cref{eqn:kernel_vector}.
\end{example}

\subsection{Computation of the PD Force State}
\label{sec:PDforce}

The methodology described in \cref{sec:manifold} is utilized to construct a local parametric domain for each material point, while the methodology to evaluate parametric gradients is developed in \cref{sec:PD_gradient}. Equipped with these ingredients, we are now able to compute the PD force state.

We compute the parametric derivatives of the current-configuration position vector as
\begin{equation} 
\begin{aligned}
x^{\rm 2D}_{i{,\xi_j}_{({\rm P})}} &= \int_{\mathcal{B}^{\rm 2D}} {x}^{\rm 2D}_{i_{\rm \an{P-Q}}} \, {\Phi}_{{\xi_j}_{\rm \an{P-Q}}} \, {\rm dQ}, \qquad j=1,2, \\
x^{\rm 2D}_{i{,\xi_j\xi_k}_{({\rm P})}} &= \int_{\mathcal{B}^{\rm 2D}} {x}^{\rm 2D}_{i_{\rm \an{P-Q}}} \, {\Phi}_{{\xi_j\xi_k}_{\rm \an{P-Q}}} \, {\rm dQ}, \qquad j=1,2,
\end{aligned}
\label{eqn:dx2Ddxi}
\end{equation}
where
\begin{equation} 
{x}^{\rm 2D}_{i_{\rm \an{P-Q}}} = x^{\rm 2D}_{i_{\rm (Q)}} - x^{\rm 2D}_{i_{\rm (P)}} .
\label{eqn:defx2Dpq}
\end{equation}
Although integral notation is employed in the above expressions, as well as in what follows, we carry out the integrals using nodal quadrature (see, e.g., \cref{eqn:grad1disc,eqn:grad2}). The normal vector $\mathbf{n}$ and the matrices $\mathbf{A}$, $\mathbf{B}^1$ and $\mathbf{B}^2$ are evaluated using \cref{eqn:normal,eqn:A,eqn:B12}. Using these objects, together with \cref{eqn:defF,eqn:x3dd,eqn:dx2Ddxi}, we compute the Jacobian $\mathbf{F}^{\rm 3D}$.

The velocity gradient is a key object in our shell formulation, and we develop it in what follows. The mid-surface parametric gradient is given by
\begin{equation} 
\begin{aligned}
v^{\rm 2D}_{i{,\xi_j}_{({\rm P})}} &= \int_{\mathcal{B}^{\rm 2D}} {v}^{\rm 2D}_{i_{\rm \an{P-Q}}} \, {\Phi}_{{\xi_j}_{\rm \an{P-Q}}} \, {\rm dQ}, \qquad j=1,2, \\
v^{\rm 2D}_{i{,\xi_j\xi_k}_{({\rm P})}} &= \int_{\mathcal{B}^{\rm 2D}} {v}^{\rm 2D}_{i_{\rm \an{P-Q}}} \, {\Phi}_{{\xi_j\xi_k}_{\rm \an{P-Q}}} \, {\rm dQ}, \qquad j=1,2,
\end{aligned}
\label{eqn:dv2Ddxi}
\end{equation}
where
\begin{equation} 
{v}^{\rm 2D}_{i_{\rm \an{P-Q}}} = v^{\rm 2D}_{i_{\rm (Q)}} - v^{\rm 2D}_{i_{\rm (P)}} .
\label{eqn:defv2Dpq}
\end{equation}
Using \cref{eqn:v3dd,eqn:dv2Ddxi}, the full parametric gradient may be expressed as
\begin{equation}
\begin{aligned}
v^{\rm 3D}_{{i,\xi_j}_{({\rm P})}} &= \int_{\mathcal{B}^{\rm 2D}} {\Phi}_{{\xi_j}_{\rm \an{P-Q}}} \, {v}^{\rm 2D}_{i_{\rm \an{P-Q}}} \,  {\rm dQ} 
+ \frac{h_{(\rm P)}}{2} \xi_3 \int_{\mathcal{B}^{\rm 2D}} \left( B^1_{{ik}_{({\rm P})}} \, {\Phi}_{{\xi_1\xi_j}_{\rm \an{P-Q}}} + B^2_{{ik}_{({\rm P})}} \, {\Phi}_{{\xi_2\xi_j}_{\rm \an{P-Q}}} \right) \, {v}^{\rm 2D}_{k_{\rm \an{P-Q}}} \, {\rm dQ} \\
& \quad + \frac{h_{(\rm P)}}{2} \xi_3 \int_{\mathcal{B}^{\rm 2D}} \left( B^1_{{ik,\xi_j}_{({\rm P})}} \, {\Phi}_{{\xi_1}_{\rm \an{P-Q}}} + B^2_{{ik,\xi_j}_{({\rm P})}} \, {\Phi}_{{\xi_2}_{\rm \an{P-Q}}} \right) \, {v}^{\rm 2D}_{k_{\rm \an{P-Q}}} \, {\rm dQ} , \quad j=1,2 \  (\text{in-plane}), \\
v^{\rm 3D}_{{i,\xi_3}_{({\rm P})}} &= \frac{h_{(\rm P)}}{2} \int_{\mathcal{B}^{\rm 2D}} \left( B^1_{{ik}_{({\rm P})}} \, {\Phi}_{{\xi_1}_{\rm \an{P-Q}}} + B^2_{{ik}_{({\rm P})}} \, {\Phi}_{{\xi_2}_{\rm \an{P-Q}}} \right) \, {v}^{\rm 2D}_{k_{\rm \an{P-Q}}} \, {\rm dQ} , \quad (\text{through-thickness}).
\end{aligned}
\label{eqn:dv3Ddxi}
\end{equation}
The spatial velocity gradient $v^{\rm 3D}_{{i,j}_{({\rm P})}}$ can be obtained using \cref{eqn:transferV,eqn:dv3Ddxi}, 
\begin{equation}
v^{\rm 3D}_{{i,k}_{({\rm P})}} = \left( \int_{\mathcal{B}^{\rm 2D}} \beta^{\rm 3D}_{{ilj}_{\rm \an{P-Q}}} \, {v}^{\rm 2D}_{l_{\rm \an{P-Q}}} \, {\rm dQ} \right) F^{\rm 3D^{-1}}_{{jk}_{({\rm P})}}, 
\label{eqn:dv3Ddx}
\end{equation}
where the auxiliary object $\boldsymbol{\beta}^{\rm 3D}_{\rm \an{P-Q}}$ is given by
\begin{equation}
\begin{aligned}
&\beta^{\rm 3D}_{{ilj}_{\rm \an{P-Q}}} = {\Phi}_{{\xi_j}_{\rm \an{P-Q}}} \, \delta_{il} + \frac{h_{(\rm P)}}{2} \xi_3 \left( B^1_{{il}_{({\rm P})}} \, {\Phi}_{{\xi_1\xi_j}_{\rm \an{P-Q}}} + B^2_{{il}_{({\rm P})}} \, {\Phi}_{{\xi_2\xi_j}_{\rm \an{P-Q}}} + B^1_{{il,\xi_j}_{({\rm P})}} \, {\Phi}_{{\xi_1}_{\rm \an{P-Q}}} + B^2_{{il,\xi_j}_{({\rm P})}} \, {\Phi}_{{\xi_2}_{\rm \an{P-Q}}} \right) , \quad j=1,2 \\
&\beta^{\rm 3D}_{{ilj}_{\rm \an{P-Q}}} = \frac{h_{(\rm P)}}{2} \left( B^1_{{il}_{({\rm P})}} \, {\Phi}_{{\xi_1}_{\rm \an{P-Q}}} + B^2_{{il}_{({\rm P})}} \, {\Phi}_{{\xi_2}_{\rm \an{P-Q}}} \right) . \quad j=3.
\end{aligned}
\label{eqn:beta}
\end{equation}

To avoid instability, we use the bond-stabilization technique~\cite{behzadinasab2020semi} and define the following bond-associated velocity gradient directly in the physical domain:
\begin{equation} 
{v}^{\rm 3D}_{i,j_{\rm \an{P-Q}}} = \frac{v^{\rm 3D}_{{i,j}_{({\rm P})}} + v^{\rm 3D}_{{i,j}_{({\rm Q})}}}{2} + \left( {v}^{\rm 3D}_{i_{\rm \an{P-Q}}} - \frac{v^{\rm 3D}_{{i,k}_{({\rm P})}} + v^{\rm 3D}_{{i,k}_{({\rm Q})}}}{2}  \, {x}^{\rm 3D}_{k_{\rm \an{P-Q}}} \right) \frac{{x}^{\rm 3D}_{j_{\rm \an{P-Q}}}}{|{x}^{\rm 3D}_{\rm \an{P-Q}}|^2} ,
\label{eqn:bondL}
\end{equation}
where
\begin{equation} 
{v}^{\rm 3D}_{i_{\rm \an{P-Q}}} = v^{\rm 3D}_{i_{\rm (Q)}} - v^{\rm 3D}_{i_{\rm (P)}} ,
\label{eqn:defv3Dpq}
\end{equation}
and
\begin{equation} 
{x}^{\rm 3D}_{i_{\rm \an{P-Q}}} = x^{\rm 3D}_{i_{\rm (Q)}} - x^{\rm 3D}_{i_{\rm (P)}} .
\label{eqn:defx3Dpq}
\end{equation}
Introducing the Kirchhoff-Love kinematics into \cref{eqn:defv3Dpq,eqn:defx3Dpq} allows us to re-express the full position and velocity vectors in terms of their mid-surface counterparts as
\begin{equation} 
\begin{aligned}
{x}^{\rm 3D}_{i_{\rm \an{P-Q}}} 
&= \left( x^{2D}_{i_{\rm (Q)}} + \xi_3 \frac{h_{\rm (Q)}}{2} n_{{i}_{({\rm Q})}} \right) - \left( x^{2D}_{i_{\rm (P)}} - \xi_3 \frac{h_{\rm (P)}}{2} n_{{i}_{({\rm P})}} \right) \\
&= {x}^{\rm 2D}_{i_{\rm \an{P-Q}}} + \frac{\xi_3}{2} \left( h_{\rm (Q)} \, n_{{i}_{({\rm Q})}} - h_{\rm (P)} \, n_{{i}_{({\rm P})}} \right)
\end{aligned}
\label{eqn:x3Dpq}
\end{equation}
and
\begin{equation} 
\begin{aligned}
{v}^{\rm 3D}_{i_{\rm \an{P-Q}}} 
&= \left( v^{2D}_{i_{\rm (Q)}} + \xi_3 \frac{h_{\rm (Q)}}{2} \dot{n}_{{i}_{({\rm Q})}} \right) - \left( v^{2D}_{i_{\rm (P)}} - \xi_3 \frac{h_{\rm (P)}}{2} \dot{n}_{{i}_{({\rm P})}} \right) \\
&= {v}^{\rm 2D}_{i_{\rm \an{P-Q}}} 
+ \left( \int_{\mathcal{B}^{\rm 2D}} \gamma^{\rm 3D}_{{ik}_{\rm \an{Q-S}}} \, {v}^{\rm 2D}_{k_{\rm \an{Q-S}}} \, {\rm dS} - \int_{\mathcal{B}^{\rm 2D}} \gamma^{\rm 3D}_{{ik}_{\rm \an{P-S}}} \, {v}^{\rm 2D}_{k_{\rm \an{P-S}}} \, {\rm dS} \right) ,
\end{aligned}
\label{eqn:bondv3D}
\end{equation}
where the auxiliary object $\boldsymbol{\gamma}^{\rm 3D}_{\rm \an{P-S}}$ is given by
\begin{equation} 
\gamma^{\rm 3D}_{{ik}_{\rm \an{P-S}}} = \frac{h_{\rm (P)}}{2} \, \xi_3 \left( B^1_{{ik}_{({\rm P})}} \, {\Phi}_{{\xi_1}_{\rm \an{P-S}}} + B^2_{{ik}_{({\rm P})}} \, {\Phi}_{{\xi_2}_{\rm \an{P-S}}} \right) .
\label{eqn:gamma}
\end{equation}
The bond-associated Kirchhoff stress is given by 
\begin{equation} 
\boldsymbol{\tau}^{\rm 3D}_{\rm \an{P-Q}} = {J}^{\rm 3D}_{\rm \an{P-Q}} \, \boldsymbol{\sigma}^{\rm 3D}_{\rm \an{P-Q}} ,
\label{eqn:tau-sigma}
\end{equation}
where ${J}^{\rm 3D}_{\rm \an{P-Q}}$ is the bond-associated Jacobian of the deformation gradient evolved using the trace of the bond-associated spatial velocity gradient
\begin{equation} 
{\dot{J}}^{\rm 3D}_{\rm \an{P-Q}} = {J}^{\rm 3D}_{\rm \an{P-Q}} \, {v}^{\rm 3D}_{i,i_{\rm \an{P-Q}}} ,
\label{eqn:Jdot}
\end{equation}
and $\boldsymbol{\sigma}^{\rm 3D}_{\rm \an{P-Q}}$ is the bond-associated Cauchy stress that is obtained from the appropriate rate-form constitutive law. The latter leads to a stress update algorithm driven by the stain increment that depends on the bond-associated symmetric spatial velocity gradient and the time step size.

Introducing the objects derived in this section into \cref{eqn:U_stress_based} for the stress power yields the following expressions: 
\begin{equation}   
\begin{aligned}
&\dot{\mathcal{U}}_{(\mathcal{B})} = \dot{\rm U}1 + \dot{\rm U}2 + \dot{\rm U}3 , \\
&\dot{\rm U}1 = \int_{\mathcal{B}^{\rm 2D}} \int_{\mathcal{B}^{\rm 2D}} \int_{-1}^{1} {a}^{\rm 3D}_{i_{\rm \an{P-Q}}} \, {v}^{\rm 3D}_{i_{\rm \an{P-Q}}} \, {\rm d}\xi_3 \, {\rm dQ} \, {\rm dP} , \\ 
&\dot{\rm U}2 = \int_{\mathcal{B}^{\rm 2D}} \int_{\mathcal{B}^{\rm 2D}} \int_{-1}^{1} {b}^{\rm 3D}_{ik_{\rm \an{P-Q}}} \, v^{\rm 3D}_{{i,k}_{({\rm P})}} \, {\rm d}\xi_3 \, {\rm dQ} \, {\rm dP} , \\
&\dot{\rm U}3 = \int_{\mathcal{B}^{\rm 2D}} \int_{\mathcal{B}^{\rm 2D}} \int_{-1}^{1} {b}^{\rm 3D}_{ik_{\rm \an{P-Q}}} \,  v^{\rm 3D}_{{i,k}_{({\rm Q})}} \, {\rm d}\xi_3 \, {\rm dQ} \, {\rm dP} ,
\end{aligned}
\label{eqn:U1}
\end{equation}
where the auxiliary objects $\mathbf{a}^{\rm 3D}_{\rm \an{P-Q}}$ and $\mathbf{b}^{\rm 3D}_{\rm \an{P-Q}}$ are given by
\begin{equation} 
\begin{aligned}  
&{a}^{\rm 3D}_{i_{\rm \an{P-Q}}} = \frac{h_{\rm (P)}}{2} \, \alpha_{\rm \an{P-Q}} \, {\tau}^{\rm 3D}_{ij_{\rm \an{P-Q}}} \frac{{x}^{\rm 3D}_{j_{\rm \an{P-Q}}}}{|{x}^{\rm 3D}_{\rm \an{P-Q}}|^2} , \\
&{b}^{\rm 3D}_{ik_{\rm \an{P-Q}}} = \frac{h_{\rm (P)}}{4} \, \alpha_{\rm \an{P-Q}} \, {\tau}^{\rm 3D}_{ij_{\rm \an{P-Q}}} \Bigg( \delta_{jk} - \frac{{x}^{\rm 3D}_{j_{\rm \an{P-Q}}} \, {x}^{\rm 3D}_{k_{\rm \an{P-Q}}}}{|{x}^{\rm 3D}_{\rm \an{P-Q}}|^2} \Bigg) .
\end{aligned}  
\label{eqn:ab}
\end{equation}
The first term in \cref{eqn:U1} is evaluated using \cref{eqn:bondv3D} resulting in
\begin{equation} 
\begin{aligned}  
\dot{\rm U}1 &= \int_{\mathcal{B}^{\rm 2D}} \int_{\mathcal{B}^{\rm 2D}} \int_{-1}^{1} {a}^{\rm 3D}_{i_{\rm \an{P-Q}}} \, {v}^{\rm 2D}_{i_{\rm \an{P-Q}}} \, {\rm d}\xi_3 \, {\rm dQ} \, {\rm dP} \\
& \quad + \int_{\mathcal{B}^{\rm 2D}} \int_{\mathcal{B}^{\rm 2D}} \int_{-1}^{1} {a}^{\rm 3D}_{i_{\rm \an{P-Q}}} \int_{\mathcal{B}^{\rm 2D}} \gamma^{\rm 3D}_{{ik}_{\rm \an{Q-S}}} \, {v}^{\rm 2D}_{k_{\rm \an{Q-S}}} \, {\rm dS} \, {\rm d}\xi_3 \, {\rm dQ} \, {\rm dP} \\
& \quad - \int_{\mathcal{B}^{\rm 2D}} \int_{\mathcal{B}^{\rm 2D}} \int_{-1}^{1} {a}^{\rm 3D}_{i_{\rm \an{P-Q}}} \int_{\mathcal{B}^{\rm 2D}} \gamma^{\rm 3D}_{{ik}_{\rm \an{P-S}}} \, {v}^{\rm 2D}_{k_{\rm \an{P-S}}} \, {\rm dS} \, {\rm d}\xi_3 \, {\rm dQ} \, {\rm dP} .
\end{aligned}
\label{eqn:U11}
\end{equation}
Interchanging the order of integration in \cref{eqn:U11} we get
\begin{equation} 
\begin{aligned}  
\dot{\rm U}1 &= \int_{\mathcal{B}^{\rm 2D}} \int_{\mathcal{B}^{\rm 2D}} \left( \int_{-1}^{1} {a}^{\rm 3D}_{i_{\rm \an{P-Q}}} \, {\rm d}\xi_3 \right) {v}^{\rm 2D}_{i_{\rm \an{P-Q}}} \, {\rm dQ} \, {\rm dP} \\
& \quad + \int_{\mathcal{B}^{\rm 2D}} \int_{\mathcal{B}^{\rm 2D}} \left( \int_{\mathcal{B}^{\rm 2D}} \int_{-1}^{1} \gamma^{\rm 3D}_{{ik}_{\rm \an{P-Q}}} \, {a}^{\rm 3D}_{i_{\rm \an{S-P}}} \, {\rm d}\xi_3 \, {\rm dS} \right) {v}^{\rm 2D}_{k_{\rm \an{P-Q}}} \, {\rm dQ} \, {\rm dP} \\
& \quad - \int_{\mathcal{B}^{\rm 2D}} \int_{\mathcal{B}^{\rm 2D}} \left( \int_{\mathcal{B}^{\rm 2D}} \int_{-1}^{1} \gamma^{\rm 3D}_{{ik}_{\rm \an{P-Q}}} \, {a}^{\rm 3D}_{i_{\rm \an{P-S}}} \, {\rm d}\xi_3 \, {\rm dS} \right) {v}^{\rm 2D}_{k_{\rm \an{P-Q}}} \, {\rm dQ} \, {\rm dP} \\
&= \int_{\mathcal{B}^{\rm 2D}} \int_{\mathcal{B}^{\rm 2D}} \left( \int_{-1}^{1} {a}^{\rm 3D}_{i_{\rm \an{P-Q}}} \, {\rm d}\xi_3 \right) {v}^{\rm 2D}_{i_{\rm \an{P-Q}}} \, {\rm dQ} \, {\rm dP} \\
& \quad + \int_{\mathcal{B}^{\rm 2D}} \int_{\mathcal{B}^{\rm 2D}} \left( \int_{\mathcal{B}^{\rm 2D}} \int_{-1}^{1} \gamma^{\rm 3D}_{{ik}_{\rm \an{P-Q}}} \left( {a}^{\rm 3D}_{i_{\rm \an{S-P}}} - {a}^{\rm 3D}_{i_{\rm \an{P-S}}} \right) \, {\rm d}\xi_3 \, {\rm dS} \right) {v}^{\rm 2D}_{k_{\rm \an{P-Q}}} \, {\rm dQ} \, {\rm dP} .
\end{aligned}
\label{eqn:Udot1}
\end{equation}
The second term in \cref{eqn:U1} can be rewritten using \cref{eqn:transferV,eqn:dv3Ddx} as
\begin{equation}   
\begin{aligned}
\dot{\rm U}2 &= \int_{\mathcal{B}^{\rm 2D}} \int_{\mathcal{B}^{\rm 2D}} \int_{-1}^{1} {b}^{\rm 3D}_{ik_{\rm \an{P-Q}}} \,  \left( \int_{\mathcal{B}^{\rm 2D}} \beta^{\rm 3D}_{{ilj}_{\rm \an{P-S}}} \, {v}^{\rm 2D}_{l_{\rm \an{P-S}}} \, {\rm dS} \right) F^{\rm 3D^{-1}}_{{jk}_{({\rm P})}} \, {\rm d}\xi_3 \, {\rm dQ} \, {\rm dP} \\
&= \int_{\mathcal{B}^{\rm 2D}} \int_{\mathcal{B}^{\rm 2D}} F^{\rm 3D^{-1}}_{{jk}_{({\rm P})}} \, \beta^{\rm 3D}_{{ilj}_{\rm \an{P-Q}}}  \left( \int_{\mathcal{B}^{\rm 2D}} \int_{-1}^{1} {b}^{\rm 3D}_{ik_{\rm \an{P-S}}} \, {\rm d}\xi_3 \, {\rm dS} \right) \, {v}^{\rm 2D}_{l_{\rm \an{P-Q}}} \, {\rm dQ} \, {\rm dP} .
\end{aligned}
\label{eqn:Udot2}
\end{equation}
Similarly, the last term in \cref{eqn:U1} may be expanded using \cref{eqn:transferV,eqn:dv3Ddx} as
\begin{equation}   
\begin{aligned}
\dot{\rm U}3 &= \int_{\mathcal{B}^{\rm 2D}} \int_{\mathcal{B}^{\rm 2D}} \int_{-1}^{1} {b}^{\rm 3D}_{ik_{\rm \an{P-Q}}} \,  \left( \int_{\mathcal{B}^{\rm 2D}} \beta^{\rm 3D}_{{ilj}_{\rm \an{Q-S}}} \, {v}^{\rm 2D}_{l_{\rm \an{Q-S}}} \, {\rm dS} \right) F^{\rm 3D^{-1}}_{{jk}_{({\rm Q})}} \, {\rm d}\xi_3 \, {\rm dQ} \, {\rm dP} \\
&= \int_{\mathcal{B}^{\rm 2D}} \int_{\mathcal{B}^{\rm 2D}} F^{\rm 3D^{-1}}_{{jk}_{({\rm P})}} \, \beta^{\rm 3D}_{{ilj}_{\rm \an{P-Q}}} \left( \int_{\mathcal{B}^{\rm 2D}} \int_{-1}^{1} {b}^{\rm 3D}_{ik_{\rm \an{S-P}}} \, {\rm d}\xi_3 \, {\rm dS} \right) \, {v}^{\rm 2D}_{l_{\rm \an{P-Q}}} \, {\rm dQ} \, {\rm dP} .
\end{aligned}
\label{eqn:Udot3}
\end{equation}

Finally, the rate of rotational kinetic energy may be expressed as
\begin{equation} 
\begin{aligned}
\int_{\mathcal{B}^{\rm 2D}} \frac{\rho^{\rm 2D}_{(\rm P)} \, h^2_{\rm (P)}}{12} \, \ddot{\mathbf{n}}_{\rm (P)} \cdot \dot{\mathbf{n}}_{\rm (P)} \, {\rm dP} &= \int_{\mathcal{B}^{\rm 2D}} \frac{\rho^{\rm 2D}_{(\rm P)} \, h^2_{\rm (P)}}{12} \, \ddot{n}_{i_{\rm (P)}} \left( B^1_{{ik}_{({\rm P})}} \, v^{\rm 2D}_{{k,\xi_1}_{({\rm P})}} + B^2_{{ik}_{({\rm P})}} \, v^{\rm 2D}_{{k,\xi_2}_{({\rm P})}} \right) \, {\rm dP} \\
&= \int_{\mathcal{B}^{\rm 2D}} \int_{\mathcal{B}^{\rm 2D}} \frac{\rho^{\rm 2D}_{(\rm P)} \, h^2_{\rm (P)}}{12} \, \ddot{n}_{i_{\rm (P)}} \left( B^1_{{ik}_{({\rm P})}} \, {\Phi}_{{\xi_1}_{\rm \an{P-Q}}} + B^2_{{ik}_{({\rm P})}} \, {\Phi}_{{\xi_2}_{\rm \an{P-Q}}} \right) \, {v}^{\rm 2D}_{k_{\rm \an{P-Q}}} \, {\rm dQ} \, {\rm dP} ,
\end{aligned}
\label{eqn:rot_energy}
\end{equation}
where $\ddot{n}_{i_{\rm (P)}}$ can be obtained analytically by taking a time derivative of \cref{eqn:ndot} or numerically using a finite difference approach.

Combining the above results and comparing with \cref{eqn:U_T_based_2,eqn:T_rot}, we deduce that the total force state that includes the internal stress and rotation contributions becomes
\begin{equation}   
\begin{aligned}
{T}_{i_{\rm \an{P-Q}}} 
& = \int_{-1}^{1} {a}^{\rm 3D}_{i_{\rm \an{P-Q}}} \, {\rm d}\xi_3 + \int_{-1}^{1} \left( \int_{\mathcal{B}^{\rm 2D}} \left( {a}^{\rm 3D}_{k_{\rm \an{S-P}}} - {a}^{\rm 3D}_{k_{\rm \an{P-S}}} \right) \, {\rm dS}\right) \gamma^{\rm 3D}_{{ki}_{\rm \an{P-Q}}} \, {\rm d}\xi_3  \\
& + \int_{-1}^{1} \left( \int_{\mathcal{B}^{\rm 2D}} \left( {b}^{\rm 3D}_{lk_{\rm \an{P-S}}} + {b}^{\rm 3D}_{lk_{\rm \an{S-P}}} \right) \, {\rm dS} \right) \left(F^{\rm 3D}_{{jk}_{({\rm P})}}\right)^{-1} \, \beta^{\rm 3D}_{{lij}_{\rm \an{P-Q}}} \, {\rm d}\xi_3 \\
& + \frac{\rho^{\rm 2D}_{(\rm P)} \, h^2_{\rm (P)}}{12} \, \ddot{n}_{k_{\rm (P)}} \left( B^1_{{ki}_{({\rm P})}} \, {\Phi}^1_{{\xi_1}_{\rm \an{P-Q}}} + B^2_{{ki}_{({\rm P})}} \, {\Phi}^1_{{\xi_2}_{\rm \an{P-Q}}} \right) .
\end{aligned}
\label{eqn:force}
\end{equation}

The through-thickness integration is carried out using Gaussian quadrature in the parametric $\xi_3$ direction. Let $N_G$ denote the number of through-thickness quadrature points and let $w_I$ denote the $I^{\rm th}$ quadrature weight. The space-discrete PD formulation of the KL shell can be thus summarized as follows:
\begin{equation} 
\rho^{\rm 2D}_{(\rm P)} \, \dot{\mathbf{v}}^{\rm 2D}_{(\rm P)} = \sum_{{\rm Q}\in\mathcal{B}^{\rm 2D}} \left( \mathbf{T}_{\rm \an{P-Q}} - \mathbf{T}_{\rm \an{Q-P}} \right) \, A_{\rm (Q)} + \rho^{\rm 2D}_{(\rm P)}~\mathbf{b}_{(\rm P)},
\end{equation}
where
\begin{equation}   
\begin{aligned}
& {T}_{i_{\rm \an{P-Q}}} = \sum_{I=1}^{N_G} w^I \, {T}_{i_{\rm \an{P-Q}}}^{{\rm int}~I} + {T}_{i_{\rm \an{P-Q}}}^{\rm rot} , \\
& {T}_{i_{\rm \an{P-Q}}}^{{\rm int}~I} = {a}^{{I}}_{i_{\rm \an{P-Q}}} + \bar{a}^{{I}}_{k_{\rm (P)}} \ \gamma^{{I}}_{{ki}_{\rm \an{P-Q}}} + \bar{b}^{{I}}_{lk_{\rm (P)}} \ \beta^{{I}}_{{lij}_{\rm \an{P-Q}}} \left(F^{{I}}_{{jk}_{({\rm P})}}\right)^{-1} , \\
& \bar{a}^{{I}}_{k_{\rm (P)}} = \sum_{{\rm S}\in\mathcal{B}^{\rm 2D}} \left( {a}^{{I}}_{k_{\rm \an{S-P}}} - {a}^{{I}}_{k_{\rm \an{P-S}}} \right) A_{\rm (S)} , \\
& \bar{b}^{{I}}_{lk_{\rm (P)}} = \sum_{{\rm S}\in\mathcal{B}^{\rm 2D}} \left( {b}^{{I}}_{lk_{\rm \an{P-S}}} + {b}^{{I}}_{lk_{\rm \an{S-P}}} \right) A_{\rm (S)} , \\
& {T}_{i_{\rm \an{P-Q}}}^{\rm rot} = \frac{\rho^{\rm 2D}_{(\rm P)} \, h^2_{\rm (P)}}{12} \, \ddot{n}_{k_{\rm (P)}} \left( B^1_{{ki}_{({\rm P})}} \, {\Phi}_{{\xi_1}_{\rm \an{P-Q}}} + B^2_{{ki}_{({\rm P})}} \, {\Phi}_{{\xi_2}_{\rm \an{P-Q}}} \right) ,
\end{aligned}
\label{eqn:force-discrete}
\end{equation}
in which $\mathbf{T}^{{\rm int}~I}$ is the internal force state at the through-thickness quadrature point $I$, and $\mathbf{T}^{\rm rot}$ is the rotational force state, which is pre-integrated through the thickness analytically. The algorithmic implementation of the above formulation is summarized in \cref{sec:implementation}. A linearized version of the formulation is provided in \cref{sec:linearization}, and may be used in quasi-static or implicit calculations.

\subsection{Stress Update}
\label{sec:stress_update}

It is necessary to ensure that the principle of material frame indifference, or objectivity, is maintained in the rate-form constitutive relations~\cite{belytschko2013nonlinear}. To this end we make use of the co-rotational formulation for the stress update. The rate of deformation $\mathbf{D}$ is defined as the symmetric part of the velocity gradient
\begin{equation}
{D}_{ij} = \frac{1}{2} \left( {v}^{3D}_{i,j} + {v}^{3D}_{j,i} \right) .
\label{eqn:rate_of_def}
\end{equation}
The rate of deformation $\hat{\mathbf{D}}$ in the unrotated coordinate configuration is related to $\mathbf{D}$ by the following transformation:
\begin{equation}
\hat{\mathbf{D}} = \mathbf{R}^{\intercal} \, \mathbf{D} \, \mathbf{R} ,
\label{eqn:unrot_rate_of_def}
\end{equation}
where $\mathbf{R}$ is the rotation tensor associated with the polar decomposition $\mathbf{F} = \mathbf{V} \, \mathbf{R}$. The numerical algorithm of~\cite{flanagan1987accurate}, summarized in \cref{alg:rotation_alg}, is used to evolve $\mathbf{R}$ in time. In the polar decomposition, $\mathbf{V}$ is the left-stretch tensor, which is equal to the identity tensor at every point in the undeformed configuration. Note that the co-rotational approach does not suffer from the spurious stress oscillations associated with the Jaumann rate in shear-dominated deformations~\cite{dienes1979analysis,flanagan1987accurate}. 

\begin{algorithm}
	\begin{algorithmic}[1]
    \State Decompose the bond-level velocity gradient into symmetric ($\mathbf{D}$) and anti-symmetric ($\mathbf{W}$) parts.
    
	\State Calculate auxiliary vector states $\mathbf{z}$, $\boldsymbol{\omega}$, and the rate of rigid-body rotation tensor $\boldsymbol{\Omega}$ such that,
    \begin{sloppypar}
    ${z}_i =\epsilon_{ikj} \, {D}_{jm} \, {V}_{mk}$
    \end{sloppypar}
    \begin{sloppypar}
    $\boldsymbol{\omega} = \mathbf{w} +  \left[ \, \mathbf{I} \, {\rm tr} (\mathbf{V}^n) - \mathbf{V}^n \right]^{-1} \mathbf{z}$ ~ where ~ ${W}_{ij} = \epsilon_{ikj} \, {w}_k$
    \end{sloppypar}
    \begin{sloppypar}
    ${\Omega}_{ij} = \epsilon_{ikj} \, {\omega}_{k}$
    \end{sloppypar}
    
	\State Solve for the bond-level rotation matrix,
    \begin{sloppypar}
     $\left(\mathbf{I} - \dfrac{1}{2} \Delta t \, \boldsymbol{\Omega}\right) \mathbf{R}^{(n+1)} =   \left(\mathbf{I}  + \dfrac{1}{2} \Delta t \, \boldsymbol{\Omega}\right) \mathbf{R}^{(n)} $
    \end{sloppypar}
    
    \State Compute the rate of change of the bond-level left-stretch tensor $\dot{\mathbf{V}}$,
    \begin{sloppypar}
    $\dot{\mathbf{V}} = \left(\mathbf{D} + \mathbf{W}\right) \mathbf{V}^{(n)} - \mathbf{V}^{(n)} \, \boldsymbol{\Omega}$
    \end{sloppypar}
    
    \State Update the left-stretch tensor state $\mathbf{V}^{(n+1)}$,
    \begin{sloppypar}
    $\mathbf{V}^{(n+1)} = \mathbf{V}^{(n)} + \dot{\mathbf{V}} \, \Delta t  $
    \end{sloppypar}
    
    \State Store the left-stretch and rotation tensor states for the next time step,
    \begin{sloppypar}
    $\mathbf{V}^{(n)} \leftarrow \mathbf{V}^{(n+1)}$
    \end{sloppypar}
    \begin{sloppypar}
    $\mathbf{R}^{(n)} \leftarrow \mathbf{R}^{(n+1)}$
    \end{sloppypar}
	\end{algorithmic} 
	\caption{Time integration of the rotation matrix~\cite{flanagan1987accurate} at the bond level} 
	\label{alg:rotation_alg}
\end{algorithm}

The stress update is performed in the unrotated configuration as
\begin{equation}
\hat{\boldsymbol{\sigma}}^{(n+1)} = \hat{\boldsymbol{\sigma}}^{(n)} + \Delta \hat{\boldsymbol{\sigma}} ,
\label{eqn:unrot_stress}
\end{equation}
using a constitutive law for the stress rate, i.e., 
\begin{equation}
\Delta \hat{\boldsymbol{\sigma}} = \dot{\boldsymbol\sigma}(\hat{\mathbf{D}}) \, \Delta t .
\label{eqn:stress_inc}
\end{equation}
The updated stress tensor is rotated back to the deformed configuration as
\begin{equation}
\boldsymbol{\sigma}^{(n+1)} = \mathbf{R}^{(n+1)} \, \hat{\boldsymbol{\sigma}}^{(n+1)} \, \mathbf{R}^{{(n+1)}^{\intercal}} .
\label{eqn:rotated_stress}
\end{equation}
Note that all the above calculations are performed at the {\bf bond level}. 

\begin{remark}

In the computation, at time $t=0$ the rotation tensor is initialized as
\begin{equation} 
\begin{aligned}
\mathbf{R}^{(0)} & = 
\begin{bmatrix}
\uparrow & \uparrow & \uparrow \\
\hat{\mathbf{e}}^1 & \hat{\mathbf{e}}^2 & \hat{\mathbf{n}} \\
\uparrow & \uparrow & \uparrow
\end{bmatrix} 
,
\end{aligned}
\label{eqn:rot_tensor}
\end{equation}
where $\hat{\mathbf{n}}$ is the mid-surface unit normal in the undeformed configuration, and $\hat{\mathbf{e}}^1$ and $\hat{\mathbf{e}}^2$ are the orthonormal vectors in the shell mid-surface tangent plane. As the rotation matrix is evolved in time according to \cref{alg:rotation_alg}, its first two columns will contain the vectors that align with the material principal axes, while its third column will contain the surface normal vector, all in the current configuration. This setup is particularly convenient for: i. Modeling of anisotropic materials, such as fiber-reinforced composites, where the unrotated-configuration stress components $\hat{\sigma}_{11}$, $\hat{\sigma}_{22}$, and $\hat{\sigma}_{12}=\hat{\sigma}_{21}$ correspond to the fiber, matrix, and in-plane shear stress, respectively, and may be used directly to drive the corresponding damage modes (see, e.g.,~\cite{hashin1980failure,matzenmiller1995constitutive,lapczyk2007progressive}); ii. Enforcing the zero through-thickness stress condition directly on the $\hat{\sigma}_{33}$ component of the unrotated stress (see \cref{sec:zero_stress}). In the present formulation, the rotation tensor is initialized at each bond ${\rm \an{P-Q}}$, for which we define the bond-associated unit normal in the undeformed configuration
\begin{equation} 
\hat{\mathbf{n}}_{\rm \an{P-Q}} = \frac{ \hat{\mathbf{n}}_{\rm (P)} + \hat{\mathbf{n}}_{\rm (Q)} }{ \left| \hat{\mathbf{n}}_{\rm (P)} + \hat{\mathbf{n}}_{\rm (Q)} \right| }.
\label{eqn:undef_normal}
\end{equation}

\end{remark}

\subsection{$J_2$ Plasticity Theory}
\label{sec:plasticity}

In this work, we use the rate-independent three-dimensional von Mises $J_2$ (isotropic) plasticity~\cite{de2011computational} with a return mapping algorithm, which is briefly recalled here. Note that all the calculations are performed in the co-rotational system detailed in the previous section. 

The rate of deformation tensor is decomposed additively into elastic and plastic parts as
\begin{equation}\label{eq:}
\hat{\mathbf{D}} = \hat{\mathbf{D}}^{el}+ \hat{\mathbf{D}}^{pl}.
\end{equation}
An elastic predictor step is assumed for computing the trial Cauchy stress state:
\begin{equation}\label{eq:trial_stress}
 \hat{{\sigma}}^{tr}_{ij} = \hat{{\sigma}}_{ij}^{(n)} + \Delta t \, \hat{\mathbb{C}}^{el}_{ijkl} \, \hat{{D}}_{kl}
\end{equation}
where $\hat{\mathbb{C}}^{el}_{ijkl} = \lambda \, \delta_{ij} \, \delta_{kl} + \mu \left(\delta_{ik} \, \delta_{jl} + \delta_{il} \, \delta_{jk}\right)$ is the standard fourth-order isotropic elasticity tensor with Lam\'e parameters $\lambda$ and $\mu$. Next, we calculate the deviatoric part of the trial Cauchy stress as
\begin{equation}
\label{eq:trial_dev_stress}
\hat{{\sigma}}_{ij}^{dev} = \hat{{\sigma}}^{tr}_{ij} - \frac{1}{3} \, \hat{{\sigma}}^{tr}_{kk} \, \delta_{ij}
\end{equation}
and the trial von Mises stress as
\begin{equation}\label{eq:trial_vm_stress}
    \hat{{\sigma}}^{vm} = \sqrt{ \frac{3}{2} \, \hat{{\sigma}}_{ij}^{dev} \, \hat{{\sigma}}_{ij}^{dev}}.
\end{equation}
The yield function $f$ is defined as
\begin{equation}
\label{eq:yield_surface}
    f\left(\hat{{\sigma}}^{vm}, {\bar{\epsilon}}^P\right) = \hat{{\sigma}}^{vm} - \sigma_{Y}\left({\bar{\epsilon}}^P\right),
\end{equation}
where the yield stress $\sigma_Y$ is a function of the equivalent plastic strain ${\bar{\epsilon}}^P$.

If $f\left(\hat{{\sigma}}^{vm}, {\bar{\epsilon}}^P\right) \leq 0$, the stress update is completed with the elastic predictor. Otherwise, plastic yielding occurs and we radially return to the yield surface such that $f\left(\hat{{\sigma}}^{vm}, {\bar{\epsilon}}^P\right) = 0$. Using Newton's method, the plastic consistency parameter $\Delta \gamma$ is computed as follows:
\begin{equation}
\begin{aligned}
\Delta\gamma^{(k+1)}  &=  \Delta\gamma^{(k)} + \frac{\bar{f}\left(\Delta\gamma^{(k)}\right)}{3  \, \mu  + H^{(k)}}, \\
\bar{f}\left(\Delta\gamma^{(k)}\right) &= \hat{{\sigma}}^{vm} - 3  \, \mu  \,  \Delta\gamma^{(k)} - \sigma_{Y}\left({\bar{\epsilon}}^P_{(n)} + \Delta\gamma^{(k)} \right), 
\end{aligned}
\label{eq:consistency}
\end{equation}
where $H^{(k)} = \dfrac{d \sigma_Y}{d {\bar{\epsilon}}^P}\Big|_{{\bar{\epsilon}}^P_{(n)} + \Delta\gamma^{(k)}}$ is the hardening modulus computed at the iteration step $k$. Once $\bar{f}\left(\Delta\gamma^{(k+1)}\right)=0$, the yield condition is met, and the step is completed by updating the stress tensor and equivalent plastic strain states using
\begin{equation}
\begin{aligned}
\hat{{\sigma}}_{ij}^{(n+1)} &= \frac{1}{3} \, \hat{{\sigma}}^{tr}_{kk} \, \delta_{ij} + \left(1-\frac{3 \, \mu \, \Delta\gamma^{(k+1)}}{\hat{\sigma}^{vm}}\right)\hat{{\sigma}}_{ij}^{dev} , \\
{\bar{\epsilon}}^P_{(n+1)} &= {\bar{\epsilon}}^P_{(n)} + \Delta\gamma^{(k+1)} .
\end{aligned}
\label{eq:plastic_stress_update}
\end{equation}

\subsection{Enforcing Zero Through-Thickness Stress and Updating the Jacobian and Thickness}
\label{sec:zero_stress}

The direct use of 3D constitutive modeling in conjunction with the KL shell kinematics may result in a non-zero through-thickness stress state and lead to an overly stiff structural response~\cite{bischoff2018models}. Following the approach in~\cite{hallquist2006ls,bazilevs2009computational,alaydin2021updated}, the condition $\hat{\sigma}_{33} = 0$ (i.e., the plane-stress condition) is enforced by adjusting the through-thickness component of the rate of deformation $\hat{D}_{33}$. This is done iteratively using Newton's method, with the details provided in \cref{alg:zero_thick_stress}.
\begin{algorithm}
	\begin{algorithmic}[1]
	\State $\left(\hat{{\sigma}}_{ij}^{(0)}, \hat{\mathbb{C}}_{ijkl}^{alg, (0)} \right) \leftarrow \textbf{Elasto-plastic stress update using} ~ \hat{{D}}^{(0)}_{ij}$ (cf. \cref{sec:plasticity})
	\State $(\nu) \leftarrow (0)$
	\While{$\left|\hat{{\sigma}}_{33}^{(\nu)}\right| \geq tol $}
    \State $\Delta \hat{{D}}_{33} = -\dfrac{\hat{{\sigma}}^{(\nu)}_{33}}{\Delta t \, \hat{\mathbb{C}}_{3333}^{alg, (\nu)}}$
    \State $\hat{{D}}_{33}^{(\nu+1)} =  \hat{{D}}_{33}^{(\nu)} + \Delta \hat{{D}}_{33}$
    \State $\left(\hat{{\sigma}}_{ij}^{(\nu+1)}, \hat{\mathbb{C}}_{ijkl}^{alg,(\nu+1)} \right) \leftarrow \textbf{Elasto-plastic stress update with} ~ \hat{{D}}_{ij}^{(\nu+1)}$
    \State $(\nu) \leftarrow (\nu+1)$
    \EndWhile
	\end{algorithmic} 
	\caption{Enforcement of zero through-thickness stress} 	
	\label{alg:zero_thick_stress}
\end{algorithm}
In the \cref{alg:zero_thick_stress}, $\hat{\mathbb{C}}_{ijkl}^{alg}$ is the consistent tangent modulus in the unrotated configuration. For $J_2$ plasticity, $\hat{\mathbb{C}}_{ijkl}^{alg}$ is given by~\cite{simo1985consistent}:
\begin{equation}\label{eq:consistent_tangent_ep}
\begin{aligned}
\hat{\mathbb{C}}_{ijkl}^{alg} &= \left(\lambda + \frac{2\mu}{3}\right) \, \delta_{ij} \, \delta_{kl} + 2 \, \mu \,  \beta \, \mathbb{I}^{dev}_{ijkl} - 2 \, \mu \, \overline{\gamma}  \, \hat{n}_{ij} \, \hat{n}_{kl} , \\
\beta = \frac{\sigma_Y}{\hat{{\sigma}}^{vm}}, & ~~~~~ \overline{\gamma} = \frac{1}{1+\left(\dfrac{H}{3\mu}\right)}-\left(1-\beta\right), ~~~~~ \hat{n}_{ij} = \frac{\hat{{\sigma}}_{ij}^{dev}}{\sqrt{\hat{{\sigma}}_{kl}^{dev} \, \hat{{\sigma}}_{kl}^{dev}}} ,
\end{aligned}
\end{equation}
where $\mathbb{I}^{dev}_{ijkl} = \dfrac{1}{2}\left(\delta_{ik} \, \delta_{jl} + \delta_{il} \, \delta_{jk}\right) - \dfrac{1}{3} \, \delta_{ij} \, \delta_{kl}$ is the fourth-order deviatoric tensor.\\

Using the computed values of the rate of deformation tensor, the Jacobian of the deformation gradient and local shell thickness are obtained from the rate form of the corresponding equations. Applying midpoint time integration to both quantities gives the following update formula for the Jacobian:
\begin{equation} 
{J}^{(n+1)} = \left( \frac{1+\Delta t/2~\hat{{D}}_{kk}}{1-\Delta t/2~\hat{{D}}_{kk}} \right)~{J}^{(n)} ,
\label{eqn:J_update}
\end{equation}
and for the local thickness~\cite{alaydin2021updated}:
\begin{equation} 
h^{(n+1)} = \left( \frac{1+\Delta t/2~\hat{D}_{33}}{1-\Delta t/2~\hat{D}_{33}} \right)~h^{(n)} .
\label{eqn:h_update}
\end{equation}
Note that because the Jacobian is used to compute the Kirchhoff stress as per~\cref{eqn:tau-sigma}, the update given by~\cref{eqn:J_update} takes place at the bond level. Conversely, because the local thickness is a material-point variable, the update given by~\cref{eqn:h_update} takes place at the node level.

\subsection{Damage Modeling}
\label{sec:damage_modeling}

The bond-associative damage modeling approach~\cite{behzadinasab2020semi,behzadinasab2020revisiting} is utilized in this work. Compared to the state-based PD correspondence damage formulation~\cite{tupek2013approach}, the bond-associative approach applies a degradation function to the bonds instead of the nodes which results in enhanced stability and significant reduction of mass loss in damaged zones~\cite{behzadinasab2020semi}. In this approach, the influence function for a bond $\an{\eta}$ is modified as
\begin{equation}
\begin{aligned}
\omega_{\an{\eta}} &= \hat{\omega}\left(|\boldsymbol\eta|,\mathbf{q}_{\an{\eta}}\right) \\
&= \omega_{\eta}\left(|\boldsymbol\eta|\right) \ \omega_D\left(\, {D}_{\an{\eta}}\right) ,
\end{aligned}
\label{eqn:damagedOmega}
\end{equation}
where $\omega_{\eta}\left(|\boldsymbol\eta|\right)$ is the conventional damage-independent influence function providing the relative strength of interaction between two material points in the undamaged reference configuration, and $\omega_D\left(\, {D}_{\an{\eta}}\right)$ is the damage-dependent part of the influence state. The bond-associated damage parameter ${D}_{\an{\eta}} \in [0,1]$ depends on the bond-associated internal variables $\mathbf{q}$ (strain, stress, stretch, temperature, etc.) and is governed using classical continuum damage models. Damage irreversibility is included by requiring $\omega_D$ to be a non-increasing function of ${D}_{\an{\eta}}$. For an undamaged bond (${D}_{\an{\eta}}=0$), we require that $\omega_D(0) = 1$. For a fully-broken bond (${D}_{\an{\eta}}=1$), on the other hand, $\omega_D(1) = 0$. For visualization purposes, a node-associated damage variable is defined as the average damage at each material point, i.e.,
\begin{equation}
D_{\rm  (P)} = \dfrac{\sum_{{\rm Q}} \, {D}_{\rm \an{P-Q}}}{N_{\rm (P)}} ,
\end{equation}
where Q belongs to the neighbor set of P with the total number of neighbors $N_{\rm (P)}$. 

In this work, a cubic B-spline kernel is used to define the radial influence function, i.e.,
\begin{align}
  \omega_{\eta}(\hat\eta) = 
  \begin{cases}
    \dfrac{2}{3} - 4\,\hat\eta^2 + 4\hat\eta^3  ~~~ & \text{for} ~~~ 0 < \hat\eta \leq \dfrac{1}{2} \\
    \dfrac{4}{3} - 4\,\hat\eta + 4\,\hat\eta^2 - \dfrac{4}{3}\hat\eta^3  ~~~ & \text{for} ~~~ \dfrac{1}{2} < \hat\eta \leq 1 \\
    0  ~~~ & \text{otherwise} 
  \end{cases}
  , 
  \label{eqn:omega}
\end{align}
where $\hat\eta$ is defined for a bond $\rm \an{P-Q}$ as
\begin{equation}
  \hat\eta_{\rm \an{P-Q}} \equiv \frac{\left \vert X^{\rm 2D}_{\rm (Q)} - X^{\rm 2D}_{\rm (P)} \right \vert}{\delta} . 
\end{equation}

To model brittle fracture, a stress-based damage criterion (e.g., Tresca or von Mises) or a stretch-based failure theory (e.g., the critical stretch fracture model~\cite{silling2005meshfree}) can be incorporated, among others. To respect the immediate crack growth and unstable crack propagation nature of brittle materials, a bond is suddenly broken once its corresponding field variable exceeds a prescribed value. Alternatively, one can gradually degrade the bond stiffness.

To simulate ductile fracture, which is a gradual process in material degradation of malleable materials, we consider a plasticity-driven fracture model and decay the load-carrying capacity of a bond as follows:
\begin{equation}
\begin{aligned}
    & \omega_D \left( {D}_{\an{\eta}} \right) = 1 - {D}_{\an{\eta}} \ , \\
    & {D}_{\an{\eta}} = D \left( {\bar{\epsilon}}^P_{\an{\eta}} \right) =
    \begin{cases}
        0 , \qquad \qquad {\bar{\epsilon}}^P_{\an{\eta}} <= \bar{\epsilon}^P_{\rm th} , \\
        \dfrac{{\bar{\epsilon}}^P_{\an{\eta}} - \bar{\epsilon}^P_{\rm th}}{\bar{\epsilon}^P_{\rm cr} - \bar{\epsilon}^P_{\rm th}} , \ \quad \bar{\epsilon}^P_{\rm th} < {\bar{\epsilon}}^P_{\an{\eta}} < \bar{\epsilon}^P_{\rm cr} , \\
        1 , \qquad \qquad {\bar{\epsilon}}^P_{\an{\eta}} >= \bar{\epsilon}^P_{\rm cr} ,
    \end{cases}
\end{aligned}
\label{eqn:plastic_strain}
\end{equation}
where ${\bar{\epsilon}}^P_{\an{\eta}}$ is the bond-associated equivalent plastic strain. The model parameters $\bar{\epsilon}^P_{\rm th}$ and $\bar{\epsilon}^P_{\rm cr}$ are calibrated for a given material. In this approach, an undamaged bond has ${\bar{\epsilon}}^P_{\an{\eta}} \leq \bar{\epsilon}^P_{\rm th}$, and a fully-damaged bond has ${\bar{\epsilon}}^P_{\an{\eta}} = \bar{\epsilon}^P_{\rm cr}$.  Other, more sophisticated failure criteria may be incorporated in the modeling, e.g., the Johnson-Cook fracture model~\cite{johnson1985fracture} (also see~\cite{behzadinasab2020semi}).

\section{Numerical Results}
\label{sec:results}

In this section, several numerical case studies including some shell benchmark problems (e.g. the shell obstacle course ~\cite{belytschko1985stress}) are provided to highlight the capabilities of the developed formulation. The examples cover small-strain elasticity, large-deformation elasto-plasticity, and crack growth problems (brittle and ductile fracture). In the following examples, we consider the asymptotic convergence of the proposed model (also known as the $\delta$-convergence in PD \cite{bobaru2009convergence}) to reference local solutions (e.g., IGA simulations). Thus, the horizon size and nodal spacing approach zero at the same rate. Quadratic RK shape functions are employed in the computations unless otherwise noted. The horizon size is chosen to be $p+1$ times the average node spacing for RK shape functions of order $p$. Following~\cite{alaydin2021updated}, we use three Gauss points in the through-thickness direction, which results in
\begin{equation}
  \left\{ \left( \xi_3^I, w^I \right) \right\} = \left\{ \left( 0, \, \frac{8}{9} \right), \left( -\sqrt{\frac{3}{5}}, \, \frac{5}{9} \right), \left( \sqrt{\frac{3}{5}}, \, \frac{5}{9} \right) \right\} .
\end{equation}
The normalized influence state $\alpha_{\rm \an{P-Q}}$ is defined as
\begin{equation}
\alpha_{\rm \an{P-Q}} = \frac{\omega_{\rm \an{P-Q}}}{\omega_{0_{(\rm P)}}}, \qquad \omega_{0_{(\rm P)}} = \int_{\mathcal{B}^{\rm 2D}} \omega_{{\eta}_{\rm \an{P-Q}}} \, {\rm dQ} .
\end{equation}
For visualization purposes, the nodal field variables (strains, stresses, etc.) are computed by applying the kinematic and constitutive relations at the nodes.

This section contains three main parts. We first study the accuracy of the local surface parameterization in approximating curved geometries in \cref{sec:geometry_tests}. Next, in \cref{sec:elastostatics}, we carry out a set of elastostatic benchmark simulations to assess the convergence behavior of the PD formulation. The computations are done using an in-house Python code involving a quasi-static solver with Backward Euler time integration and a linearized version of the thin shell formulation detailed in \cref{sec:linearization}. The remainder of the verification simulations involving plasticity and fracture in \cref{sec:plasticity_problems} are carried out using the explicit dynamics solver (Velocity Verlet integration) of an extended version of the open-source PD code Peridigm~\cite{parks2012peridigm}. \cref{sec:fragmentation} shows two demonstrative examples that involve severe material damage and fragmentation. For the explicit calculations, the stability analysis of~\cite{silling2005meshfree} gives a good estimate of the critical time step. Some quasi-static problems are solved using the explicit solver. In these cases, in order to speed up the computations, the loading rates are scaled up while respecting the fact that the inertial forces must remain small to achieve quasi-static response. In addition, to suppress the elastic transients in these cases, the velocity boundary conditions are ramped up smoothly from zero to the desired values over a short time period.

\subsection{Geometry Tests}
\label{sec:geometry_tests}

\begin{figure*}[!hbpt]
    \centering
    \subfloat[][]{\includegraphics[width=0.5\columnwidth]{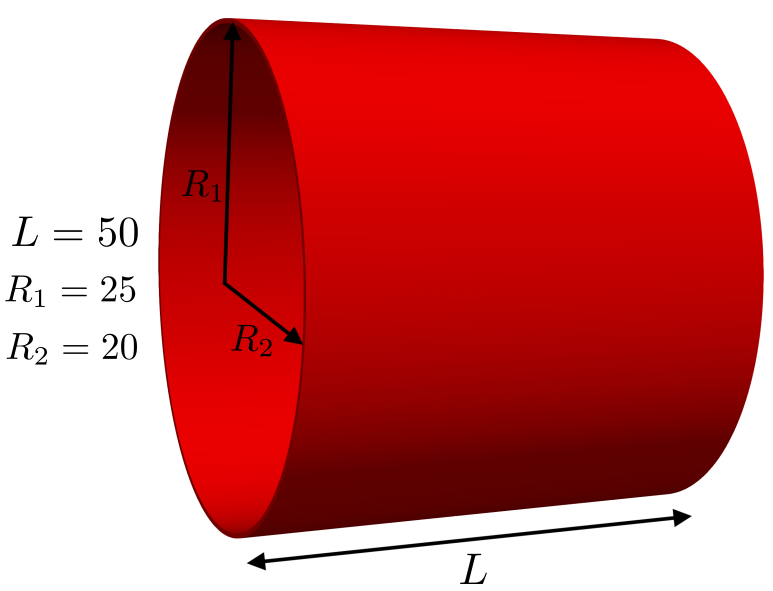}}
    \hfill
    \subfloat[][]{\includegraphics[width=0.45\columnwidth]{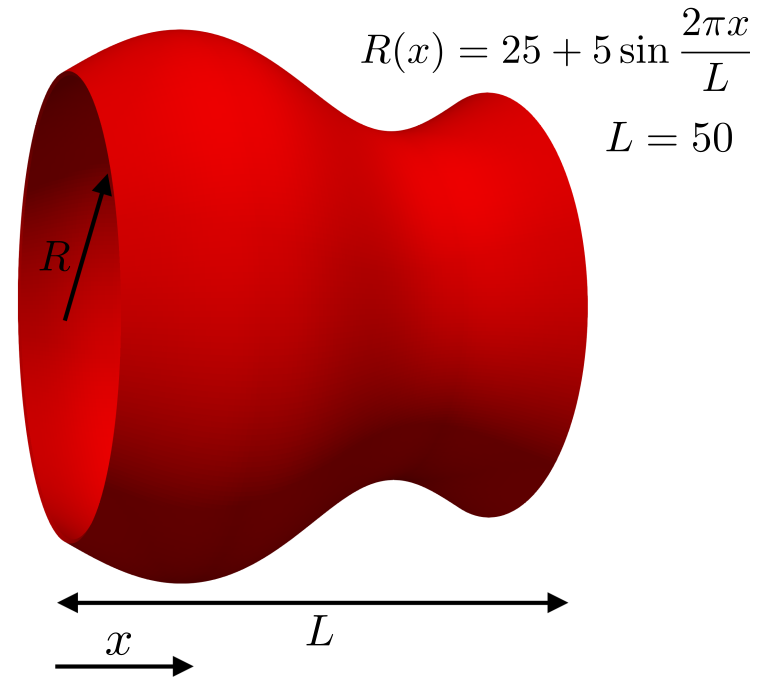}}
    \caption{The geometry tests. (a): Cylinder with an elliptical cross-section area. (b): Vase with a circular cross-section area of varying diameter.}
    \label{fig:geometries}
\end{figure*}

In this section we assess the accuracy of our approach in representing curved geometries. Note that we do not solve the equations of motion here, but only calculate the geometric quantities of interest, i.e., the normal vector and principal curvatures. An elliptical cylinder and a vase-like shape are considered (see~\cref{fig:geometries})
where the normal and curvatures may be computed analytically and used to assess the approximation error. In each case, the local manifolds are reconstructed using the PCA algorithm as described in \cref{sec:manifold}. Then, the parametric derivatives of the reference positions are computed using the PD gradient operator as in \cref{sec:PD_gradient}. Quadratic, cubic, and quartic RK shape functions are considered. The normal vector is computed using \cref{eqn:normal}, which involves only the first-order parametric derivatives. The principal curvatures are calculated as detailed in \cite[Section 8.2]{pressley2010elementary}, which involves both first- and second-order gradients. The normal and curvature errors are calculated using the root-mean-square (RMS) norm
\begin{equation}
  ||e||_{2} = \sqrt{\frac{\sum_i^N e_i^2}{N}}.
\end{equation}
%

\begin{figure*}[!hbpt]
    \centering
    \subfloat[][RMS error in normals]{\includegraphics[width=0.49\columnwidth]{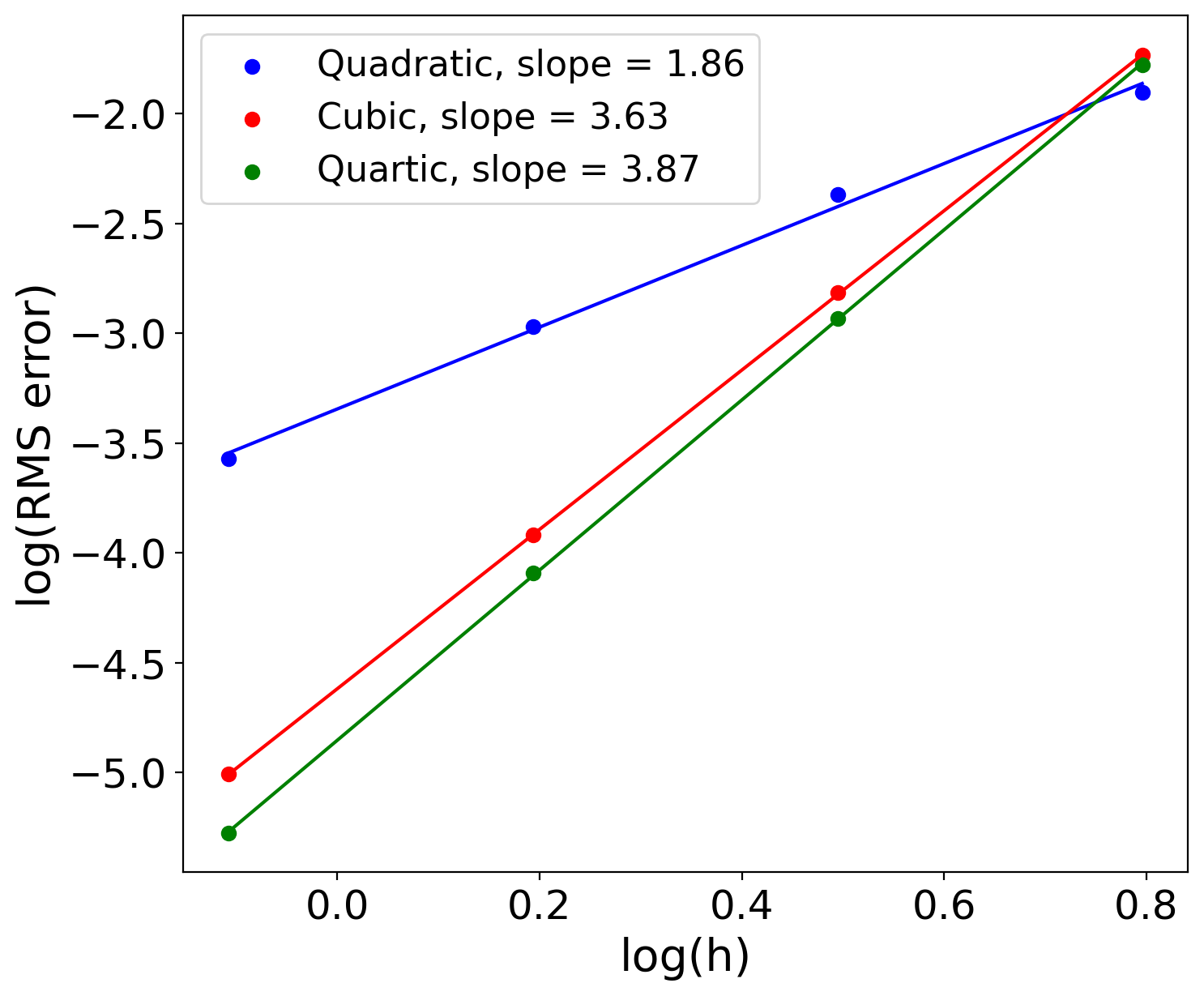}}
    \hfill
    \subfloat[][RMS error in curvatures]{\includegraphics[width=0.49\columnwidth]{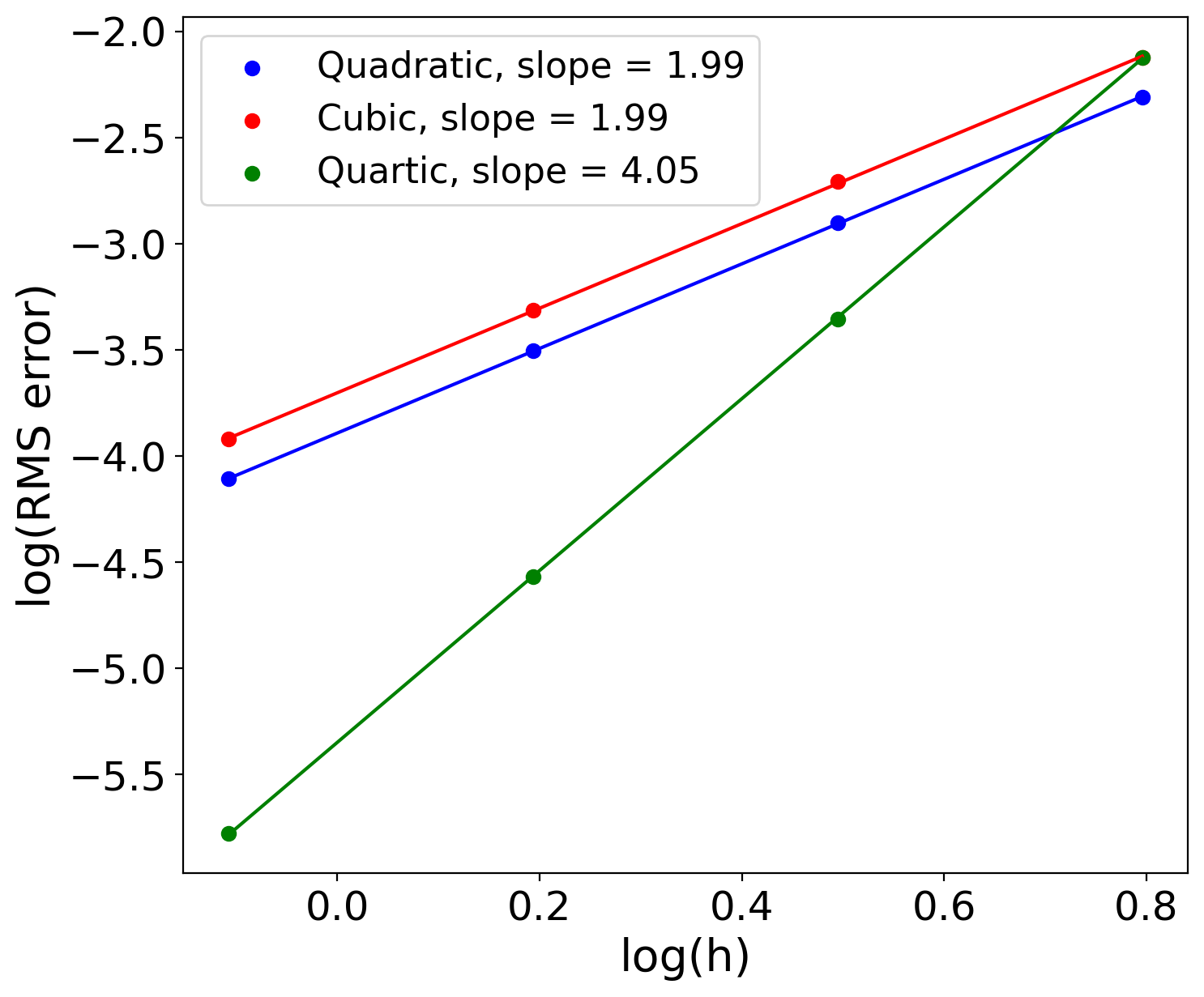}}
    \caption{Geometry test on the elliptical cylinder. RMS error analysis of (a): normals and (b): curvatures.}
    \label{fig:geometry_test_cylinder}
\end{figure*}

\begin{figure*}[!hbpt]
    \centering
    \subfloat[][RMS error in normals]{\includegraphics[width=0.49\columnwidth]{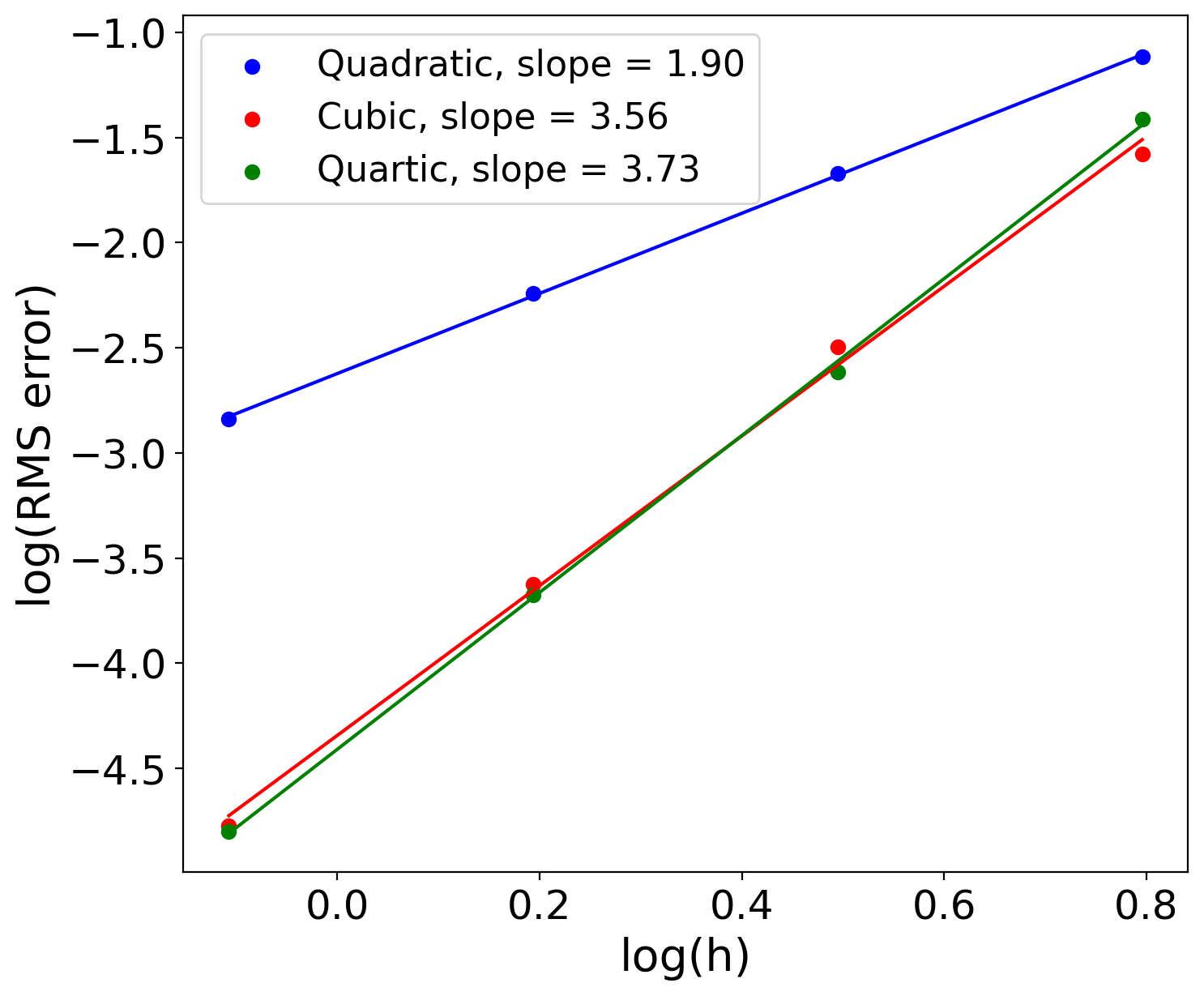}}
    \hfill
    \subfloat[][RMS error in curvatures]{\includegraphics[width=0.49\columnwidth]{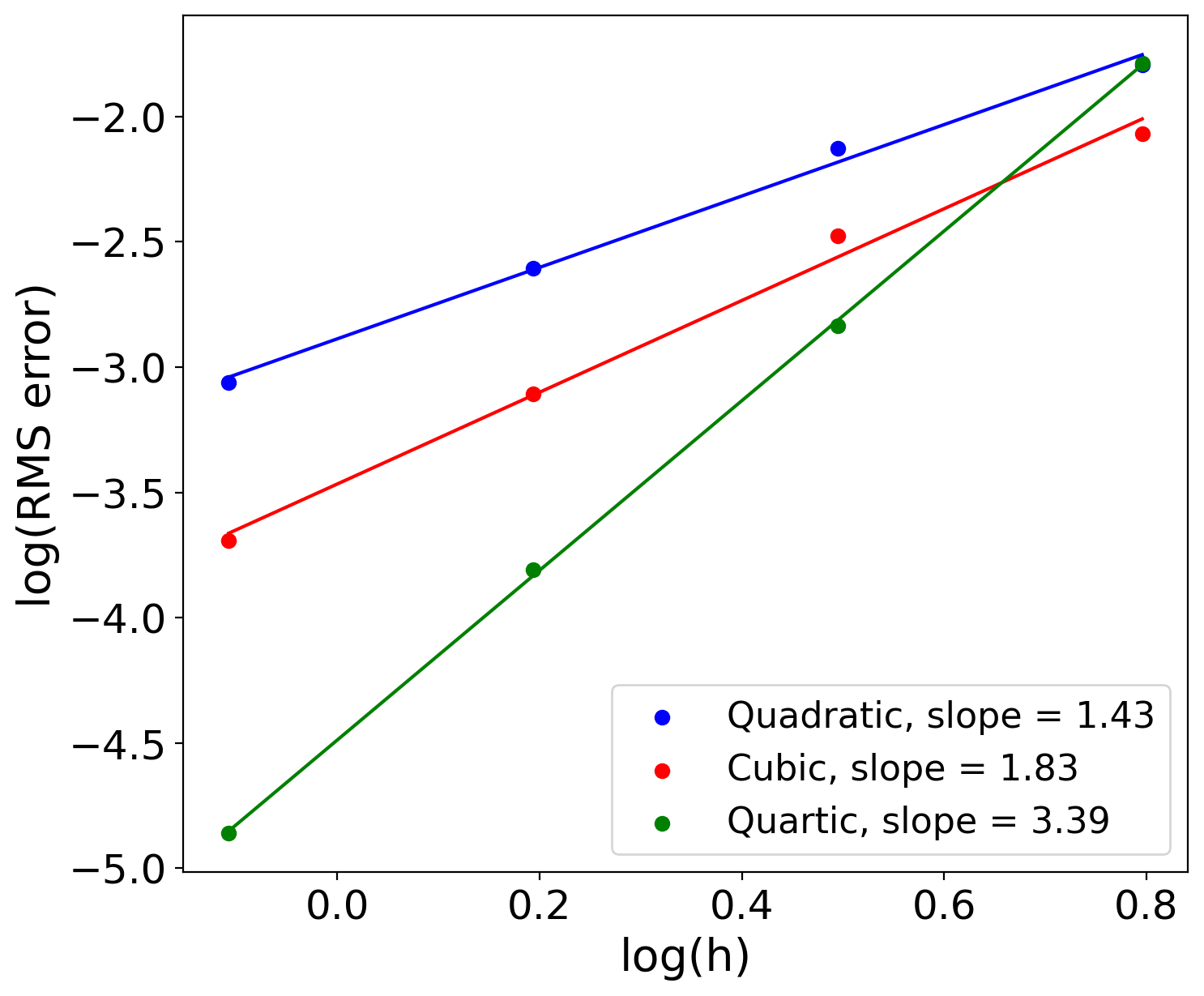}}
    \caption{Geometry test on the vase. RMS error analysis of (a): normals and (b): curvatures.}
    \label{fig:geometry_test_vase}
\end{figure*}

Both geometries are meshed in a quasi-uniform fashion such that the number of nodes along each length is 8, 16, 32, and 64 at each level of discretization. \cref{fig:geometry_test_cylinder}
shows the RMS error in the normal vector and curvatures plotted against the average nodal spacing for the elliptical cylinder. \cref{fig:geometry_test_vase} shows the same data for the vase-like shape. For the normal vector, the quadratic discretization shows a nearly second-order convergence rate, while the cubic and quartic discretizations converge at a rate between 3.5 and 4. For the curvatures, the quadratic and cubic discretizations converge at the rate between 1.5 and 2, while the quartic dicretizations converge at the rate between 3.5 and 4. The results indicate the presence of super-convergence, which is not uncommon for meshfree schemes~\cite{leng2019super,trask2019asymptotically,hillman2020generalized,behzadinasab2021unifiedI}. In particular, fwhen calculating the normal vector, which involves only the first-order derivatives, the super-convergence behavior happens for the odd orders (i.e., cubic). On the other hand, the super-convergence behavior shifts to the even orders (i.e., quadratic and quartic) for computing the curvatures, in which both the first- and second-order derivatives are involved.  

\subsection{Elastostatic Problems}
\label{sec:elastostatics}

Several elastostatic benchmark problems, including those from the shell obstacle course~\cite{belytschko1985stress}, are studied in this section. 

\subsubsection{Clamped Circular Plate under Uniform Transverse Loading}
\label{sec:clamped_circular_plate}

\begin{figure*}[!hbpt]
    \begin{minipage}{.325\columnwidth}
    \centering
    \subfloat[][PD Solution]{\includegraphics[width=\columnwidth]{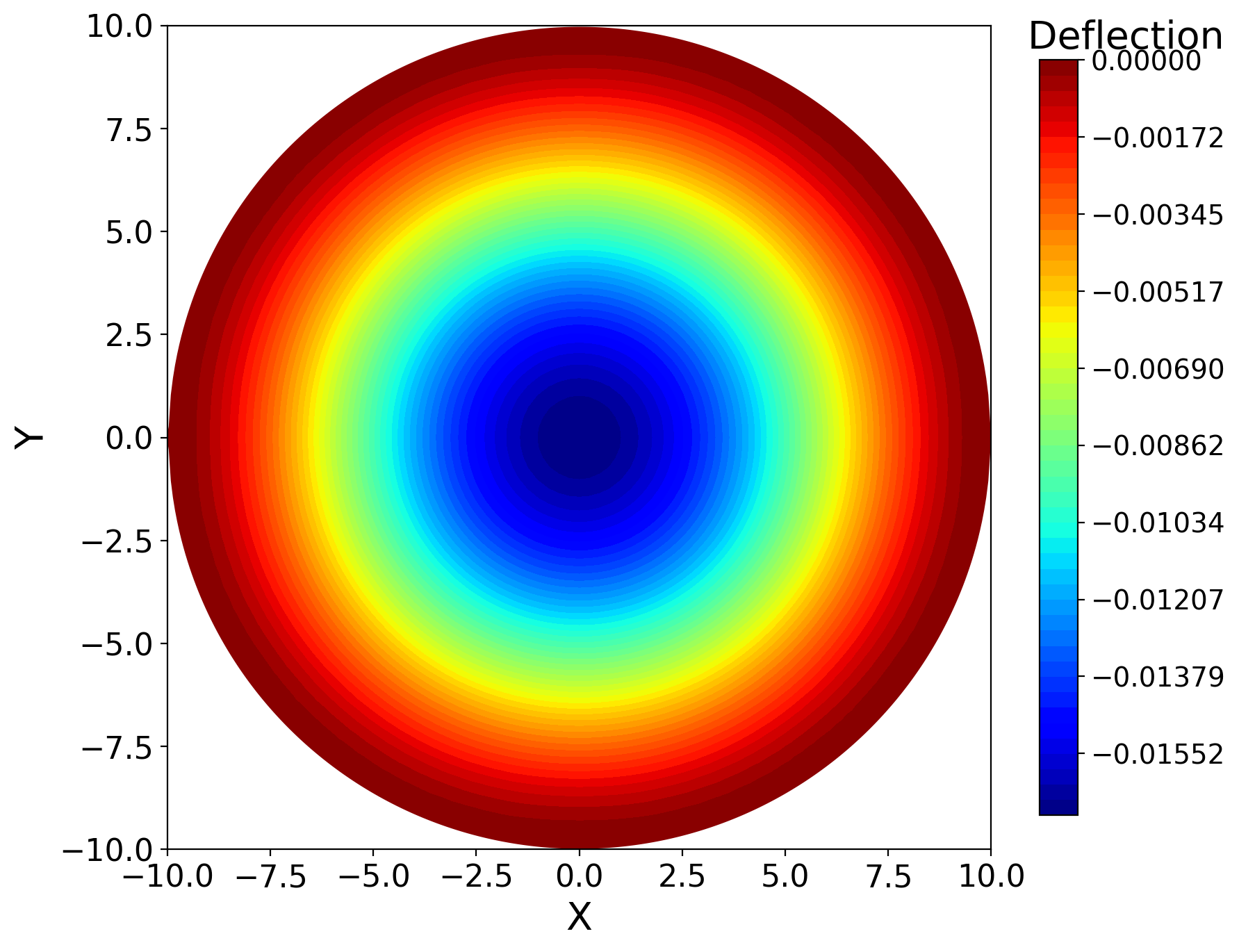}}
    
    \subfloat[][Exact Solution]{\includegraphics[width=\columnwidth]{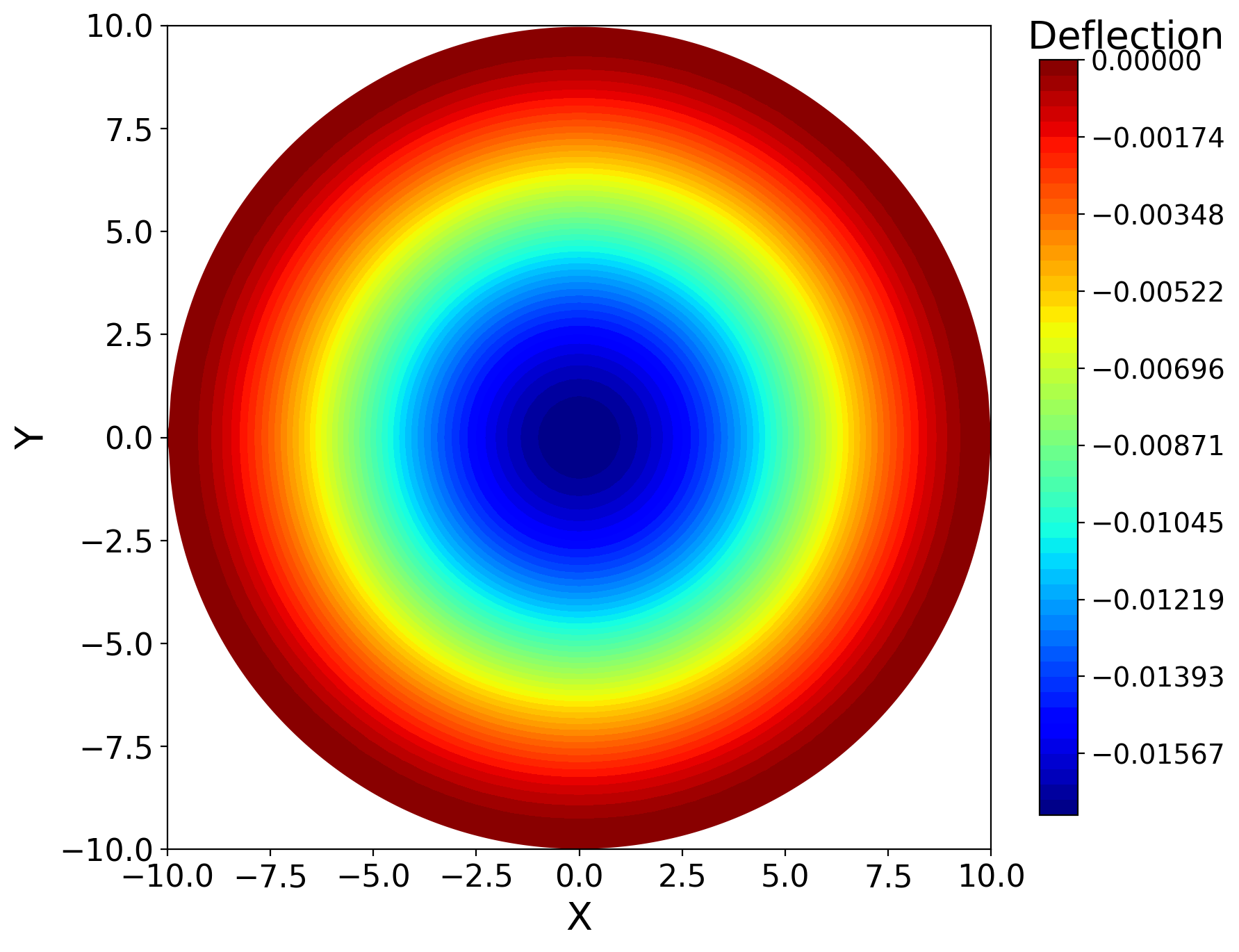}}
    \end{minipage}
    \hfill
    \begin{minipage}{.59\columnwidth}
    \centering
    \subfloat[][Convergence study]{\includegraphics[width=\columnwidth]{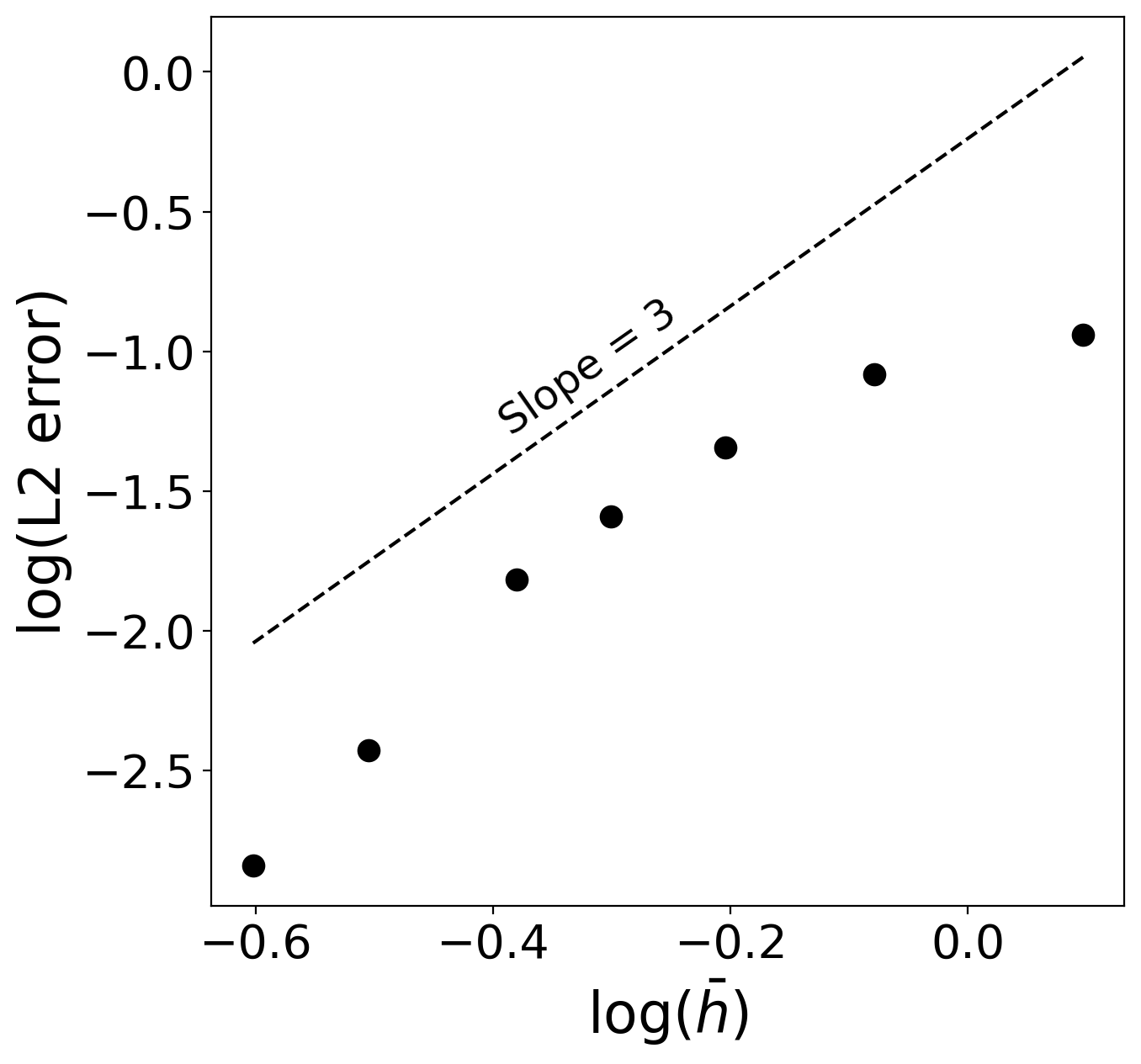}}
    \end{minipage}
    \caption{Clamped circular plate under uniform transverse loading. (a): Transverse displacement of the plate using 40 PD nodes along the radius, (b): Transverse displacement of the plate corresponding to the exact solution, and (c): $L_2$ norm of the displacement error vs. the average nodal spacing.}
    \label{fig:clamped_circular_plate}
\end{figure*}

For the first example, we study a problem with an existing analytical solution. Consider a circular plate (flat shell) that is clamped at its edges and subjected to a uniform transverse load $q$. The exact deflection using the classical theory is given by~\cite{kelly2021solid}:
\begin{equation}
    w = - \frac{q}{64D}(x^2+y^2-R^2),
\end{equation}
where $x$ and $y$ are in-plane coordinates, $R$ is the plate radius, and $D = \dfrac{E h^3}{12 \, (1-\nu^2)}$ is the flexural rigidity of the plate with thickness $h$, Young's modulus $E$, and Poisson's ratio $\nu$. In our setup, $R=10$, $h=0.1$, $E=10^5$, $\nu=0.3$, and $q=10^{-3}$. Quadratic discretizations with 8, 12, 16, 20, 24, 32, and 40 nodes along the radius of the plate are considered. Clamped boundary conditions are enforced by fixing two rows of PD nodes around the plate circumference. As shown in~\cref{fig:clamped_circular_plate}, the PD refined solution is in good agreement with the analytical solution. A displacement convergence rate higher than three is achieved in the asymptotic regime.

\subsubsection{Simply Supported Plate under Sinusoidal Pressure Loading}
\label{sec:plate_under_sinusoidal_pressure}

\begin{figure*}[!hbpt]
    \centering
    \subfloat[][]{\includegraphics[width=0.47\columnwidth]{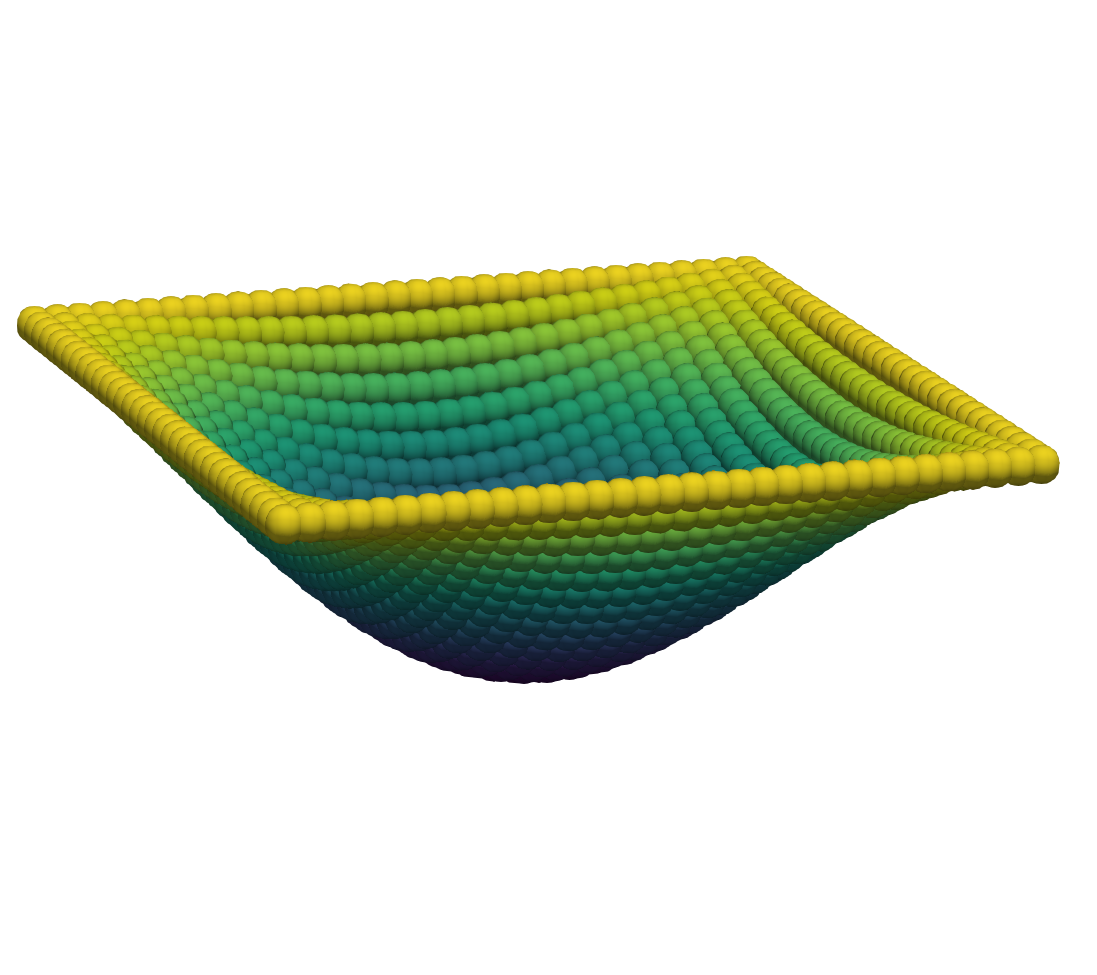}}
    \hfill
    \subfloat[][]{\includegraphics[width=0.48\columnwidth]{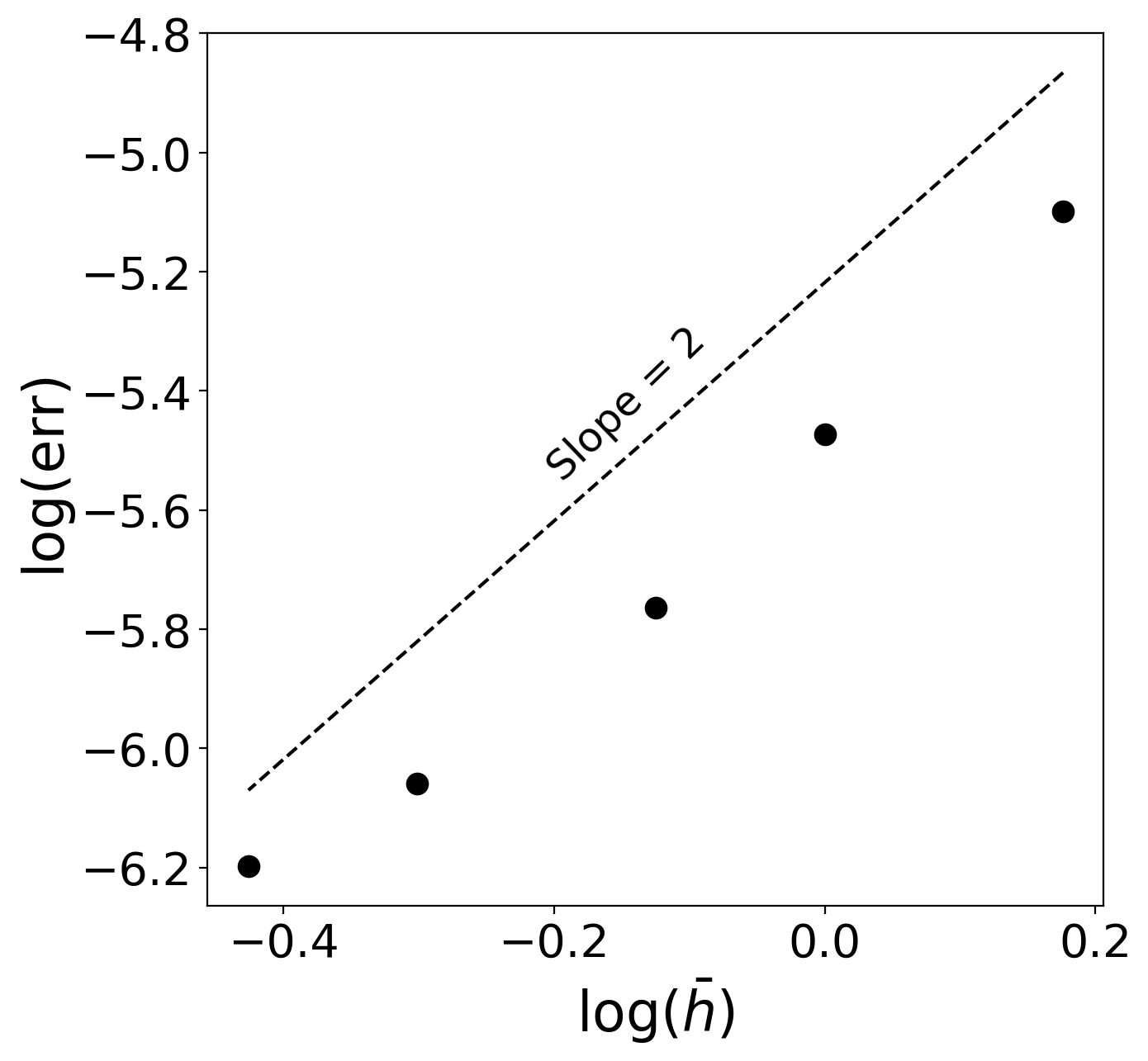}}
    \caption{Plate under sinusoidal pressure load. (a): Deformed configuration for the $32 \times 32$ PD node discretization (displacements are amplified by a factor of $2 \times 10^5$) and (b): Maximum deflection error vs. the nodal spacing.}
    \label{fig:plate_sinusoidal_error}
\end{figure*}

Next, we consider a problem outlined in~\cite{duong2017new} in which a plate with size $(0,L) \times (0,L)$, thickness $h=0.375$, Young's modulus $E = 4.8 \times 10^5$, and Poisson's ratio $\nu = 0.38$ is subjected to a sinusoidal pressure loading $p(X,Y) = p_0 \, \sin(\pi X/L) \, \cos(\pi Y/L)$. The plate is simply supported along its edges. According to the Navier solution~\cite{ugural2009stresses}, the maximum displacement occurs at the plate center and is given by
\begin{equation}
    w_{\max} = \frac{p_0 L^4}{4 \pi^4 D},
\end{equation}
where $D = \dfrac{E h^3}{12 \, (1-\nu^2)}$ is the flexural rigidity of the plate. In our setup, the plate is discretized using $N \times N$ peridynamic nodes. Zero displacements are prescribed along the edge PD nodes. To stay in the small displacement regime, a small value of $p_0 = 10^{-3}$ is chosen. Refinements corresponding to $N = 8, 12, 16, 20, 32$ are considered. The final configuration of the plate and the error in $w_{\max}$ are shown in \cref{fig:plate_sinusoidal_error}. A near second-order convergence rate in the maximum norm is obtained. 

\subsubsection{Clamped Plate Under Transverse Loading at the Free Edges}
\label{sec:clamped_plate}

\begin{figure}[!hbpt]
    \centering
    \includegraphics[width=0.5\columnwidth]{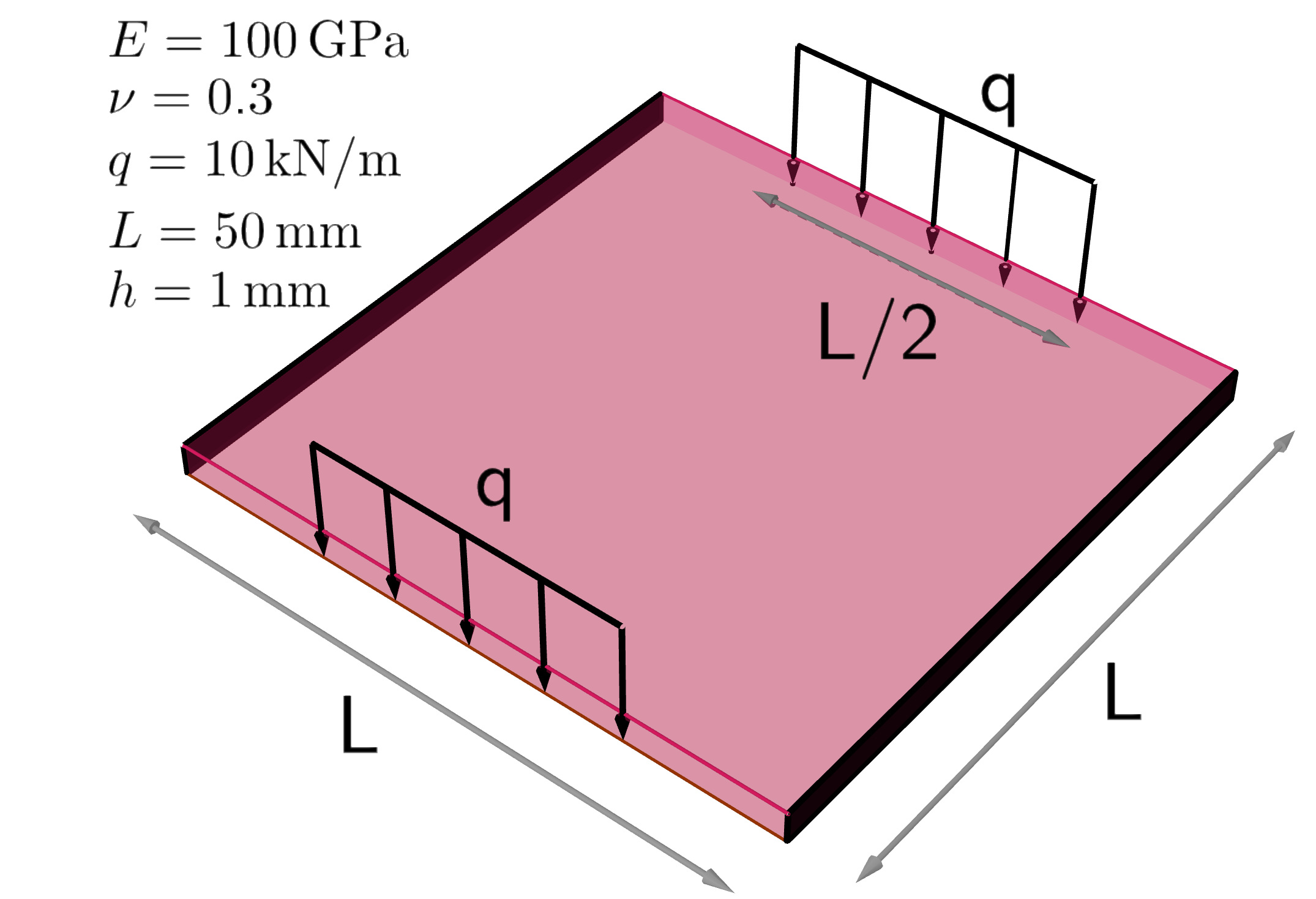}
    \caption{Clamped plate under transverse loading at the free edges. Problem setup.}
    \label{fig:clamped_plate_setup}
\end{figure}

\begin{figure*}[!hbpt]
    \centering    
    \subfloat[][IGA]{\scalebox{-1}[1]{\includegraphics[trim=15cm 7cm 12cm 14cm, clip=true, width=0.48\columnwidth]{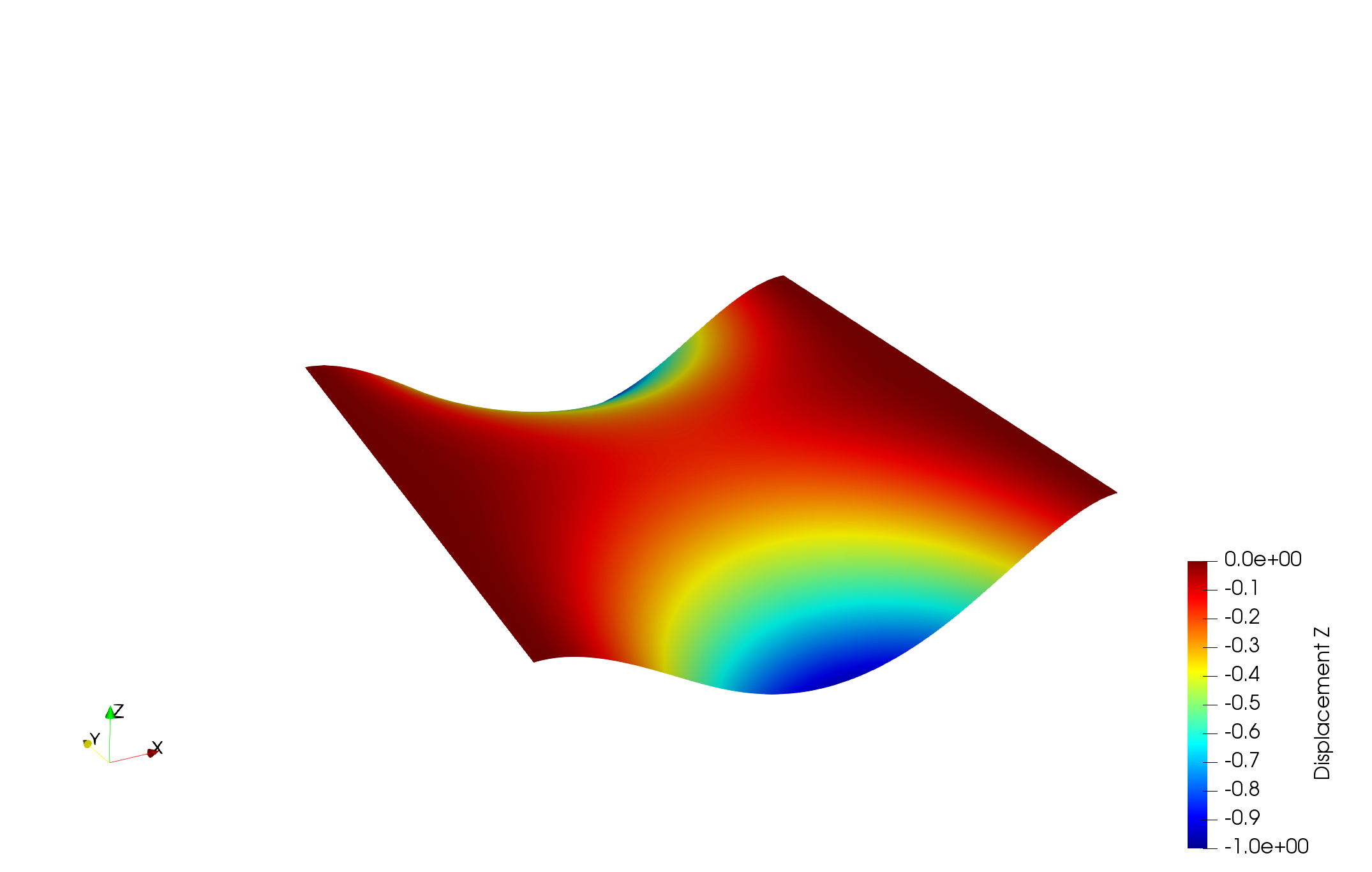}}}
    \subfloat[][PD]{\includegraphics[width=0.48\columnwidth]{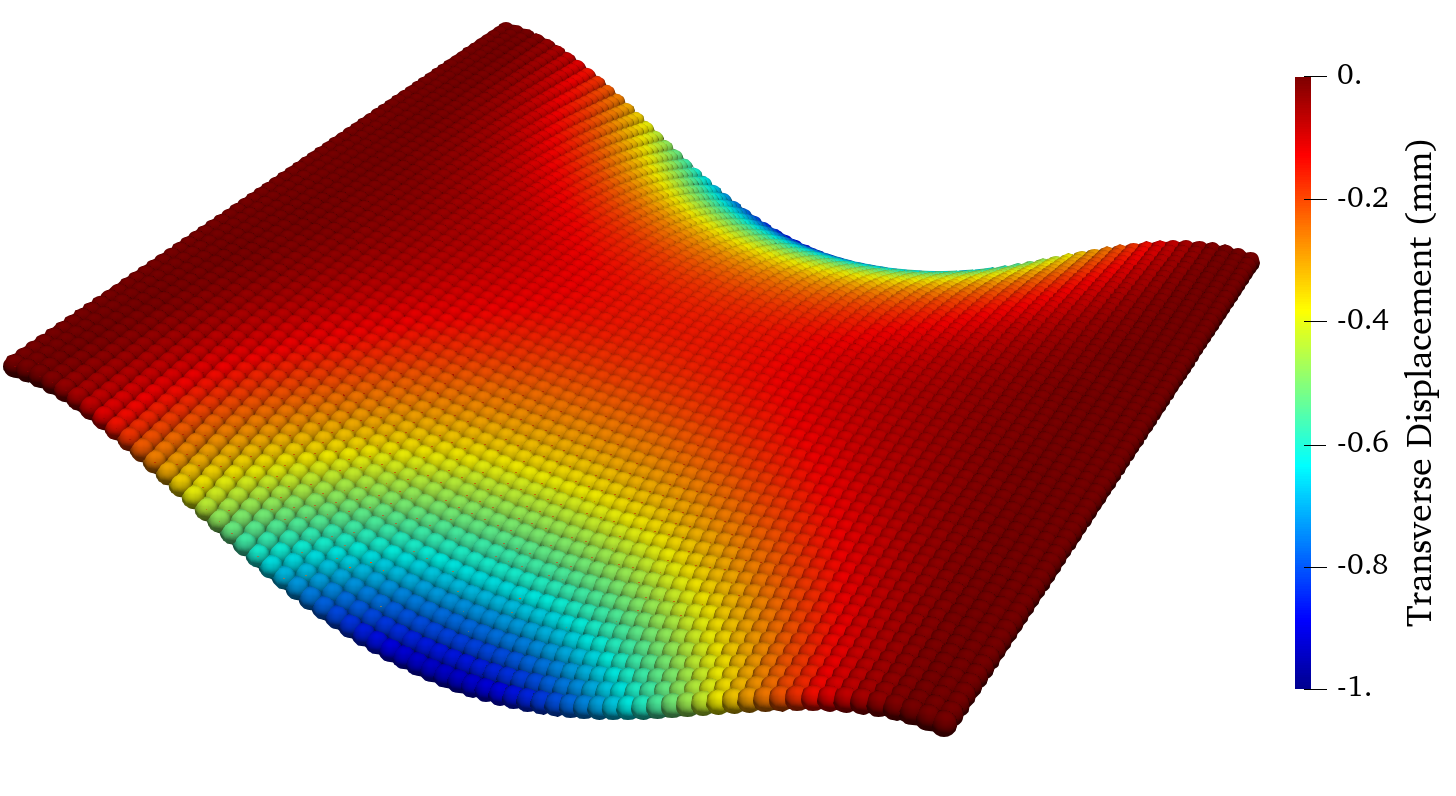}}
    
    \subfloat[][IGA]{\includegraphics[height=0.47\columnwidth]{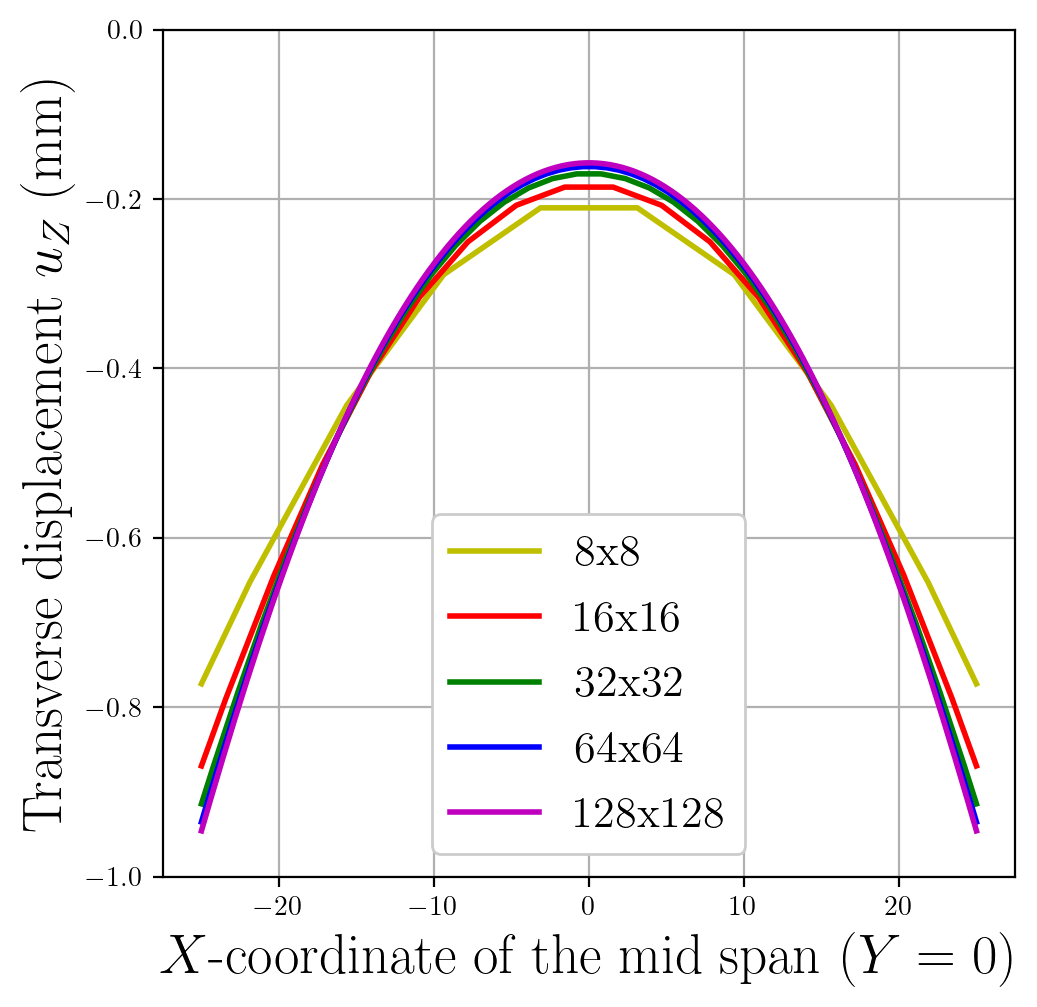}}
    \hspace*{4pt}
    \subfloat[][PD]{\includegraphics[trim = 2.51cm 0 0 0, clip, height=0.47\columnwidth]{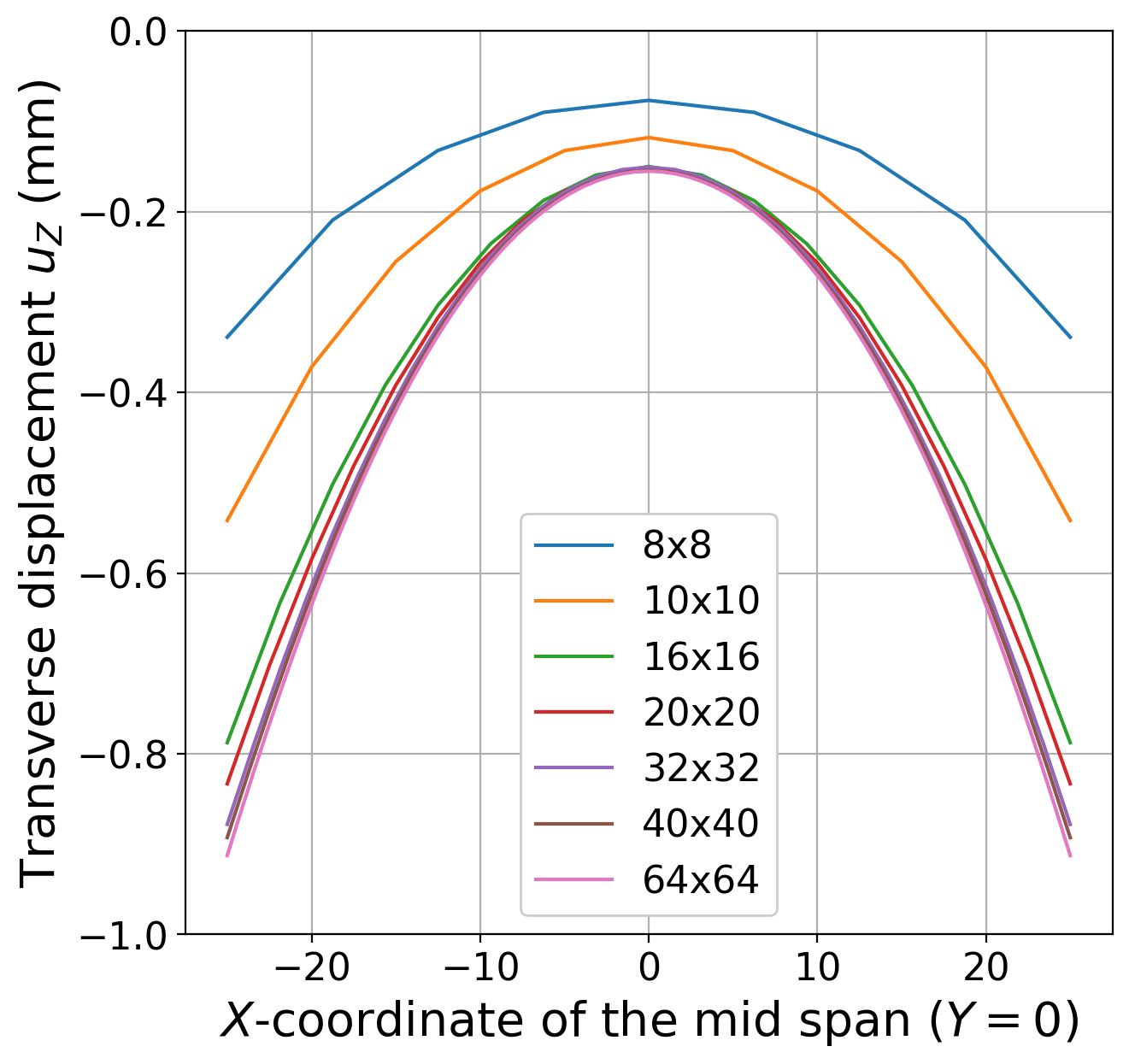}}
    \caption{Clamped plate under transverse loading at the free edges. Deformed shape of the plate using (a): 64x64 IGA discretization and (b): 64x64 PD nodal discretization (a scaling factor of 10 is used for the displacements). Convergence behavior of the mid-span deflections for (c): IGA, and (d): PD.}
    \label{fig:clamped_plate}
\end{figure*}

In this section, we compute a case with no analytical solution, instead using a highly resolved IGA solution for reference. The problem setup is shown in~\cref{fig:clamped_plate_setup}, where a square-shaped flat plate is clamped at two of its edges and subjected to a (non-uniform) transverse load applied to the mid-section of the free edges. The applied lateral force is treated as a body force for the PD computation. Clamped boundary conditions are enforced by fixing two rows of nodes along the corresponding edges. Quadratic IGA and PD discretizations with varying mesh resolution are employed for this problem. \cref{fig:clamped_plate}(a--b) present the deformed shape of the plate using the converged IGA and PD solutions. \cref{fig:clamped_plate}(c) provides the mid-span deflections for different discretization levels. The converged solutions for both methods are in good agreement.

\subsubsection{Scordelis–Lo Roof}
\label{sec:scordelis_roof}

\begin{figure}[!hbpt]
    \centering
    \includegraphics[width=0.5\columnwidth]{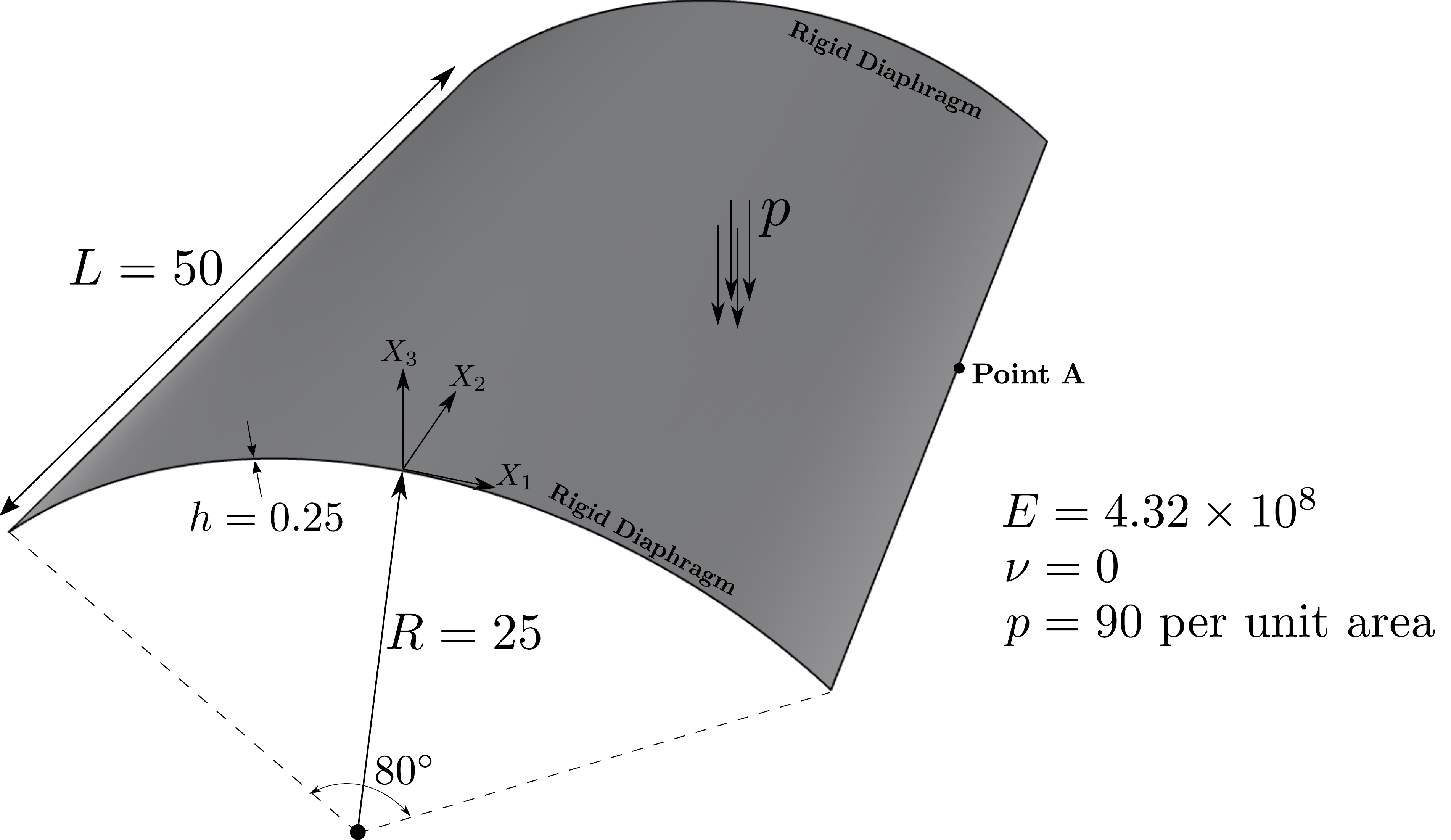}
    \caption{The setup for the Scordelis-Lo roof problem. Point A corresponds to the midpoint of the free edge.}
    \label{fig:scordelis_setup}
\end{figure}

\begin{figure*}[!hbpt]
    \begin{minipage}{.35\columnwidth}
    \centering
    \subfloat{\includegraphics[width=0.8\columnwidth]{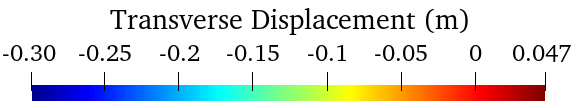}}
    
    \vspace*{-0.1cm}
    \setcounter{subfigure}{0}
    \subfloat[][IGA solution]{\includegraphics[trim=0 0 6cm 0, clip=true,width=\columnwidth,height=0.7\columnwidth]{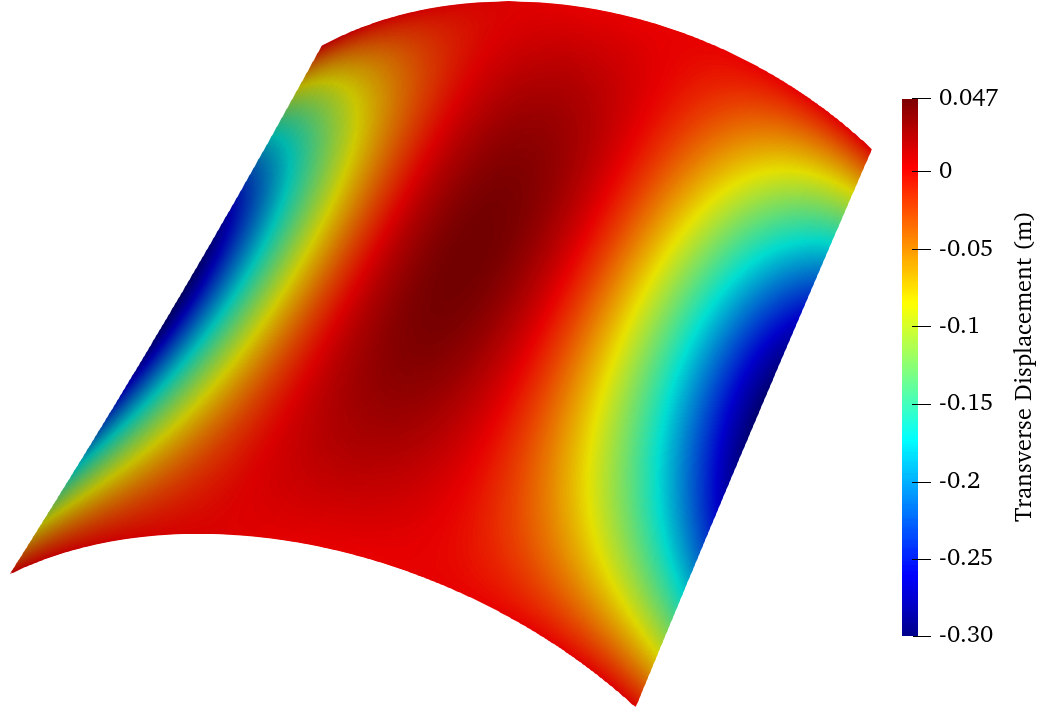}}
    
    \subfloat[][PD solution]{\includegraphics[trim=0 0 9cm 0, clip=true,width=\columnwidth,height=0.7\columnwidth]{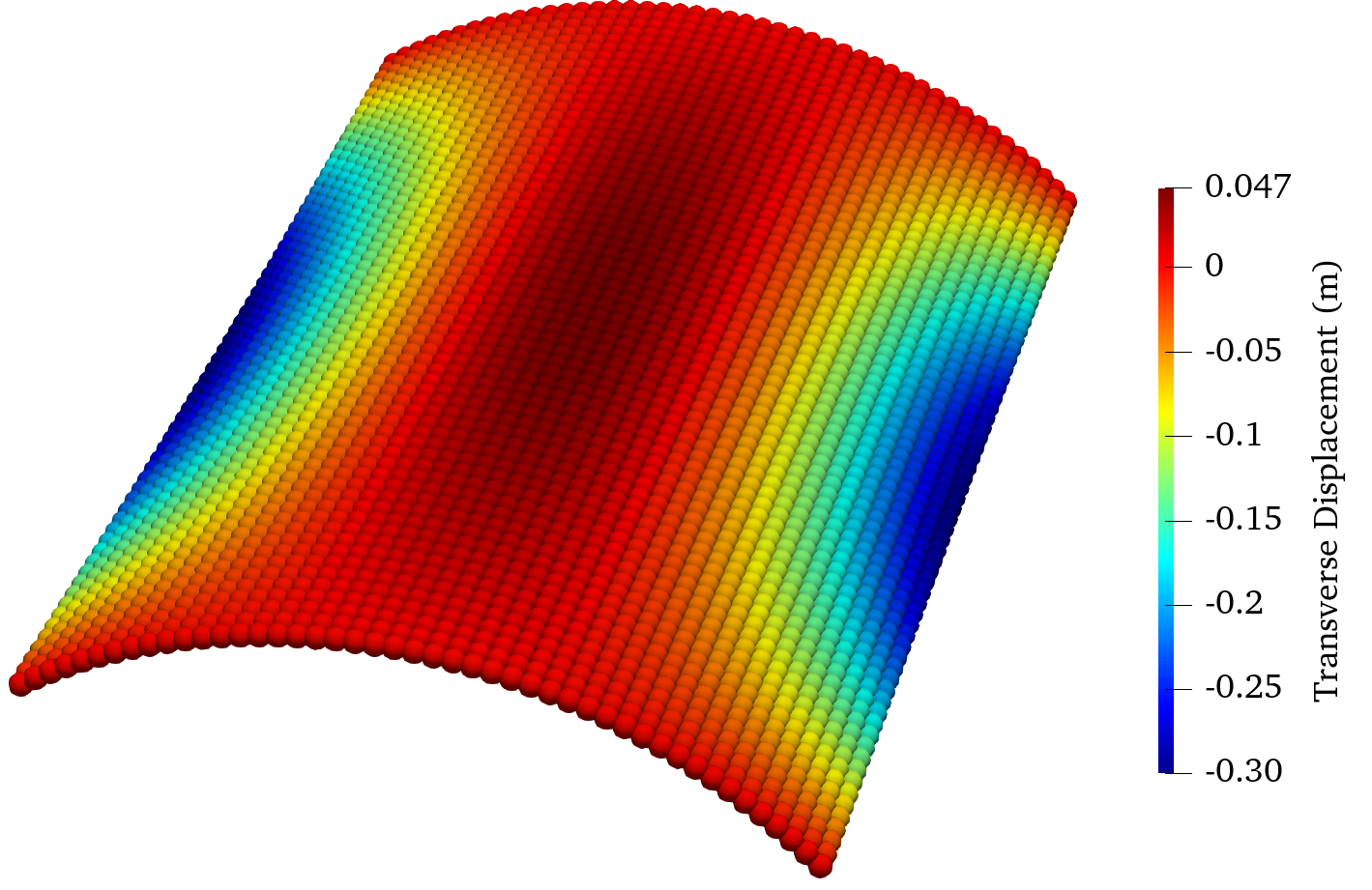}}
    \end{minipage}
    \hfill
    \begin{minipage}{.6\columnwidth}
    \centering
    \subfloat[][Convergence study]{\includegraphics[width=\columnwidth]{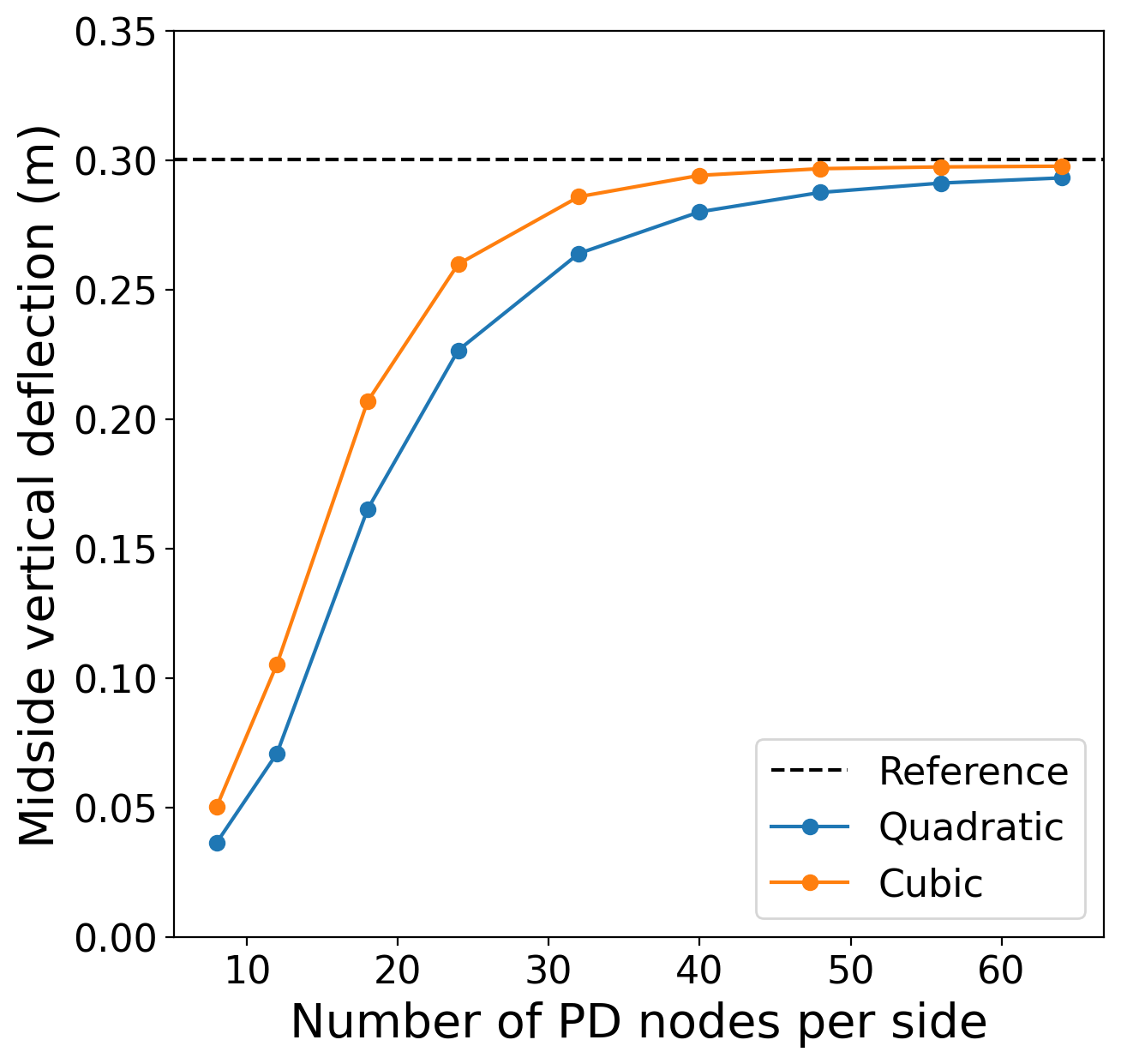}}
    \end{minipage}
    \caption{Scordelis–Lo roof under gravity. Deformed shape of the roof using (a): IGA and (b): PD. (c): Convergence of the free-edge midpoint displacement using quadratic and cubic RK shape functions in the PD discretization.}
    \label{fig:scordelis}
\end{figure*}

The Scordelis–Lo roof problem was proposed as part of the so-called {\em shell obstacle course}~\cite{belytschko1985stress} to challenge shell formulations in their ability to capture complicated membrane stresses and has been widely used as a benchmark problem in the realm of shells~\cite{macneal1985proposed,kiendl2009isogeometric,kiendl2010bending,coox2017flexible,maurin2018isogeometric,zhang2021peridynamic}. The Scordelis-Lo roof is part of a cylindrical shell, which is supported by rigid end diaphragms (constrained in-plane displacements) and has traction-free side edges. The roof is subjected to a uniform gravity load in the vertical direction. The problem description is shown in \cref{fig:scordelis_setup}. 
The vertical displacement at the midpoint of the free edge is reported as 0.3006 m~\cite{hughes2005isogeometric,kiendl2009isogeometric,coox2017flexible,maurin2018isogeometric}. Good agreement between the PD and IGA solutions of the deformed shape of the roof is obtained and shown in \cref{fig:scordelis}(a--b). As shown in \cref{fig:scordelis}(c), cubic RK shape functions exhibit faster convergence to the reference solution than their quadratic counterparts.

\subsubsection{Pinched Cylinder}
\label{sec:pinched_cylinder}

\begin{figure}[!hbpt]
    \centering
    \includegraphics[width=0.45\columnwidth]{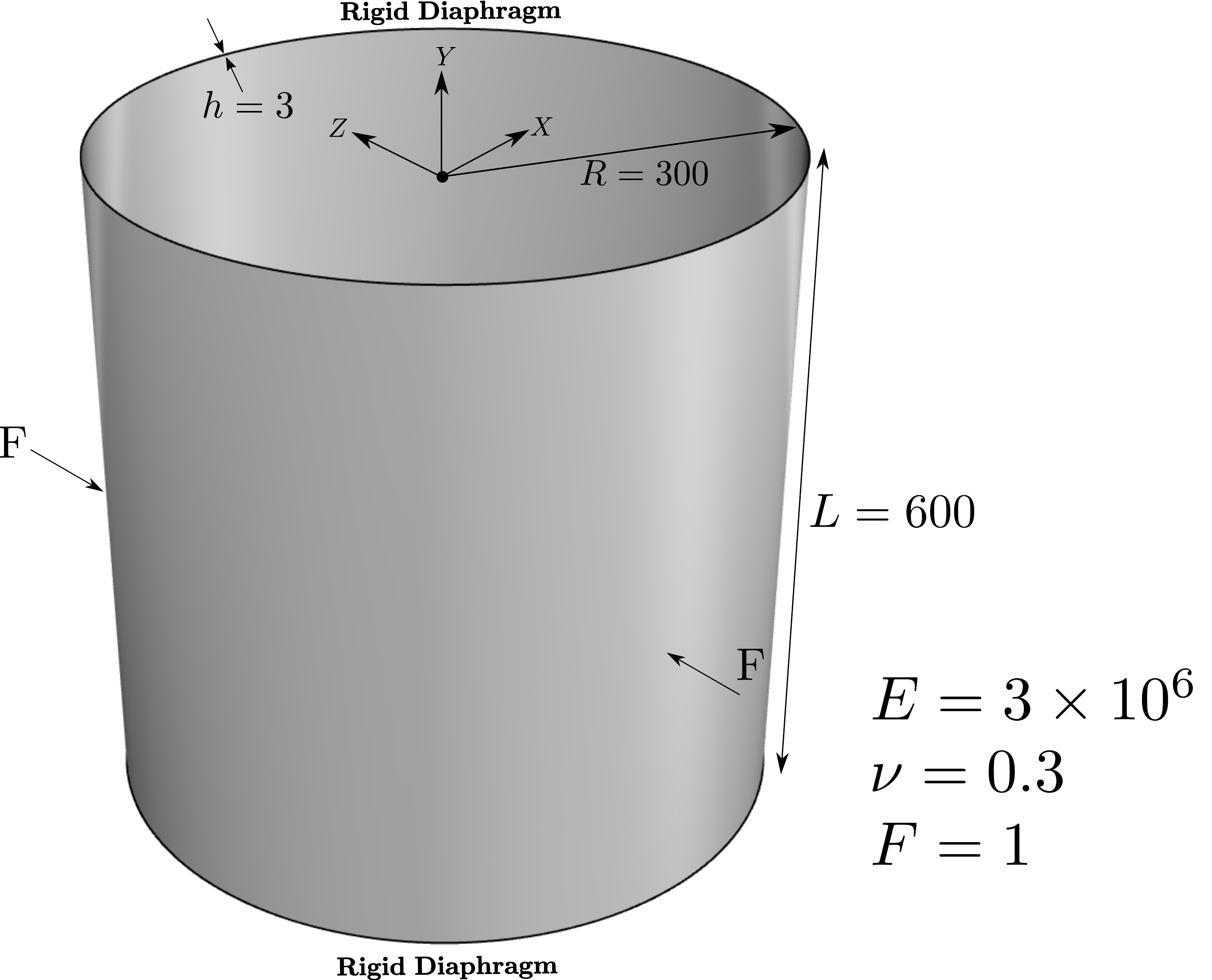}
    \caption{Setup for the pinched cylinder problem. The point loads are applied in the $Z$ direction at the points $X=0, \ Y=-L/2, \ Z=\pm R$.}
    \label{fig:pinched_cylinder_config}
\end{figure}

\begin{figure*}[!hbpt]
    \begin{minipage}{.25\columnwidth}
    \centering
    \subfloat{\includegraphics[width=\columnwidth]{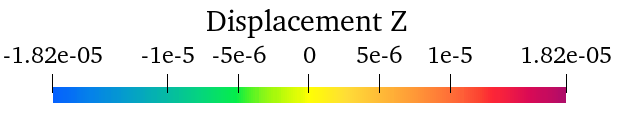}}
    
    \vspace*{-0.25cm}
    \setcounter{subfigure}{0}
    \subfloat[][IGA solution]{\includegraphics[trim=0 0 6.5cm 0, clip=true, width=\columnwidth]{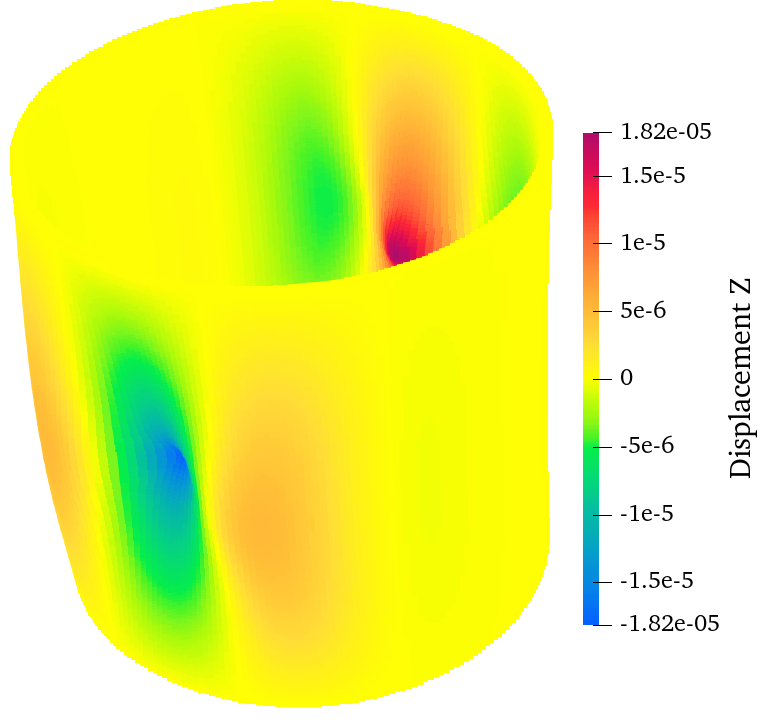}}
    
    \subfloat[][PD solution]{\includegraphics[trim=0 0 6.5cm 0, clip=true,width=\columnwidth]{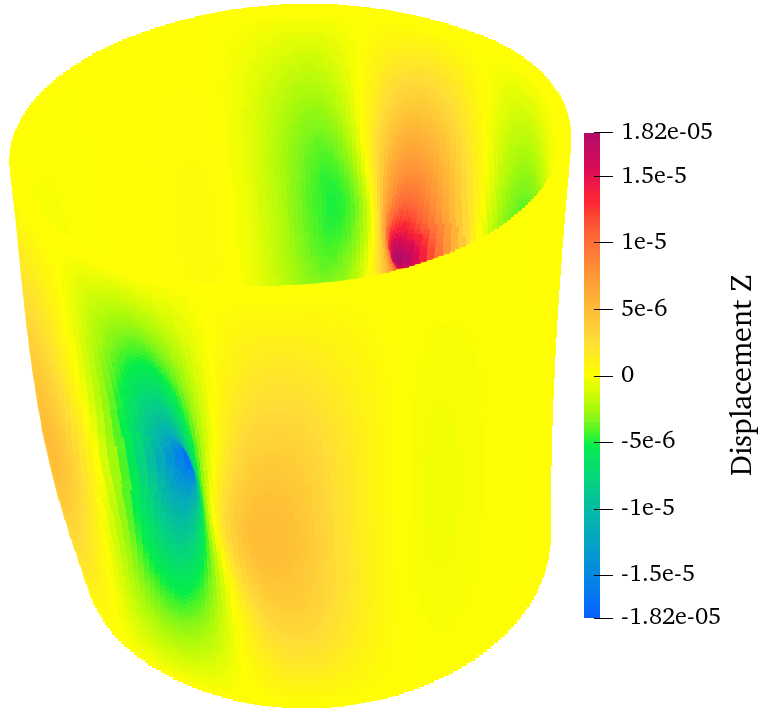}}
    \end{minipage}
    \hfill
    \begin{minipage}{.7\columnwidth}
    \centering
    \subfloat[][Convergence study]{\includegraphics[width=\columnwidth]{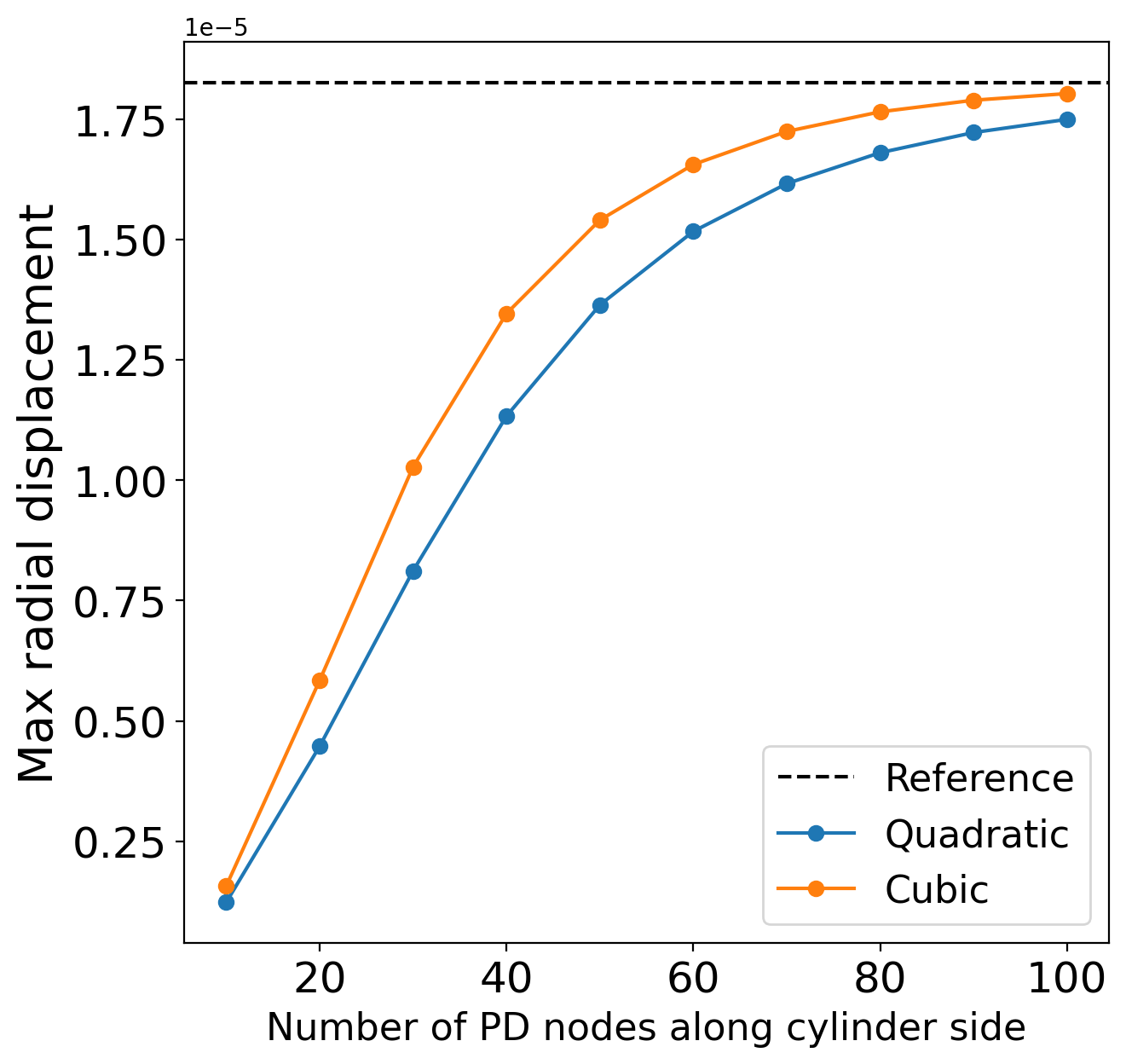}}
    \end{minipage}
    \caption{Pinched cylinder problem. Deformed shape of the cylinder using (a): IGA and (b): PD (displacements are amplified by a factor of $4 \times 10^6$). (c): Convergence behavior of the PD shell using quadratic and cubic RK shape functions.}
    \label{fig:pinched_cylinder}
\end{figure*}

The pinched cylinder is another classical benchmark problem and part of the shell obstacle course~\cite{belytschko1985stress}, which involves a curved geometry with point loads. This example includes a cylindrical shell supported by rigid end diaphragms and subjected to two opposing pinching forces of the same magnitude. This numerical test is demanding due to the severe bending and out-of-plane warping that the cylinder experiences. The problem description is provided in \cref{fig:pinched_cylinder_config}. The radial deflection under the point load of $1.8248 \times 10^{-5}$ is commonly used as a reference value in the literature~\cite{hughes2005isogeometric,kiendl2009isogeometric,coox2017flexible,duong2017new,maurin2018isogeometric}. Good agreement between the PD and IGA solutions is obtained and shown in \cref{fig:pinched_cylinder}(a--b). As shown in \cref{fig:pinched_cylinder}(c), the cubic RK shape functions in the PD discretization achieve faster convergence to the reference solution than their quadratic counterparts.

\subsubsection{Beam Subjected to Pure Bending}
\label{sec:cantilever_beam}

As an example of finite elastic deformation, we consider a beam subjected to pure bending. Using the classical beam theory~\cite{timoshenko1959theory}, the following relation holds for an Euler--Bernouilli beam subjected to a constant bending moment $M$:
\begin{equation}
    M = \frac{EI}{\rho} ,
\end{equation}
where $E$ is the Young's modulus, $I$ is the second moment of area, and $\rho$ is the radius of curvature. Therefore, a clamped beam with length $L$ and subjected to bending moment $\dfrac{EI\theta}{L}$ has the deformed shape of an arc with radius of curvature $\dfrac{L}{\theta}$ and angle $\theta$. 

\begin{figure*}[!hbpt]
    \centering
    \subfloat[][]{\includegraphics[width=0.48\columnwidth]{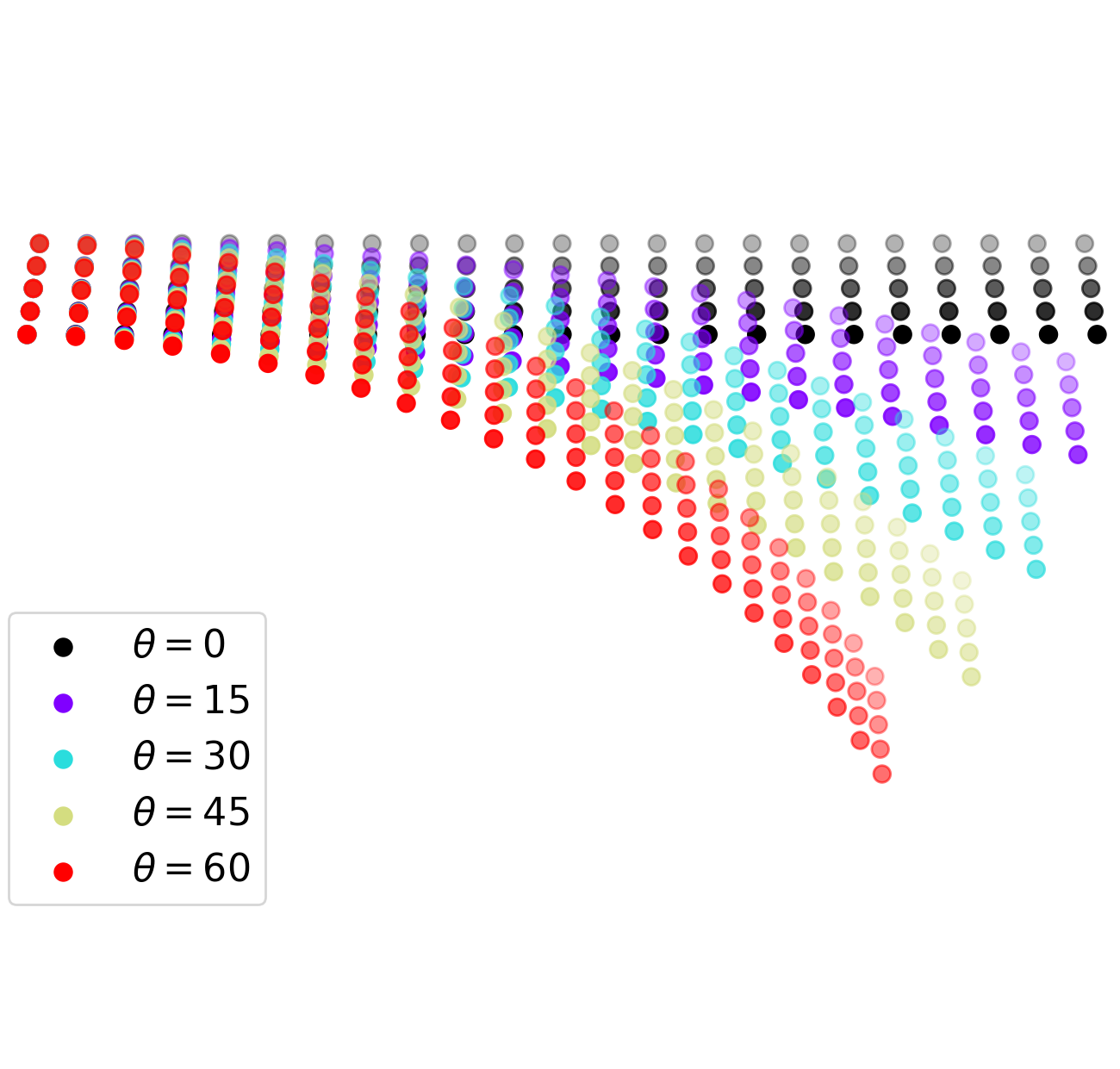}}
    \hfill
    \subfloat[][]{\includegraphics[width=0.48\columnwidth]{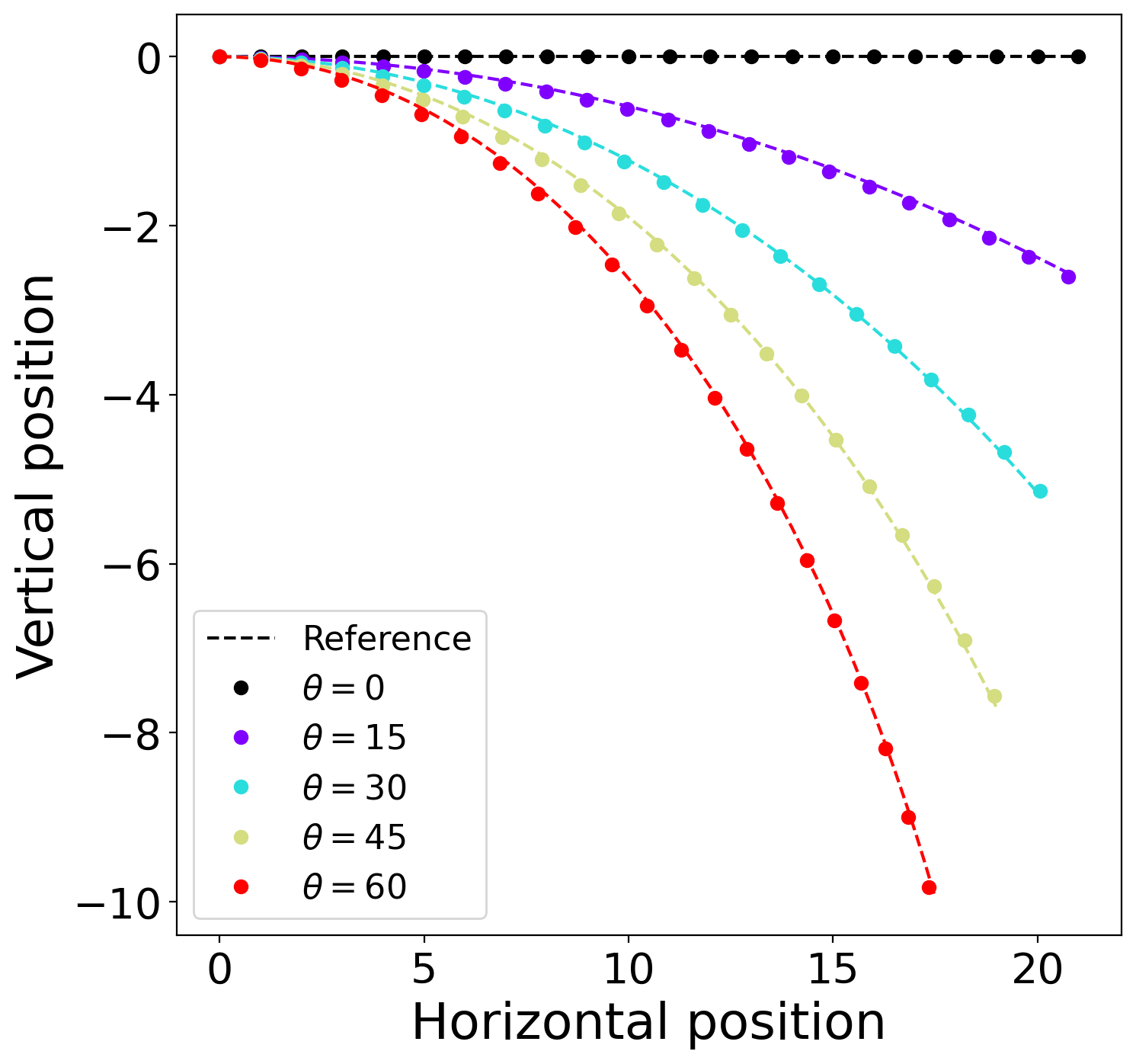}}
    \caption{Beam subjected to pure bending. (a): Deformed configuration of the beam at different loading stages and (b): PD shell deformed configuration superposed on the analytical result. The angle $\theta$ is given in degrees.}
    \label{fig:beam_under_bending}
\end{figure*}

In our setup, the beam is modeled as a narrow plate. The plate is given a Poisson ratio of zero, is clamped on the left, and subjected to pure bending at its right end. The plate has length $L=21$, width $b=5$, and thickness $h=0.25$, and is discretized with $ 21 \times 5$ nodes along $L \times b$. The clamped condition is enforced by adding additional rows of fixed nodes to the left of the plate. Following~\cite{kiendl2009isogeometric}, the bending moment is modeled by a pair of opposing forces that act on the last two rows of PD nodes on the right, and follow the deformation. For a given $\theta$, the force
\begin{equation}
\begin{aligned}
    \mathbf{F} = - & F \left( \sin\theta \ \hat{i} + \cos\theta \ \hat{k} \right) , \\
    & F = \dfrac{EI\theta}{Ld} ,
\end{aligned}
\end{equation}
is applied to the last row of nodes (and $-\mathbf{F}$ is applied to the adjacent row of nodes). $\hat{i}$ and $\hat{k}$ are the coordinate axes along the length and thickness in the undeformed configuration. $d$ is the distance between the last two rows of nodes, and $I=bh^3/12$. The forces are converted to body force per unit area $F_b$ applied uniformly to the nodes in each row, i.e., 
\begin{equation}
    F_b \, A_{\rm row} = F , \qquad A_{\rm row} = b \, \Delta X
\end{equation}
where $\Delta X$ is the spacing between each row. In a series of load increments, $\theta$ is increased from 0 to $\pi/3$ and the problem is solved in a quasi-static manner. As shown in \cref{fig:beam_under_bending}, good agreement between the numerical and analytical results is achieved.

\subsection{Inelastic Problems}
\label{sec:plasticity_problems}

In this section, we study several examples involving plasticity and/or damage phenomena. 

\subsubsection{Necking Tests}
\label{sec:necking}

\begin{figure}[!hbpt]
    \centering
    \includegraphics[width=0.4\columnwidth]{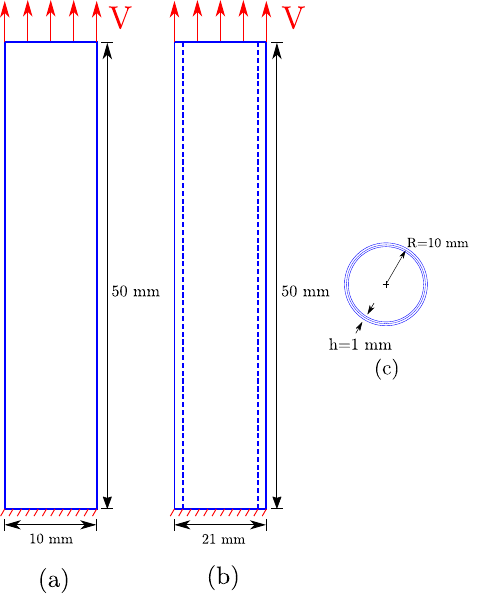}
    \caption{Setup for the necking problems. The constant velocity field V is applied at the top boundary. (a): Necking configuration for the flat specimen. (b): Necking configuration for the hollow cylinder. (c): Cross-section of the hollow cylinder with a mid-surface radius of 10 mm and 1 mm wall thickness.}
    \label{fig:necking_setup}
\end{figure}

Two different elasto-plastic problems involving specimens undergoing necking are considered from the work of~\cite{ambati2018isogeometric}. The tests include a flat specimen and a hollow cylinder under tensile loading (displacement-controlled boundary conditions). The problem geometry and boundary conditions are shown in \cref{fig:necking_setup}. 
Both cases involve an isotropic hardening material described by the following properties given numerically: Young's modulus E $= 189 \times 10^{3}$, Poisson's ratio $\nu = 0.3$, and a nonlinear hardening law $\sigma_Y(\bar{\epsilon}^P) = 343 + (680 - 343) \, (1 - \exp(-16.93 \,  \bar{\epsilon}^P)) + 300 \, \bar{\epsilon}^P$. A displacement-vector norm ${\rm U_{ Norm}}$ is calculated for each velocity increment to compare with the IGA KL shell solutions using the total Lagrangian approach~\cite{ambati2018isogeometric} and the updated Lagrangian approach~\cite{alaydin2021updated}. The displacement-vector norm ${\rm U_{ Norm}}$ is defined as
\begin{equation}
    {\rm U_{ Norm}} = \sqrt{\left[ \sum_{\rm P \in \mathcal{B}^{\rm 2D}}{\mathbf{u}_{\rm (P)} \cdot \mathbf{u}_{\rm (P)}} \right] \bigg / {\rm N_P}} ,
    \label{eqn:Unorm}
\end{equation}
where $\mathbf{u}_{\rm (P)}$ is the displacement vector of P, and ${\rm N_P}$ is the total number of nodes in the body $\mathcal{B}^{\rm 2D}$. In each case, the PD solutions are compared with their IGA counterparts. We simulate these tests using an explicit solver with a moderate loading rate.

\begin{figure*}[!hbpt]
    \centering
    \subfloat{\includegraphics[width=0.75\columnwidth]{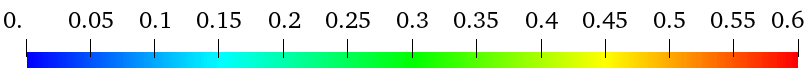}}
    
    \setcounter{subfigure}{0}
    \subfloat[][]{\includegraphics[width=0.15\columnwidth]{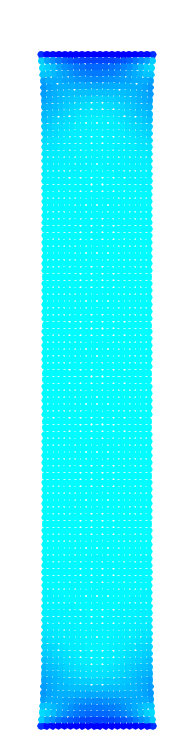}}
    \hspace*{0.35cm}
    \subfloat[][]{\includegraphics[width=0.15\columnwidth]{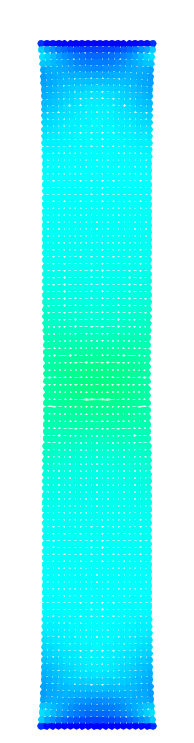}}
    \hspace*{0.35cm}
    \subfloat[][]{\includegraphics[width=0.15\columnwidth]{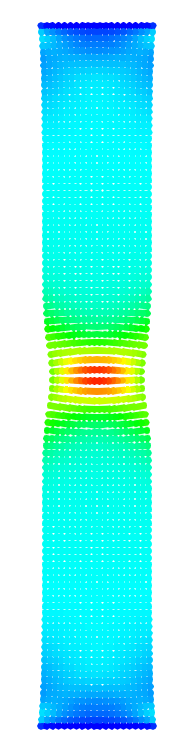}}
    \hspace*{0.35cm}
    \subfloat[][]{\includegraphics[width=0.15\columnwidth]{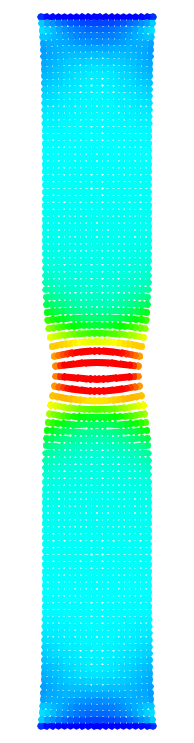}}
    \caption{Necking of a rectangular strip. Equivalent plastic strain ($\bar{\epsilon}^P$) distribution at different loading stages. (a)--(d): ${\rm U_{ Norm}}=$ 4.4, 5, 6, 6.5 mm. Good agreement between these snapshots and those reported in~\cite{ambati2018isogeometric,alaydin2021updated} is observed.}
    \label{fig:necking_plate_eqps}
\end{figure*}

\begin{figure}[!hbpt]
    \centering
    \includegraphics[width=0.9\columnwidth]{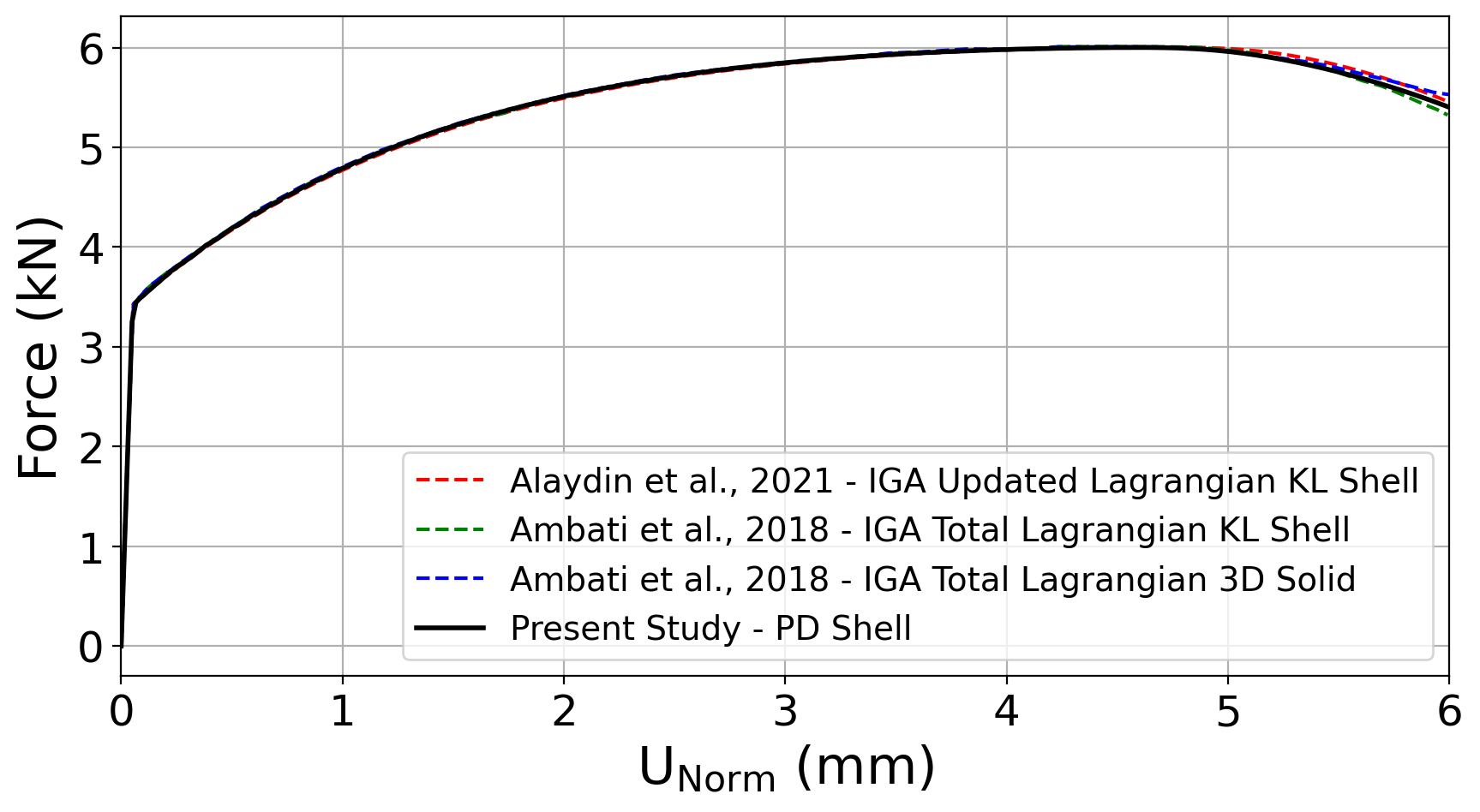}
    \caption{Necking of a rectangular strip. Load-displacement results for the proposed PD shell formulation compared with results reported in~\cite{ambati2018isogeometric,alaydin2021updated}.}
    \label{fig:necking_plate_load}
\end{figure}

In the first necking problem, we study a rectangular strip subjected to uniaxial tension. The specimen is dicretized with $20 \times 50$ PD nodes. The bottom row of nodes is fixed in all directions while the top row of nodes is assigned a constant velocity in the vertical direction and zero velocity in the other two directions. The evolution of equivalent plastic strains is shown in \cref{fig:necking_plate_eqps}. 
As expected, plastic deformation starts from the free edges and plastic instability (necking) develops in the middle section where the stress surpasses the yield strength. \cref{fig:necking_plate_load} provides the load-displacement curve of this problem, involving an elastic regime followed by work hardening and post necking. There is good agreement between the results from our proposed PD shell formulation and the results from~\cite{ambati2018isogeometric,alaydin2021updated}. 

\begin{figure*}[!hbpt]
    \centering
    \subfloat{\includegraphics[width=0.9\columnwidth]{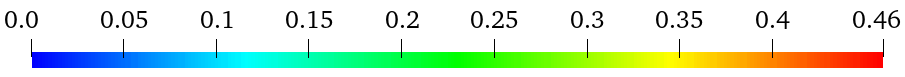}}
    
    \setcounter{subfigure}{0}
    \subfloat[][]{\includegraphics[width=0.32\columnwidth]{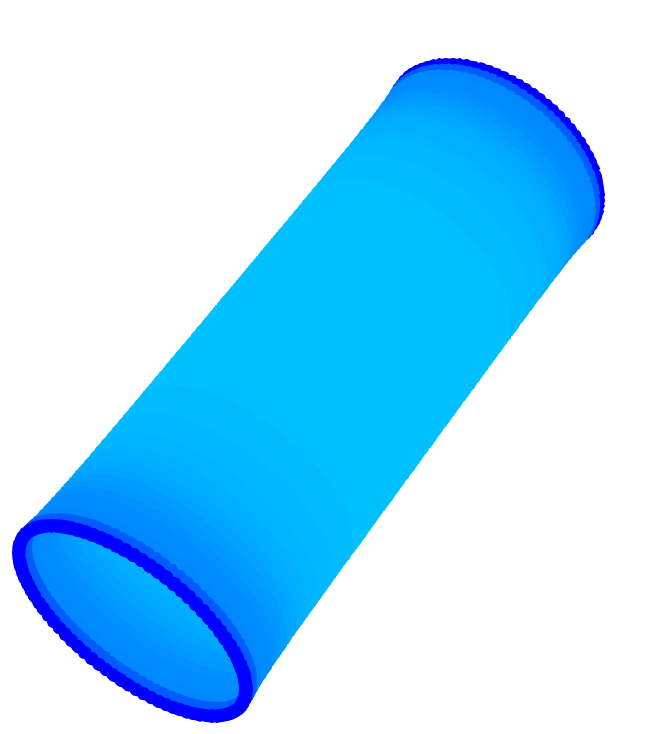}}
    \subfloat[][]{\includegraphics[width=0.32\columnwidth]{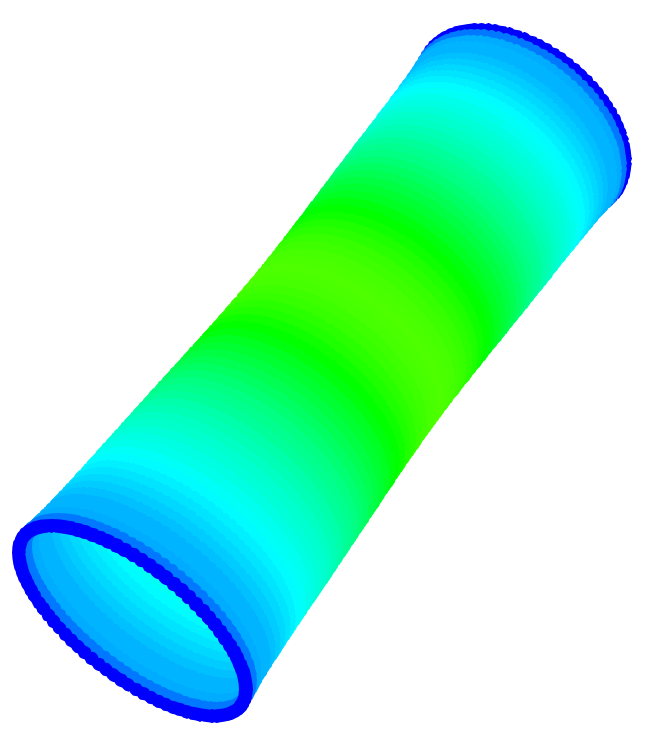}}
    \subfloat[][]{\includegraphics[width=0.32\columnwidth]{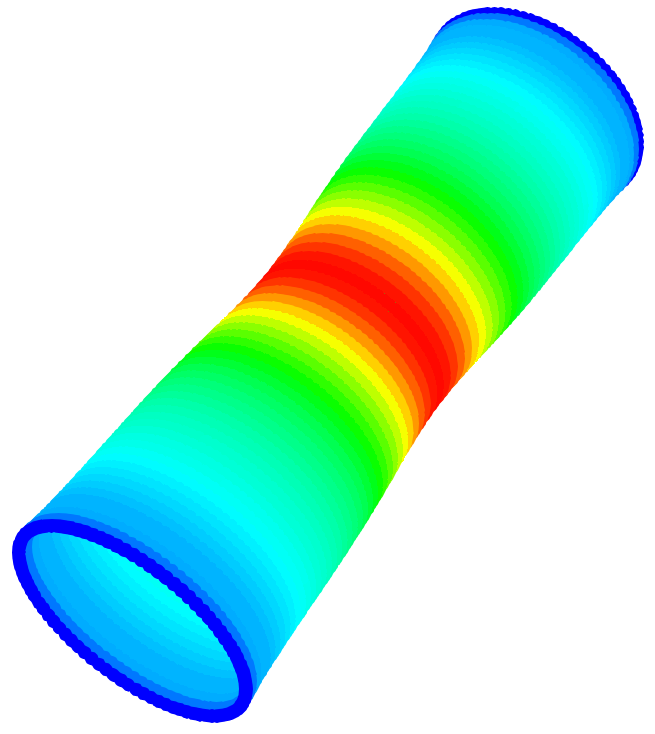}}
    \caption{Necking of a hollow cylinder. Equivalent plastic strain ($\bar{\epsilon}^P$) distribution at different loading stages. (a)--(c): ${\rm U_{ Norm}}=$ 2.5, 5.5, 7 mm. Good agreement between these snapshots and those reported in~\cite{ambati2018isogeometric,alaydin2021updated} is observed.}
    \label{fig:necking_cylinder_eqps}
\end{figure*}

\begin{figure}[!hbpt]
    \centering
    \includegraphics[width=0.9\columnwidth]{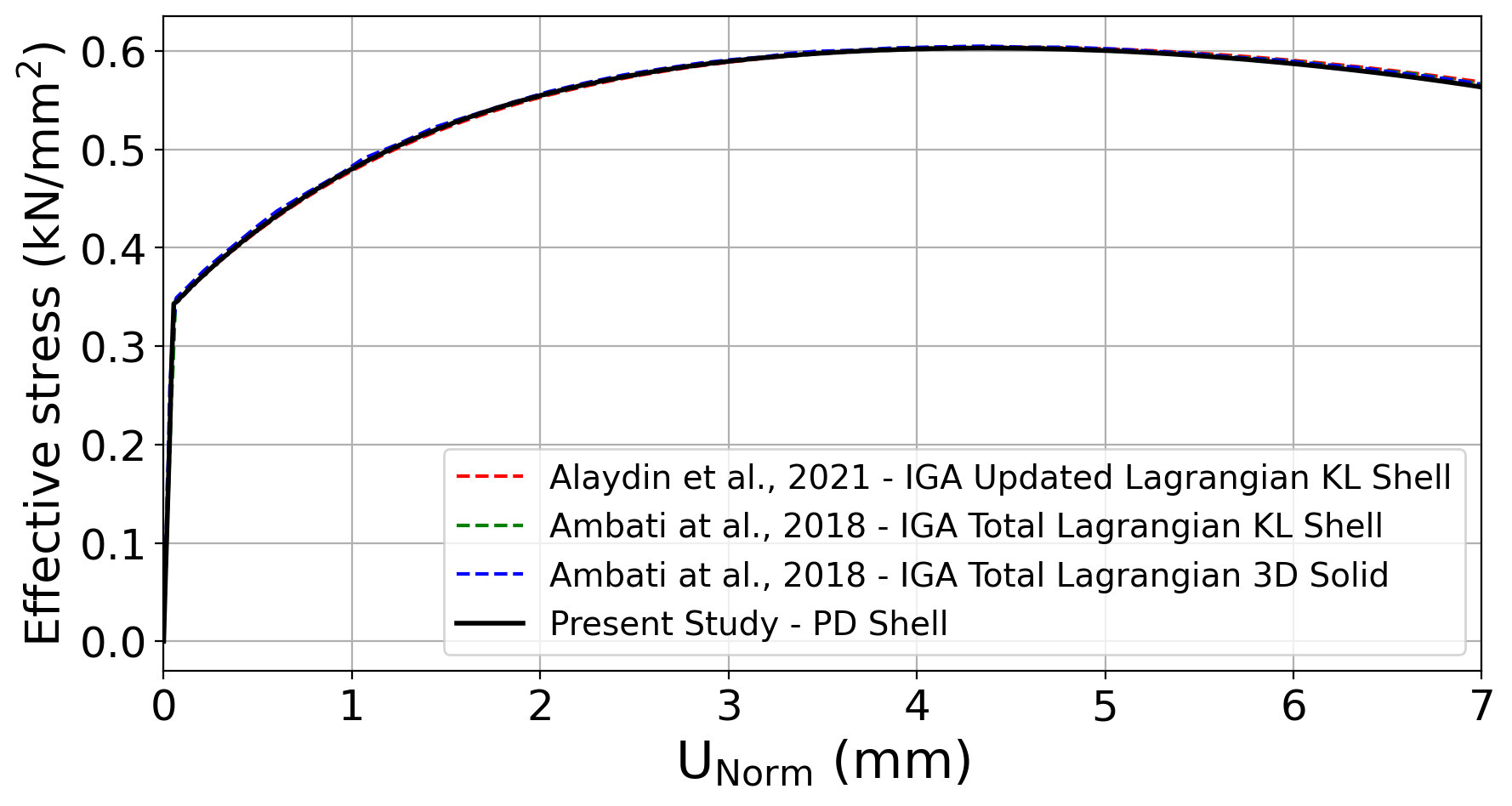}
    \caption{Necking of a hollow cylinder. Load-displacement results for the proposed PD shell formulation compared with results reported in~\cite{ambati2018isogeometric,alaydin2021updated}.}
    \label{fig:necking_cylinder_load}
\end{figure}

The second necking problem considers a cylindrical hollow tube subjected to uniaxial tension. The boundary conditions are applied in a similar fashion to the previous case, i.e., a fixed row of PD nodes at the bottom and a top row of nodes with prescribed velocity boundary conditions in the axial direction and fixed in the other two directions. The specimen is discretized with an average nodal spacing of 0.667 mm ($75 \times 90$ nodes). The distribution of equivalent plastic strain on the deformed geometry is shown in \cref{fig:necking_cylinder_eqps}, indicating typical characteristics of necking phenomena in hollow tubes. The load-displacement results of the proposed PD shell model is compared with the results of~\cite{ambati2018isogeometric,alaydin2021updated} in \cref{fig:necking_cylinder_load}. Good agreement between our results and prior work is obtained. 

\subsubsection{Twisted Hollow Cylinder}
\label{sec:twist}

\begin{figure*}[!hbpt]
    \centering
    \includegraphics[width=\columnwidth]{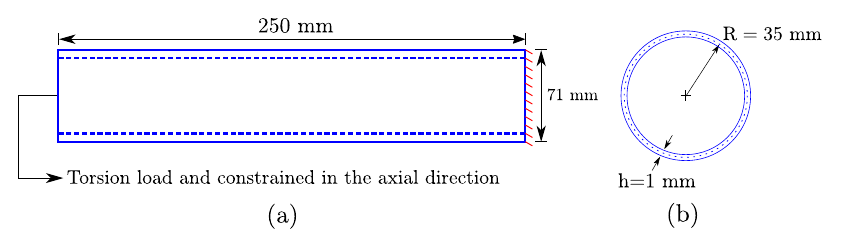}
    \caption{The setup for the twisted cylinder problem. (a): Problem geometry and BCs. (b): Cross-section of the hollow cylinder with a mid-surface radius of 35 mm and 1 mm wall thickness.}
    \label{fig:twisted_configuration}
\end{figure*}

In this example, simulation of a large elasto-plastic twist deformation of a cylindrical shell is conducted. The problem geometry and boundary conditions are shown in \cref{fig:twisted_configuration}. The material properties are the same as in \cref{sec:necking}. The specimen is discretized with an average nodal spacing of 2.5 mm ($86 \times 100$ PD nodes). The problem is driven by displacement control applied to one row of nodes on each end of the structure. 

\begin{figure*}[!hbpt]
    \centering
    \subfloat{\includegraphics[width=0.75\columnwidth]{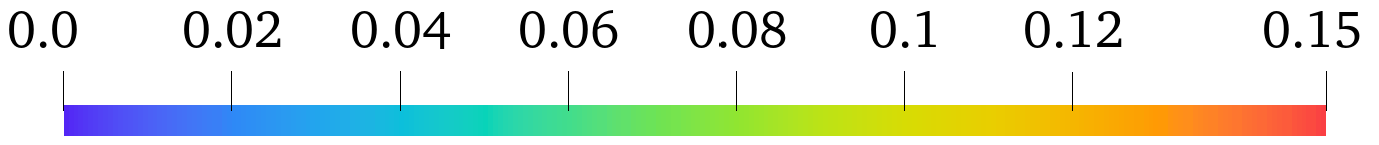}}
    
    \setcounter{subfigure}{0}
    \subfloat[][]{\includegraphics[width=0.12\columnwidth]{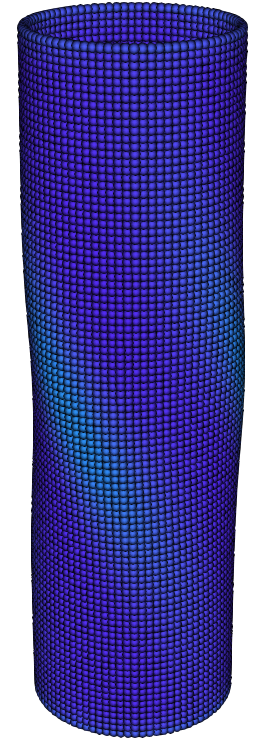}}
    \hspace*{0.3cm}
    \subfloat[][]{\includegraphics[width=0.12\columnwidth]{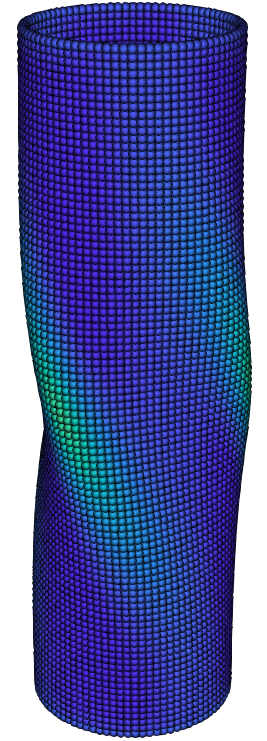}}
    \hspace*{0.3cm}
    \subfloat[][]{\includegraphics[width=0.12\columnwidth]{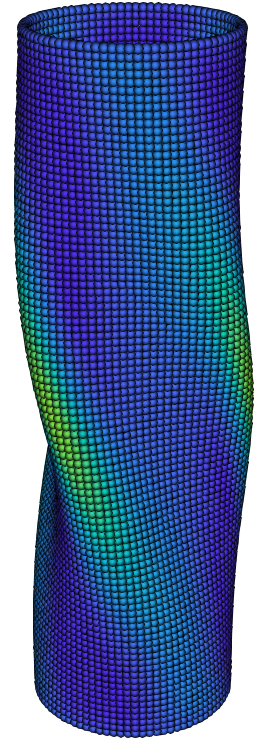}}
    \hspace*{0.3cm}
    \subfloat[][]{\includegraphics[width=0.12\columnwidth]{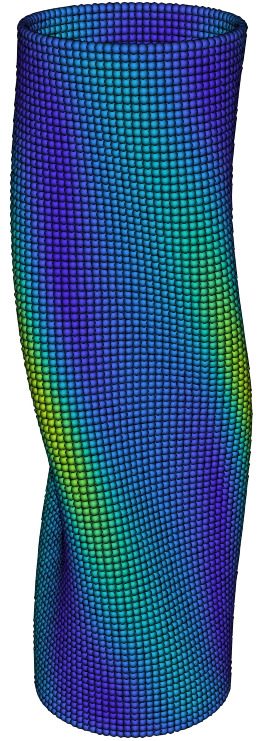}}
    \hspace*{0.3cm}
    \subfloat[][]{\includegraphics[width=0.12\columnwidth]{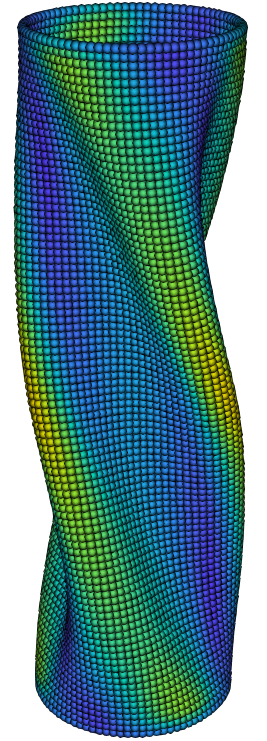}}
    \hspace*{0.3cm}
    \subfloat[][]{\includegraphics[width=0.12\columnwidth]{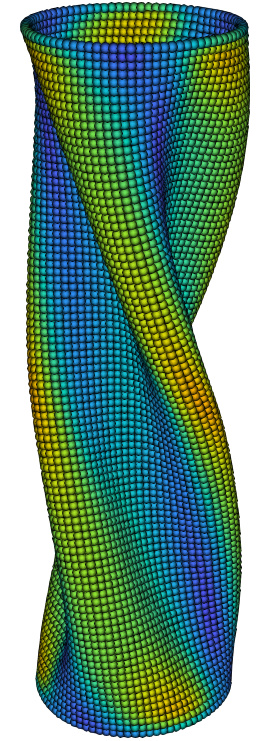}}
    \hspace*{0.3cm}
    \subfloat[][]{\includegraphics[width=0.12\columnwidth]{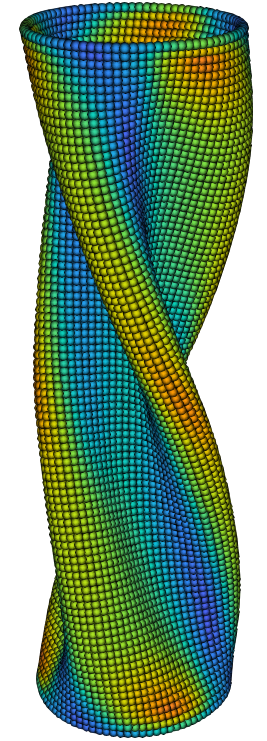}}
    \caption{Large plastic deformation of a twisted cylinder. Equivalent plastic strain ($\bar{\epsilon}^P$) distribution at different loading stages. (a)--(g): Twist angle of $5^{\circ}$, $10^{\circ}$, $20^{\circ}$, $30^{\circ}$, $60^{\circ}$, $90^{\circ}$, and $110^{\circ}$.}
    \label{fig:twist_eqps}
\end{figure*}

\begin{figure}[!hbpt]
    \centering
    \includegraphics[width=0.9\columnwidth]{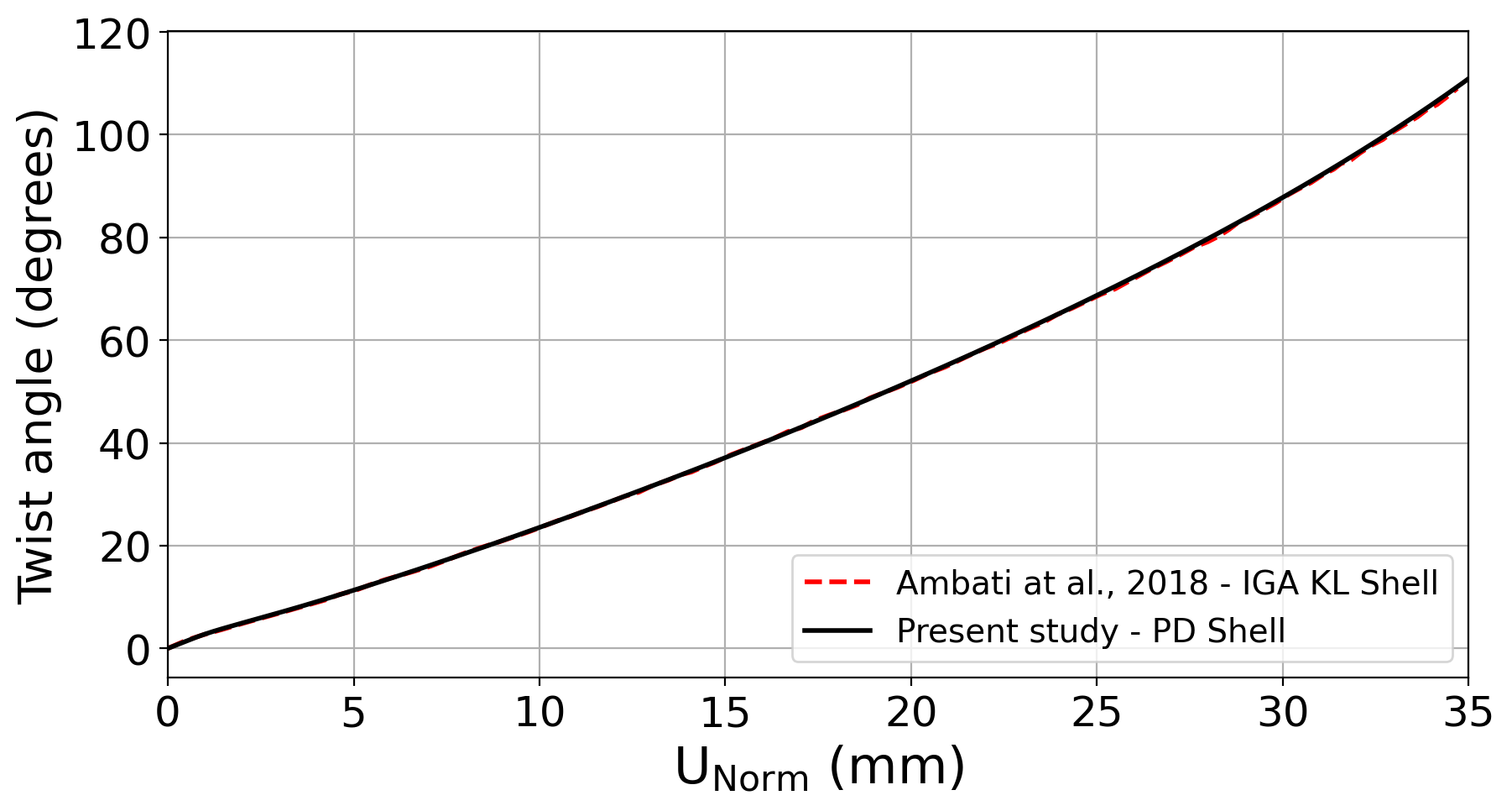}
    \caption{Twisted cylinder under large elasto-plastic deformations. Angle of twist-displacement curve for the proposed PD shell formulation compared with the IGA results reported in~\cite{ambati2018isogeometric}.}
    \label{fig:twist_load}
\end{figure}

To compare our results with those of~\cite{ambati2018isogeometric}, ${\rm U_{ Norm}}$, as defined in \cref{eqn:Unorm}, is computed for each load increment. \cref{fig:twist_eqps} shows a series of snapshots with the equivalent plastic strain plotted on the deformed configuration of the structure at multiple loading stages. 
The angle of twist versus the total displacement norm is plotted in \cref{fig:twist_load}. Good agreement between our results and the reference solution~\cite{ambati2018isogeometric} is obtained. Note that the applied loading in this problem is carried out in such a way that no self-contact occurs~\cite{ambati2018isogeometric}.

\subsubsection{Dynamic Crack Branching}
\label{sec:crack_branching}

\begin{figure*}[!hbpt]
    \centering
    \includegraphics[width=0.7\columnwidth]{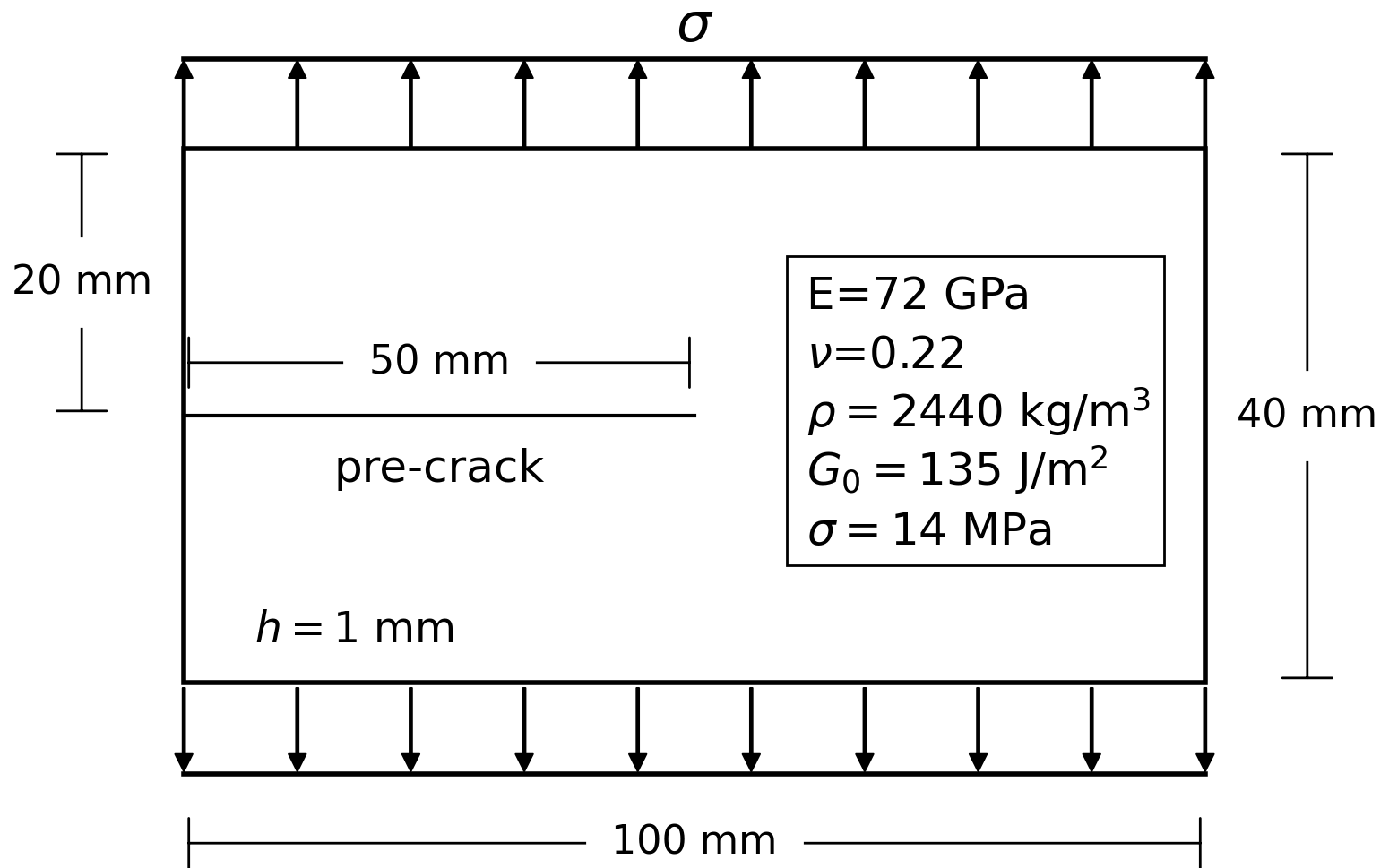}
    \caption{The problem description for the dynamic crack branching example.}
    \label{fig:crack_branching_setup}
\end{figure*}

A plane stress problem involving dynamic crack branching phenomena is studied in this example. As shown in \cref{fig:crack_branching_setup},a pre-notched thin rectangular elastic plate is subjected to impulsive mode-I loading in which the interaction of the reflected stress waves with the propagating fracture causes crack branching~\cite{bowden1967controlled,ravi1984experimental}.
The material (soda-lime glass) properties, geometry, and boundary conditions are chosen as described in~\cite{ha2010studies}. To simulate brittle fracture in this setup, the critical stretch failure criterion~\cite{silling2005meshfree} is utilized, with the following relationship between the critical bond stretch ($s_0$) and the energy release rate ($G_0$)~\cite{ha2010studies}:
\begin{equation}
    s_0 = \sqrt{\frac{5 \pi G_0}{9 E \delta}} .
\end{equation}
To examine the response under mesh refinement, four different discretizations are considered with the nodal spacing of $\Delta x = 1.0, 0.5, 0.25, 0.125$~mm. The horizon size is chosen to be $\delta = 4 \Delta x$. The bonds that cross the pre-cracked surface are broken in the undeformed configuration. The traction tensile loads are applied abruptly along a layer of PD nodes at the top and bottom of the plate and then held fixed for the remainder of the simulations. 

\begin{figure*}[!hbpt]
    \subfloat[][]{\includegraphics[width=0.49\columnwidth]{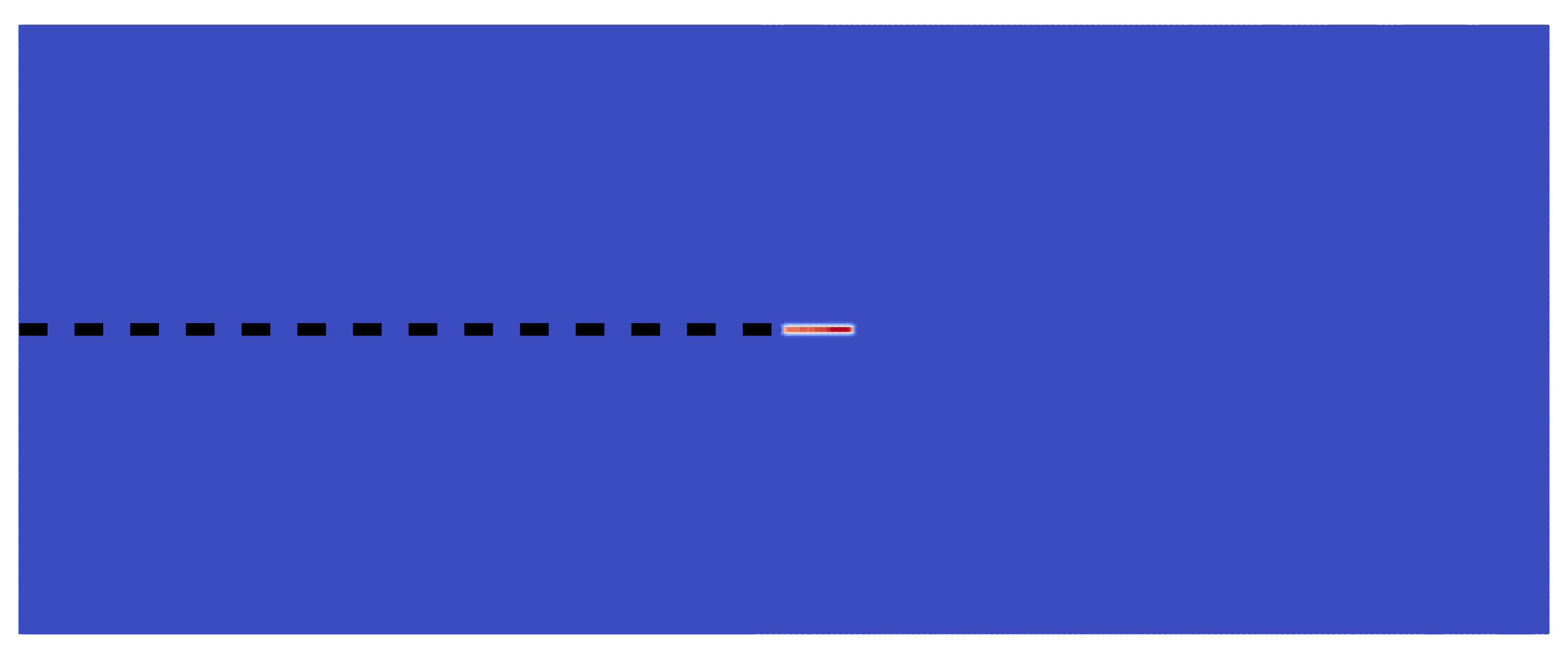}}
    \hfill
    \subfloat[][]{\includegraphics[width=0.49\columnwidth]{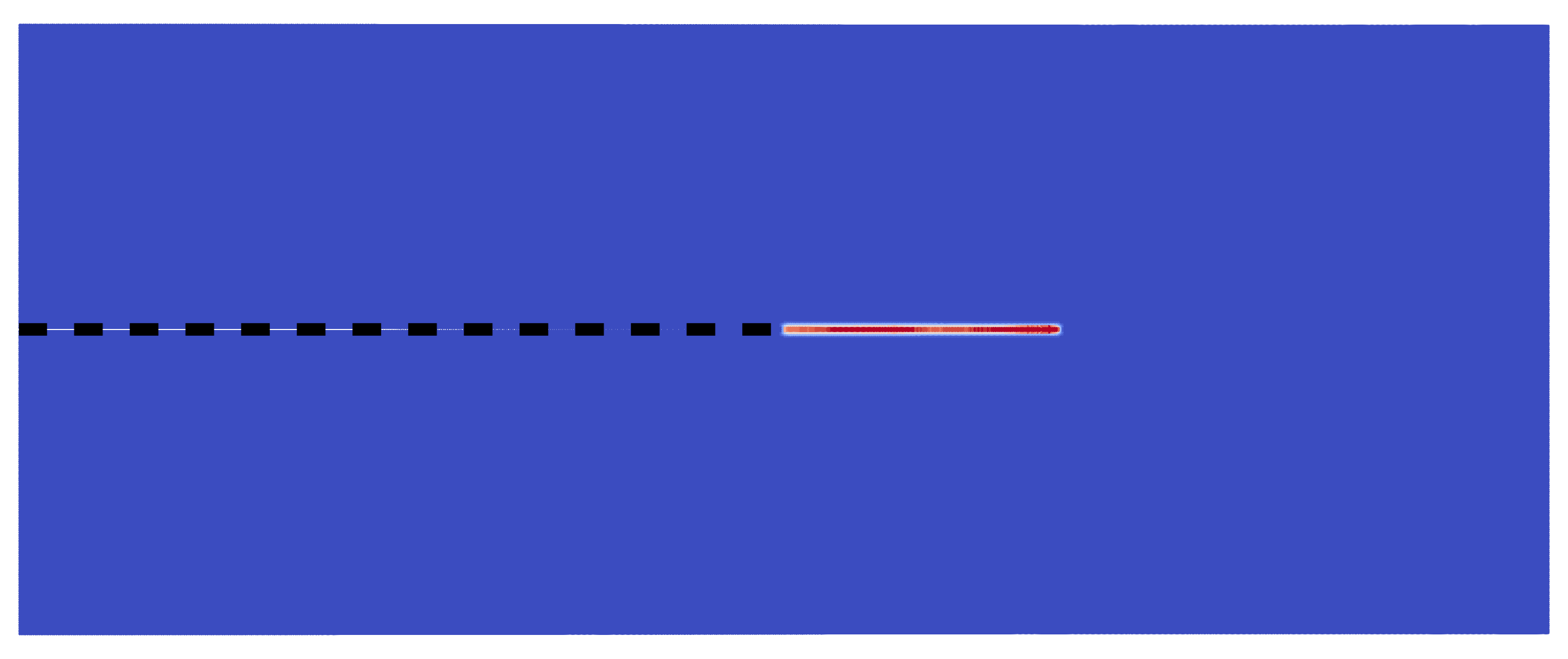}}
    
    \subfloat[][]{\includegraphics[width=0.49\columnwidth]{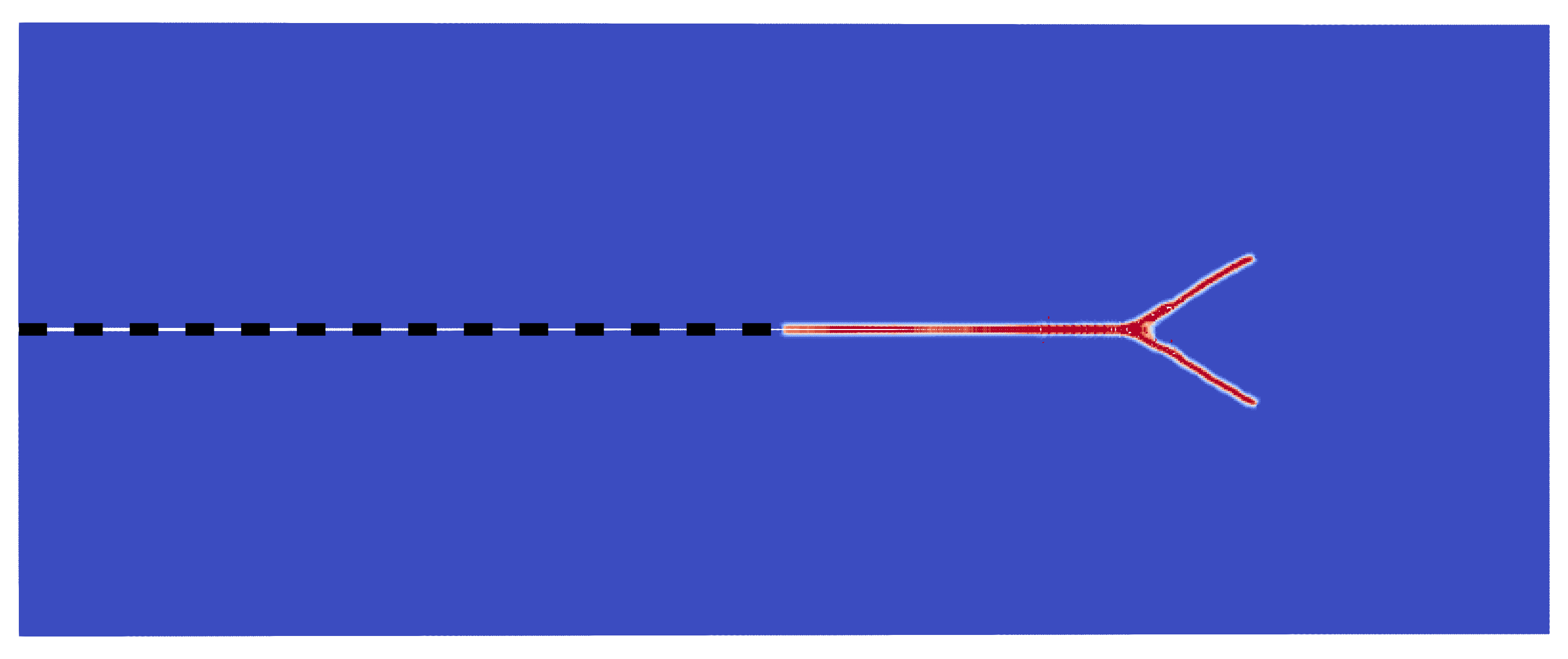}}
    \hfill
    \subfloat[][]{\includegraphics[width=0.49\columnwidth]{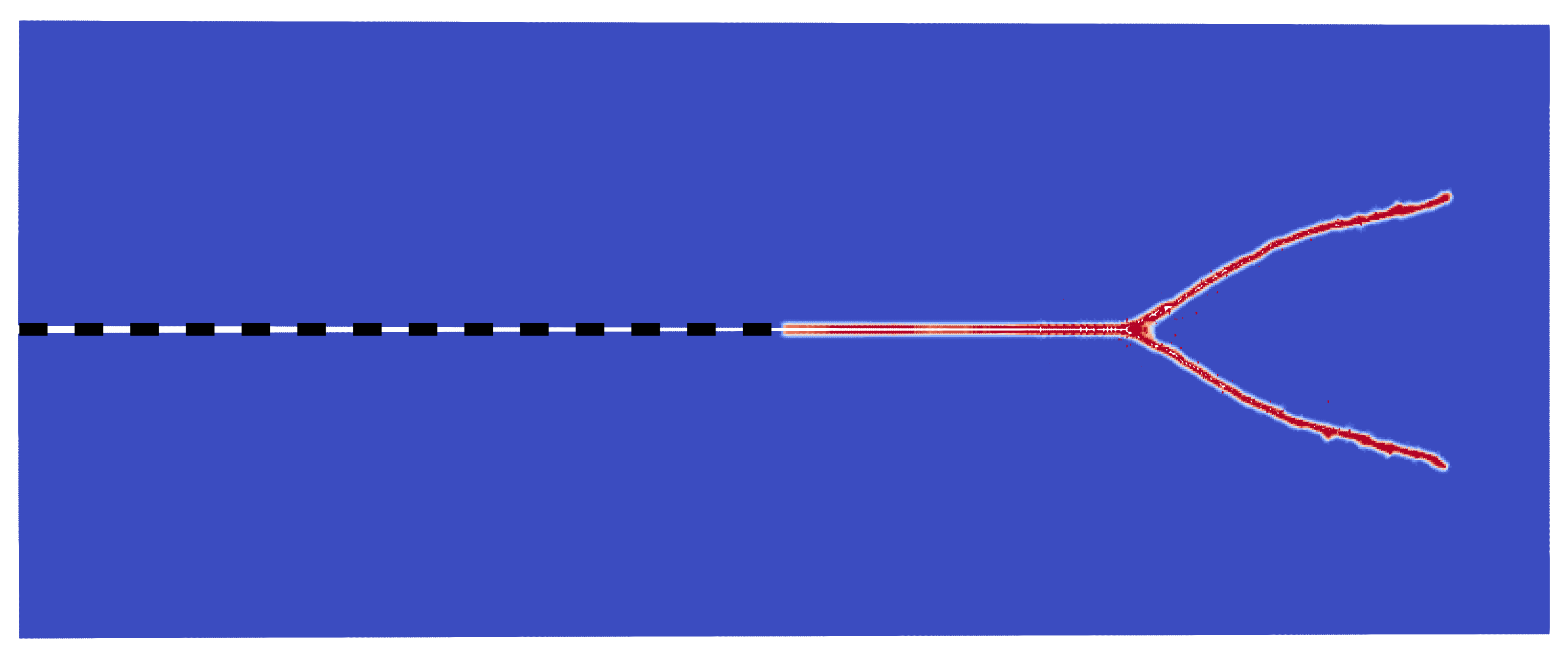}}
    \caption{Dynamic crack branching problem. Crack propagation path for the finest mesh ($\Delta x = \SI{0.125}{\milli\meter}$). (a): $t=\SI{10}{\micro\second}$, (b): $t=\SI{20}{\micro\second}$, (c): $t=\SI{30}{\micro\second}$, (d): $t=\SI{40}{\micro\second}$. Colors indicate the damage index (0 to 0.5).}
    \label{fig:crack_branching_damage_propagation}
\end{figure*}
\begin{figure*}[!hbpt]
    \subfloat[][]{\includegraphics[width=0.49\columnwidth]{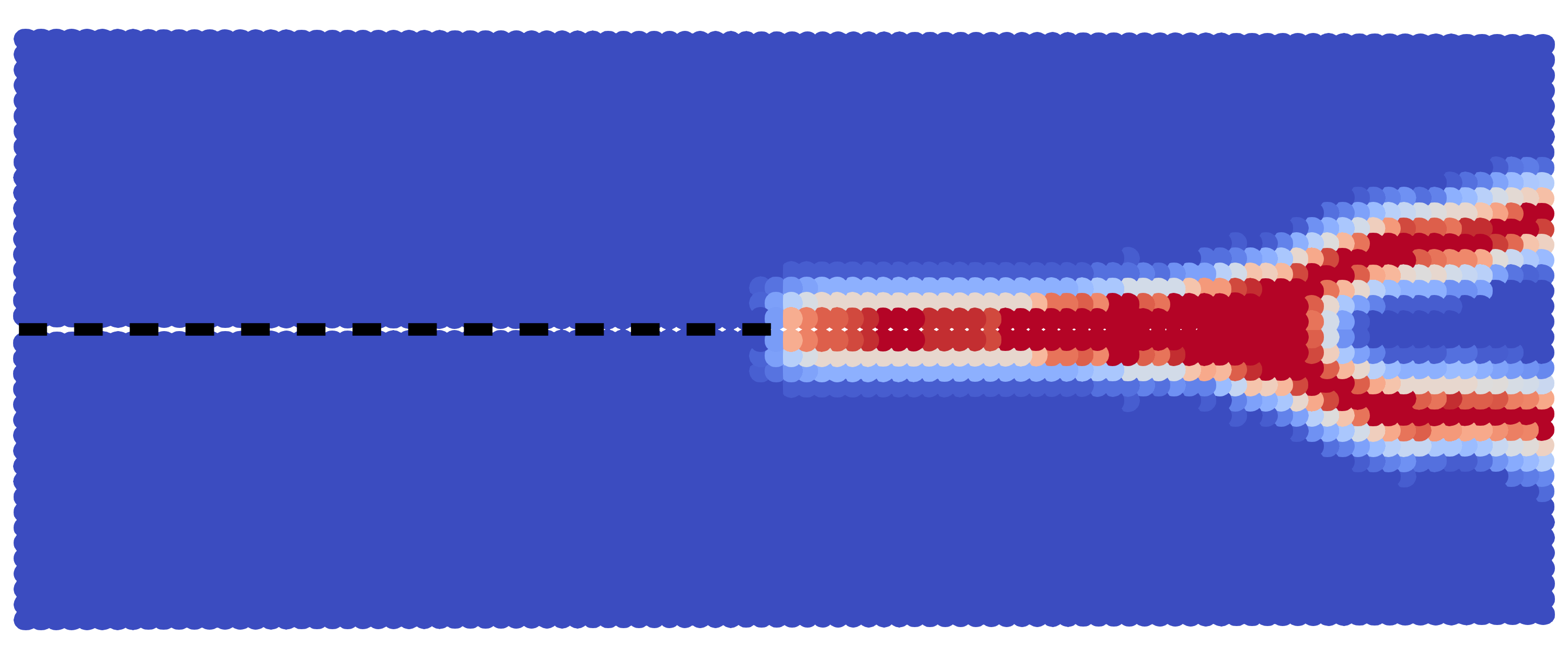}}
    \hfill
    \subfloat[][]{\includegraphics[width=0.49\columnwidth]{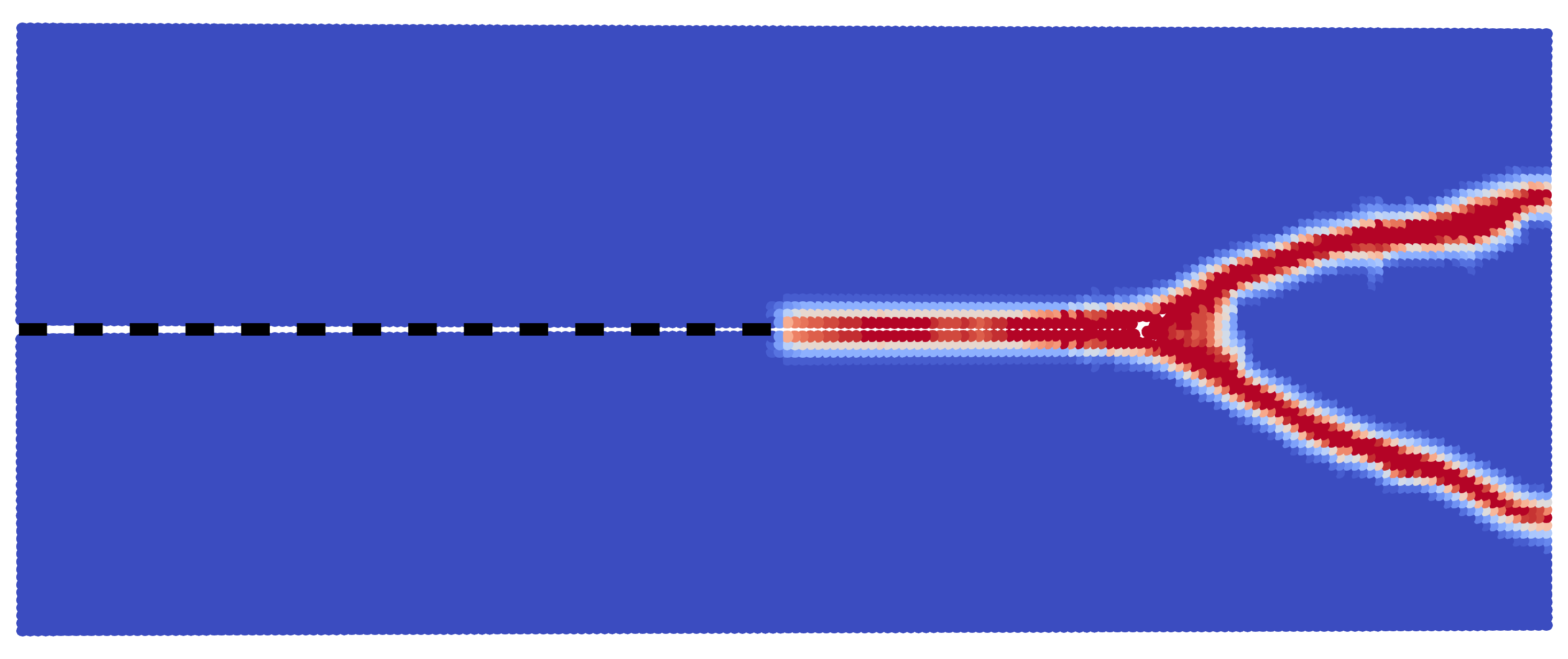}}
    
    \subfloat[][]{\includegraphics[width=0.49\columnwidth]{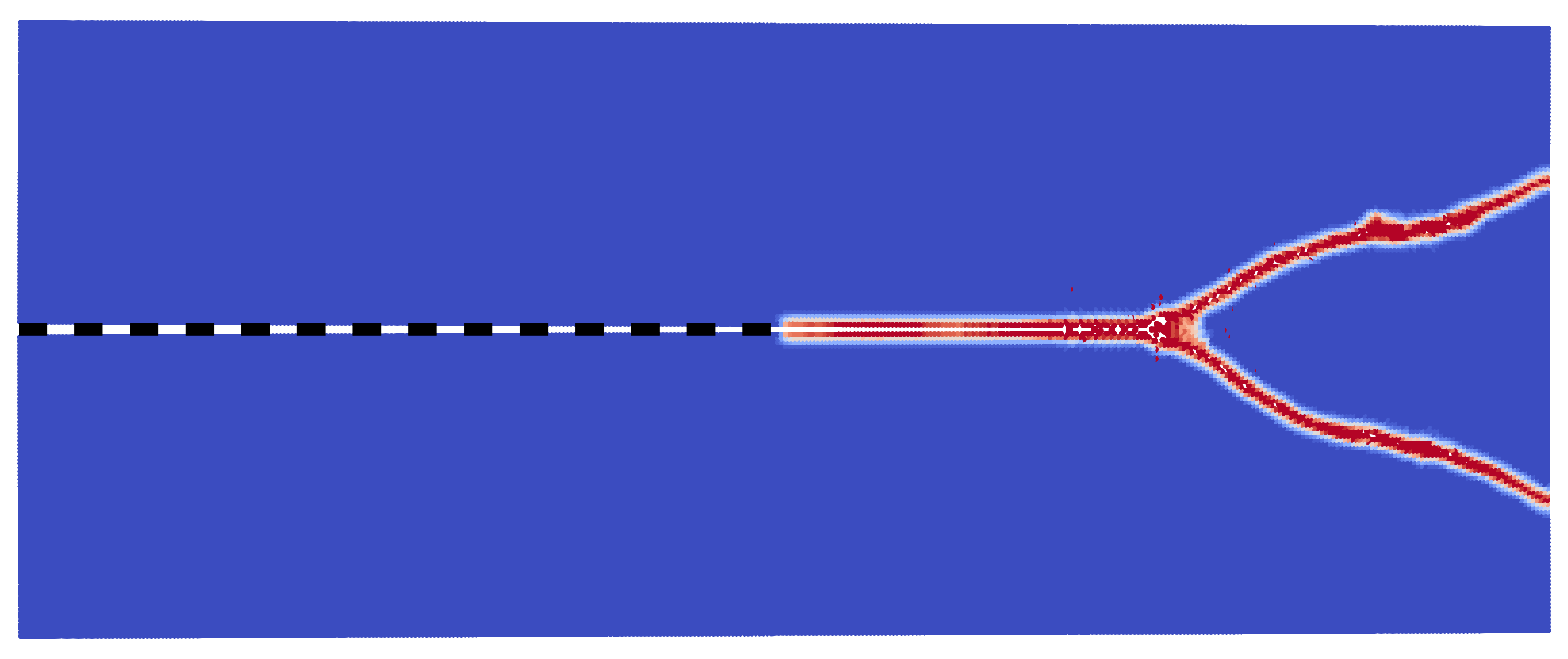}}
    \hfill
    \subfloat[][]{\includegraphics[width=0.49\columnwidth]{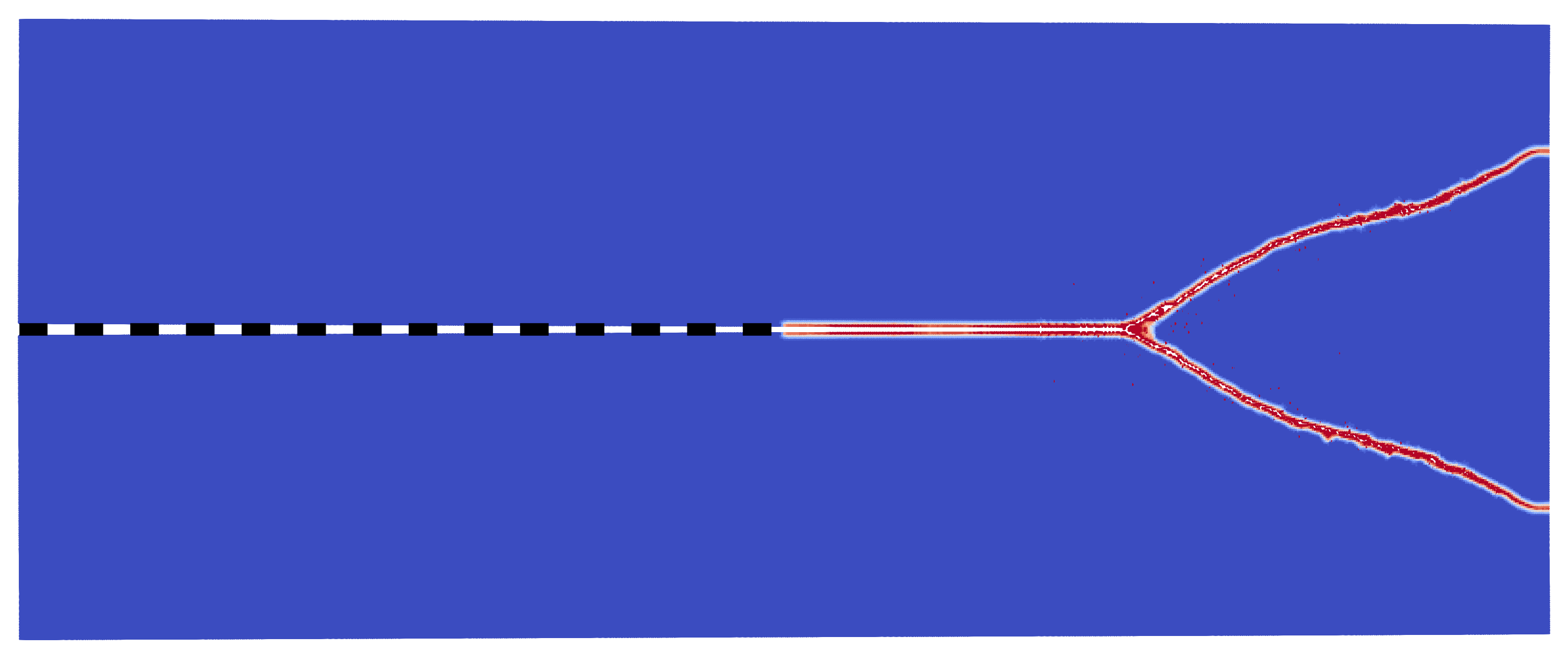}}
    \caption{Dynamic crack branching problem. State of damage at the final configuration ($t=\SI{48}{\micro\second}$) for different particle discretizations. (a): $\Delta x = \SI{1}{\milli\meter}$, (b): $\Delta x = \SI{0.5}{\milli\meter}$, (c): $\Delta x = \SI{0.25}{\milli\meter}$, (d): $\Delta x = \SI{0.125}{\milli\meter}$.}
    \label{fig:crack_branching_final_damage}
\end{figure*}
\begin{figure*}[!hbpt]
    \centering
    \includegraphics[width=0.9\columnwidth]{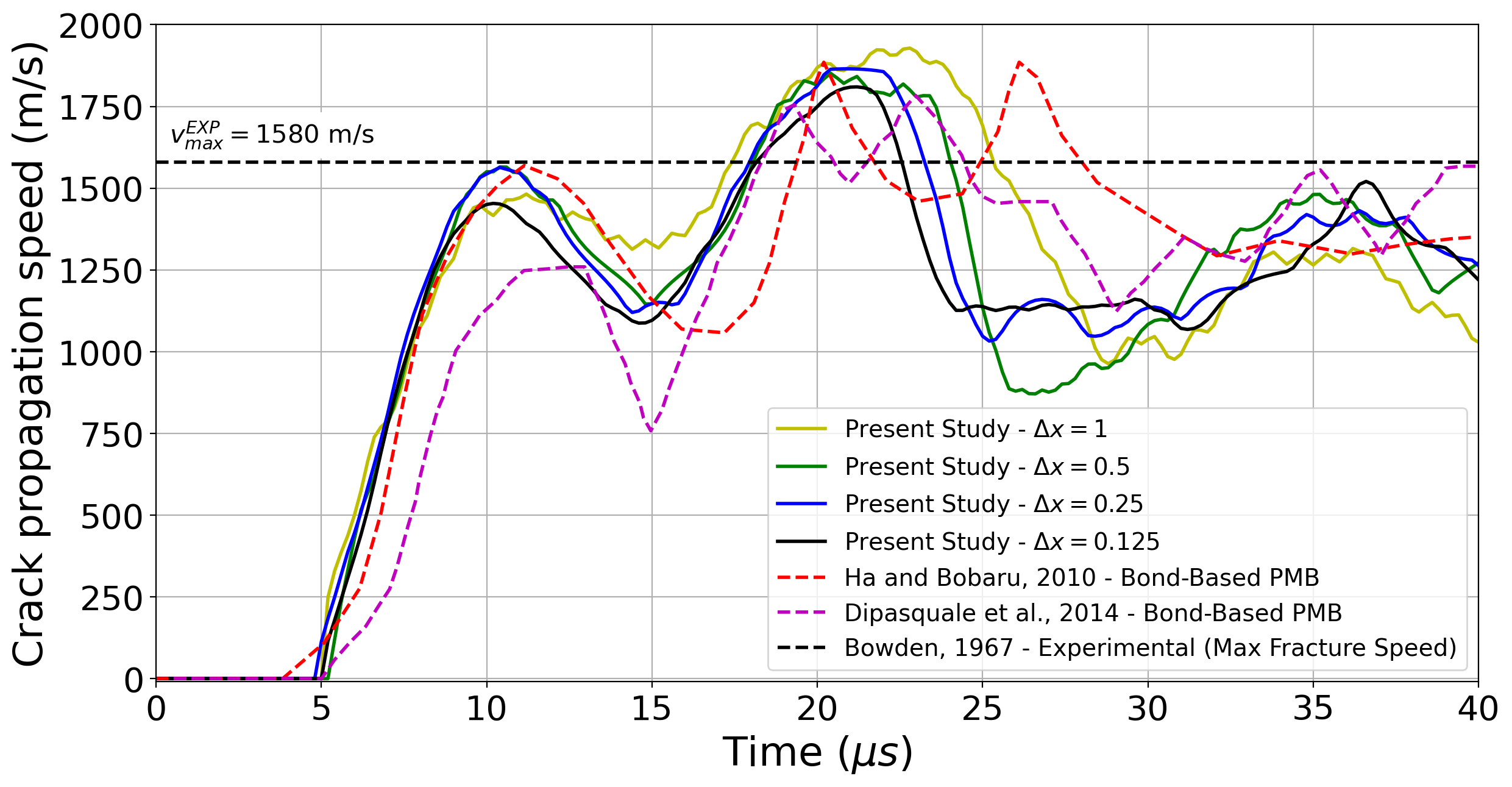}
    \caption{Dynamic crack branching problem. Crack propagation speed for different discretizations of the proposed PD shell model compared with the PMB simulations~\cite{ha2010studies,dipasquale2014crack} and experimental data~\cite{bowden1967controlled}.}
    \label{fig:crack_branching_speed}
\end{figure*}

\cref{fig:crack_branching_damage_propagation} shows the evolution of the crack and its branching on the finest mesh, which is in good agreement with the Prototype Micro Brittle (PMB) simulations of~\cite{ha2010studies} and~\cite{dipasquale2014crack}. The damage state at the final configuration is compared between different discretizations in \cref{fig:crack_branching_final_damage}, where the crack pattern is consistent for all cases. In \cref{fig:crack_branching_speed}, the crack propagation speed is compared between our results, the PMB simulations~\cite{ha2010studies,dipasquale2014crack}, and experimental data~\cite{bowden1967controlled} for maximum crack speed. Good agreement is achieved between our PD results for all the discretizations employed and the reference data. 

\subsubsection{Plasticity-Driven Fracture}
\label{sec:plastic_fracture}

\begin{figure*}[!hbpt]
    \centering
    \includegraphics[width=0.8\columnwidth]{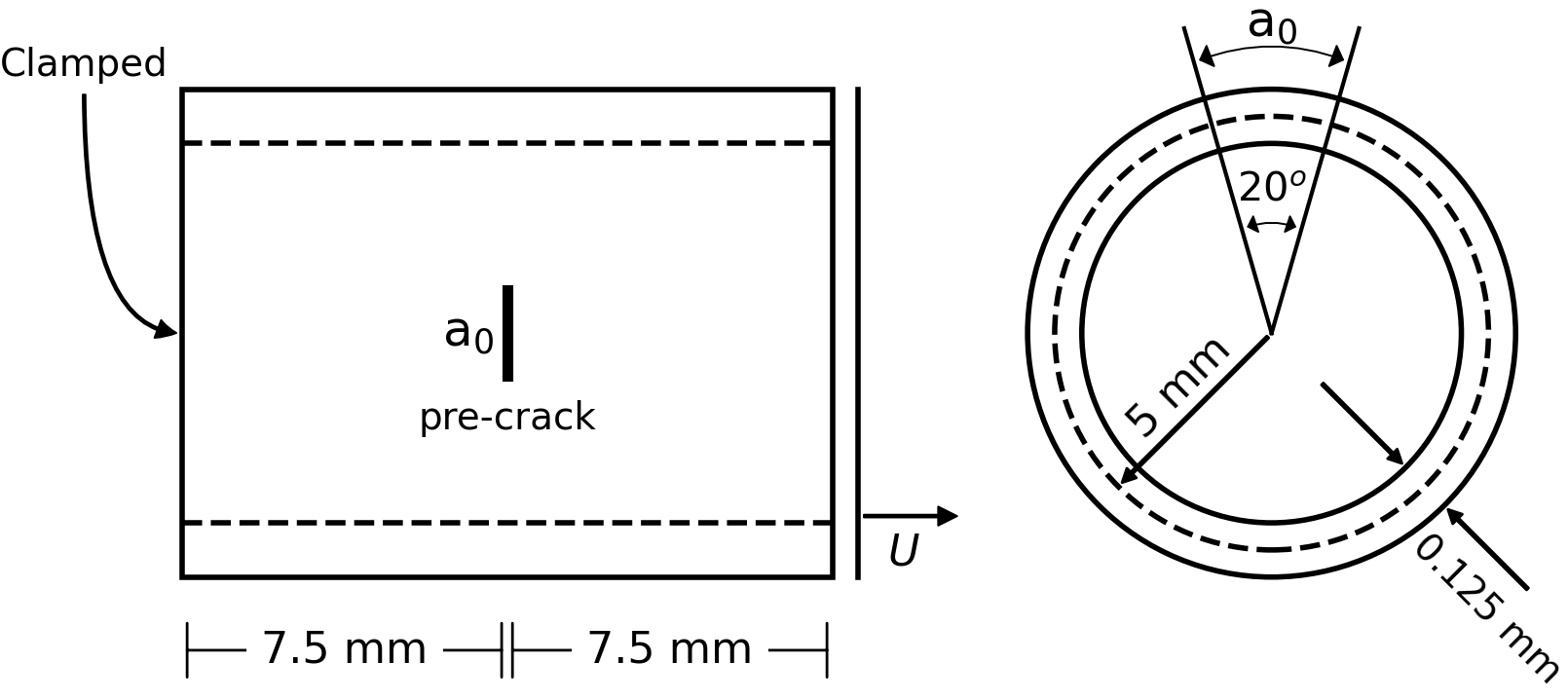}
    \caption{Circumferentially notched cylinder under tensile loading. Description of geometry and boundary conditions.}
    \label{fig:notched_cylinder_setup}
\end{figure*}

A ductile fracture problem is considered in this section. As described in~\cite{ambati2016phase}, an elasto-plastic thin cylinder with a initial circumferential through-thickness notch is subjected to the axial tensile loading. The geometry and boundary conditions are shown in \cref{fig:notched_cylinder_setup}. 
An isotropic hardening material with the following properties is used: Young's modulus E $=72$ GPa, Poisson's ratio $\nu = 0.32$, and a linear hardening law $\sigma_Y(\bar{\epsilon}^P) = 345 + 250 \, \bar{\epsilon}^P$ MPa. The plasticity-driven failure approach is used as in \cref{eqn:plastic_strain} with $\bar{\epsilon}^P_{\rm th}=0.1$ and $\bar{\epsilon}^P_{\rm cr}=0.3$. 

\begin{figure*}[!hbpt]
    \centering
    
    
    
    
    \includegraphics[width=0.8\columnwidth]{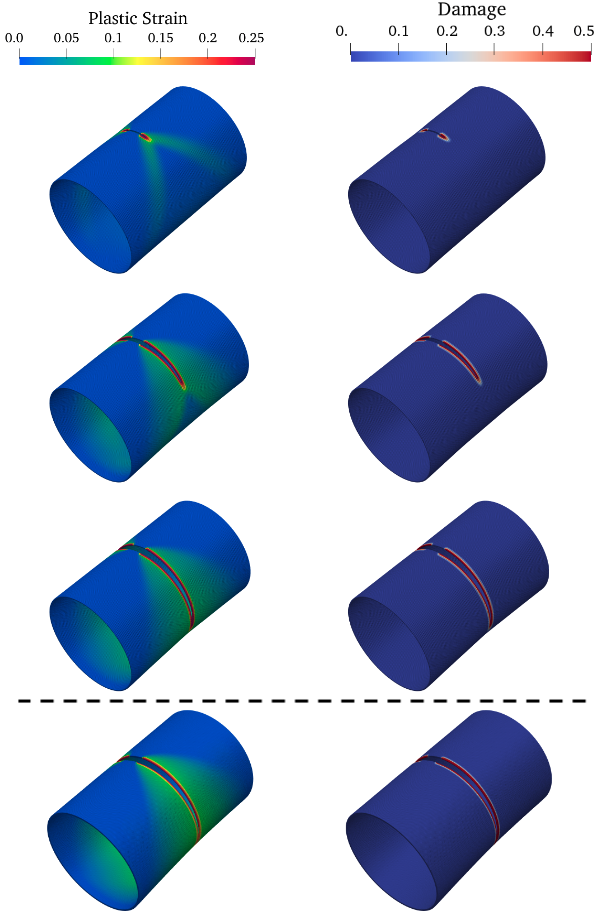}
    \caption{Circumferentially notched cylinder under tensile loading. Distribution of (left): equivalent plastic strain ($\bar{\epsilon}^P$) and (right): damage. Top three rows - PD Shell: the three stages correspond to the loading displacements $U=0.2$, 0.5 and 0.73 mm, respectively. Bottom row - PD Solid: $U=0.73$ mm.}
    \label{fig:notched_cylinder_contours}
\end{figure*}

\begin{figure}[!hbpt]
    \centering
    \includegraphics[width=0.9\columnwidth]{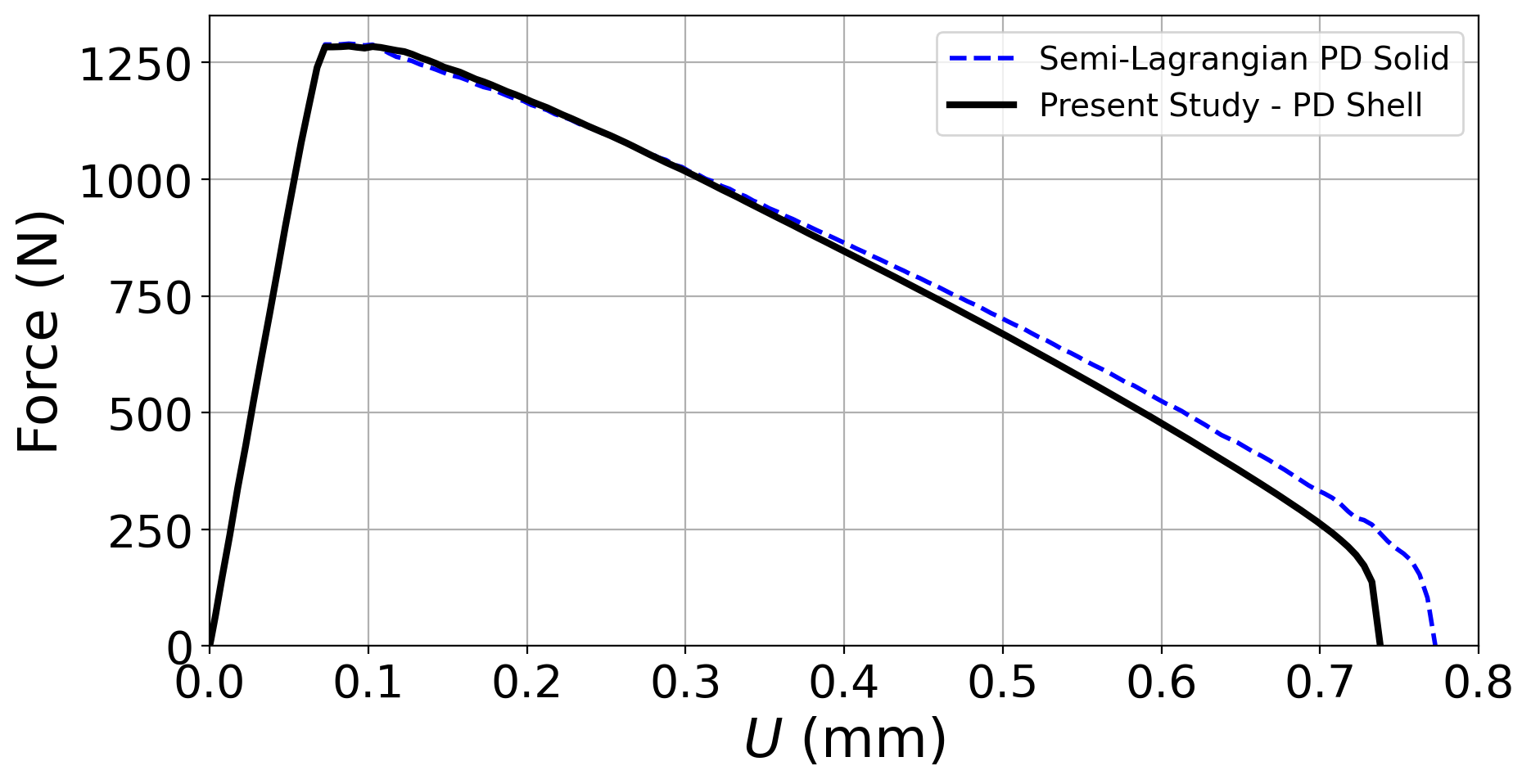}
    \caption{Circumferentially notched cylinder under tensile loading. Load-displacement results for the proposed PD shell formulation are compared with the semi-Lagrangian PD solid simulation.}
    \label{fig:notched_cylinder_load}
\end{figure}

The full cylinder is discretized with an average nodal spacing 0.1~mm, resulting in 48,000 PD nodes. The simulation is conducted using displacement control on the right boundary until failure occurs. The evolution of the plastic strain and damage fields are shown in \cref{fig:notched_cylinder_contours} at several loading stages.
Damage is localized and a stable crack propagation is obtained, which is a common feature of ductile fracture. Our results are qualitatively in good agreement with the simulations of~\cite{ambati2016phase}. Since the force-displacement data in the IGA simulations are not reported in~\cite{ambati2016phase}, we simulate this problem using the semi-Lagrangian, bond-associated PD solid~\cite{behzadinasab2020semi}. The mesh of the solid has five layers of nodes in the thickness direction and the mid-surface layer is coincident with the shell mesh. The macroscopic response of the structure is compared between the shell and solid models in \cref{fig:notched_cylinder_load}, where a good agreement is observed between the two approaches. The simulation time for the solid case is about 15 times longer than the shell case.

\subsubsection{Mode-III Tearing Fracture}
\label{sec:tearing_fracture}

\begin{figure*}[!hbpt]
    \centering
    \subfloat[][]{\includegraphics[width=0.7\columnwidth]{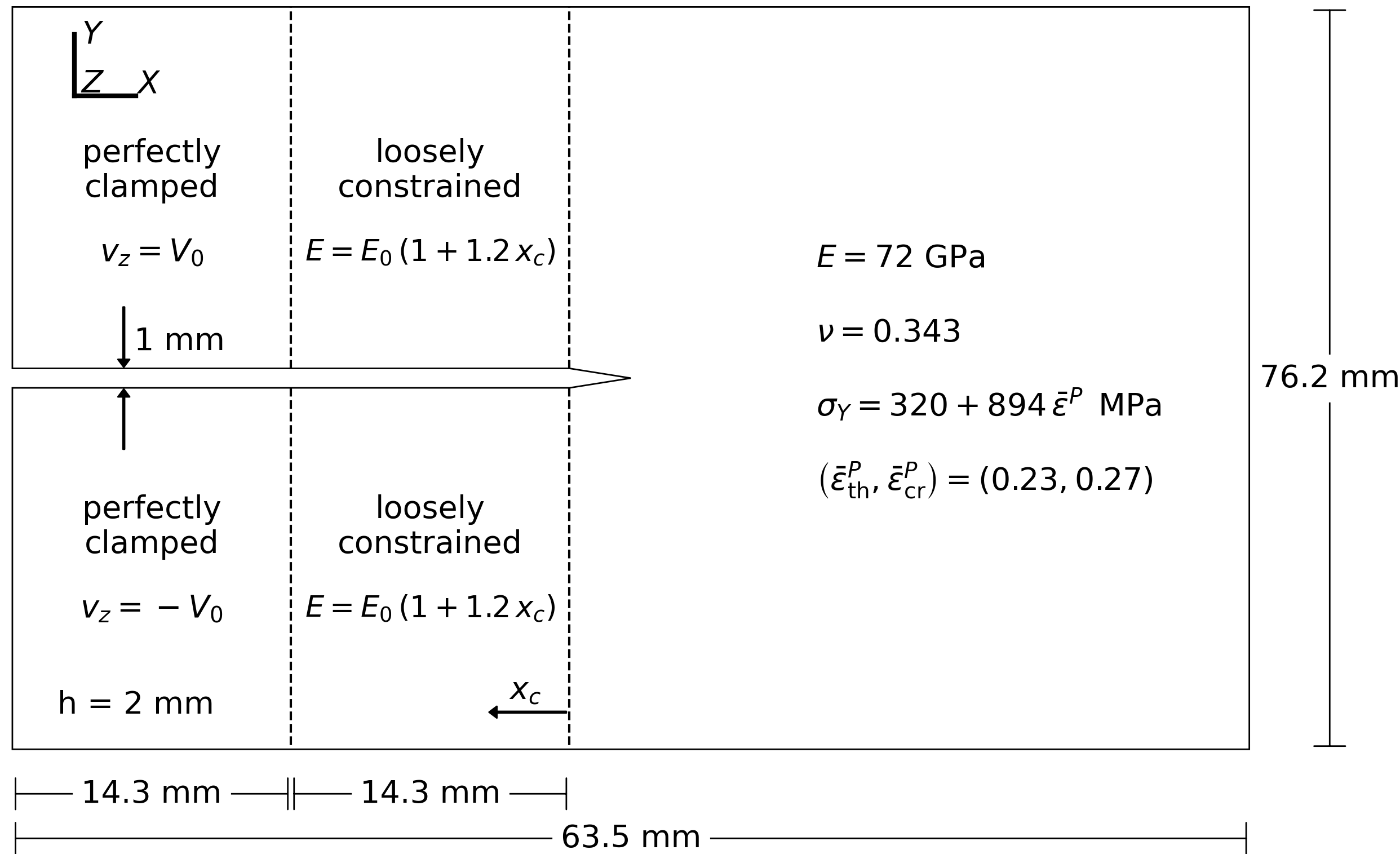}}
    
    \subfloat[][]{\includegraphics[width=0.75\columnwidth]{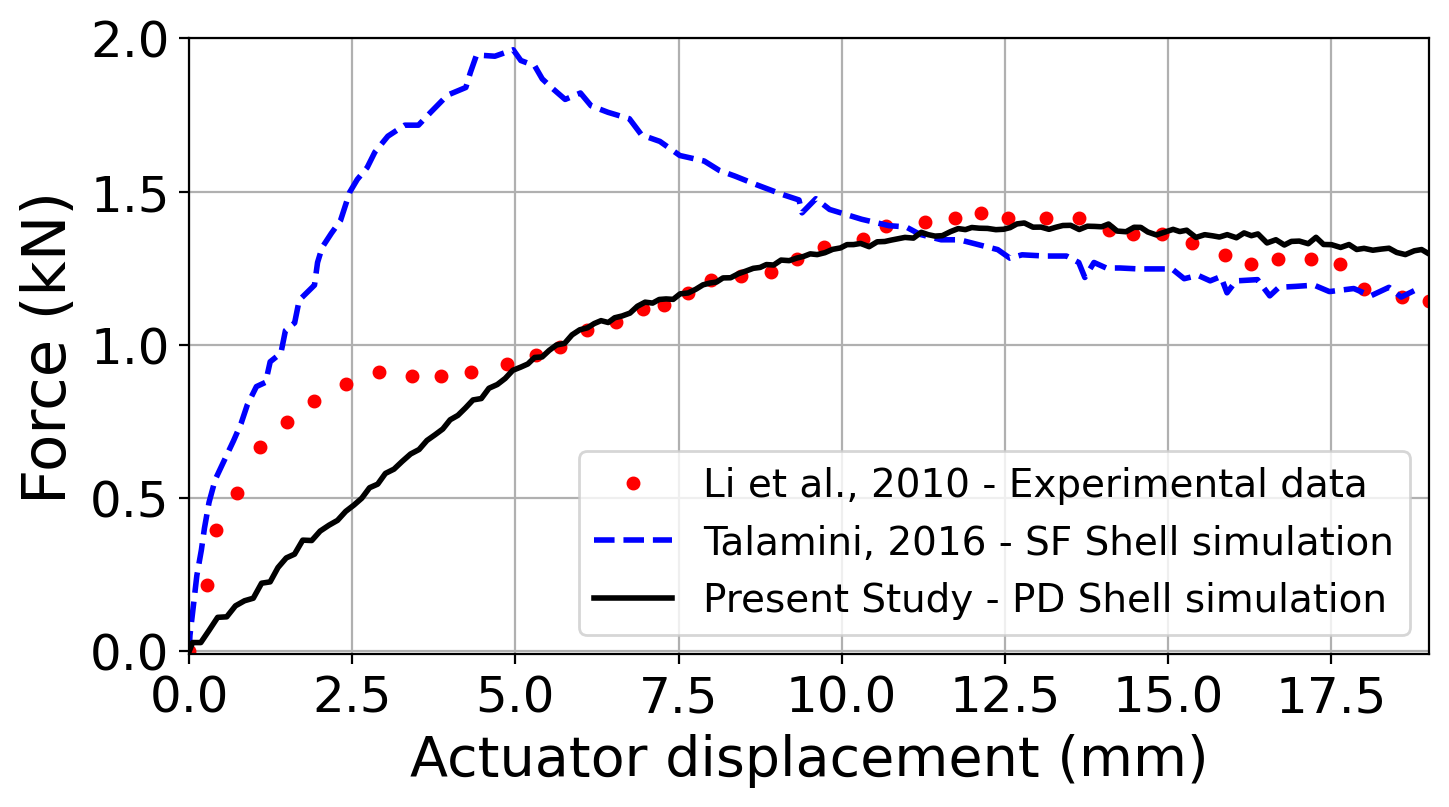}}
    
    \subfloat[][]{\includegraphics[width=0.9\columnwidth]{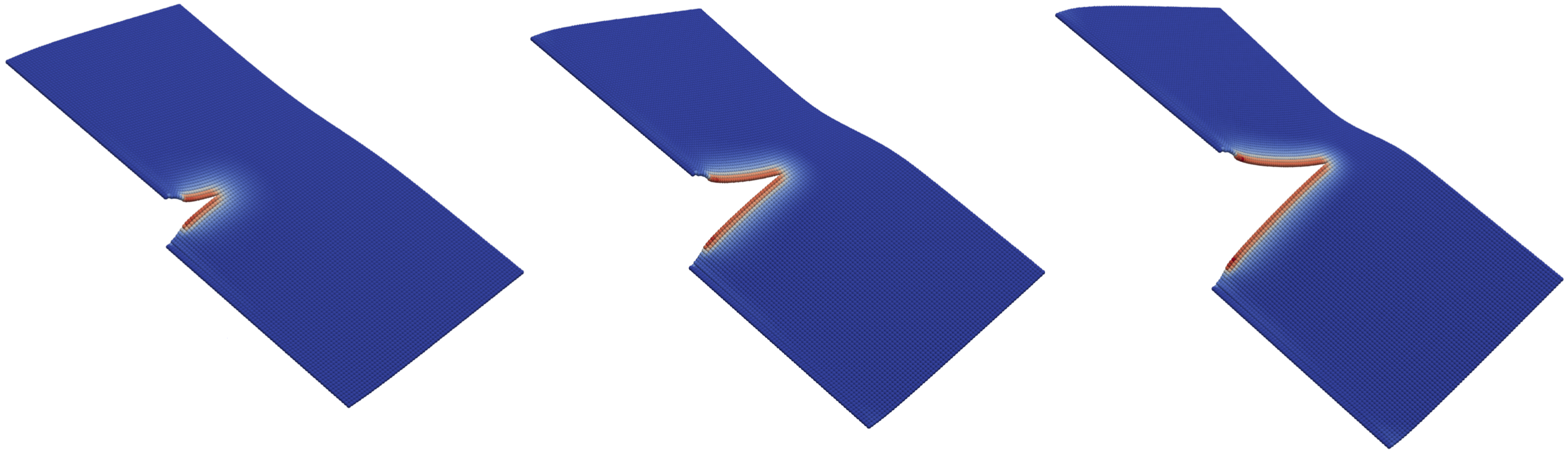}}
    \caption{Plate undergoing mode-III tearing fracture. (a): Description of the geometry and boundary conditions employed in the computation. (b): Comparison of the vertical load-deflection curves of the proposed PD shell with the experimental data~\cite{li2011mixed} and the Shear-Flexible (SF) Shell formulation~\cite{talamini2017parallel}. (v): Evolution of the plastic strain in the PD shell simulation.}
    \label{fig:tearing_problem}
\end{figure*}

We simulate a plate tearing experiment~\cite{li2011mixed} to show that the PD thin shell model is capable of capturing out-of-plane shear failure despite the fact that the through-thickness shear strains are not considered in the theory. The problem setup is shown in \cref{fig:tearing_problem}a, where a notched plate is subjected to a nominally mode-III loading condition, i.e., the upper left and lower left parts of the plate (above and below the notch) are constrained and displaced orthogonally to the plate in the opposite directions.
In the experimental setup, a gripping system involving multiple pins in the constrained region was used to clamp the left side of the plate. While the experimental intent was to achieve a perfectly rigid condition in the clamped region, it was found that there was substantial slippage near the pins that resulted in a fair amount of flexibility in that area. Thus, assuming a perfectly clamped BC in the modeling can degrade the results~\cite{li2011mixed}, which is also evident from the computations presented in~\cite{talamini2017parallel}. To take the grip elasticity into account, we consider a part of the clamped region to be {\em loosely constrained} with the stiffness of the plate linearly increasing from from $E_0$ to $18E_0$ in that region. The details of the computational setup are shown in~\cref{fig:tearing_problem}a. 

The plate is made of 6061-T6 aluminum alloy, which is modeled as an isotropic hardening material with the following properties: Young's modulus $E_0=72$ GPa, Poisson's ratio $\nu = 0.343$, and a linear hardening law $\sigma_Y(\bar{\epsilon}^P) = 320 + 894 \, \bar{\epsilon}^P$ MPa. A typical fracture strain for this alloy is around 0.2-0.25~\cite{zhu2011characterization}. The plasticity-driven failure approach is utilized as in \cref{eqn:plastic_strain} with the following threshold and critical values:
\begin{equation}
    \bar{\epsilon}^P_{\rm th}=0.23 , \ \ \bar{\epsilon}^P_{\rm cr}=0.27 .
\end{equation}
The computed vertical force-displacement is plotted in \cref{fig:tearing_problem}b and compared with both experimental data~\cite{li2011mixed} and the Shear-Flexible (SF) Shell simulation data~\cite{talamini2017parallel}. In the latter, the fracture is represented using a cohesive zone model. Our approach captures the peak force and the steady-state response of the structure well. The SF simulations, however, assume a perfectly rigid clamping boundary condition, which resulted in an over-prediction of the peak force.

\subsection{Beyond Inelasticity: Fragmentation}
\label{sec:fragmentation}

The two numerical examples included in this section are meant to demonstrate the ability of the proposed formulation to handle large-deformation phenomena with fragmentation.

\subsubsection{Impact of an Egg-Shaped Thin-Walled Object on a Rigid Wall}
\label{sec:impact_egg}

\begin{figure*}[!hbpt]
    \centering
    \includegraphics[width=1\columnwidth]{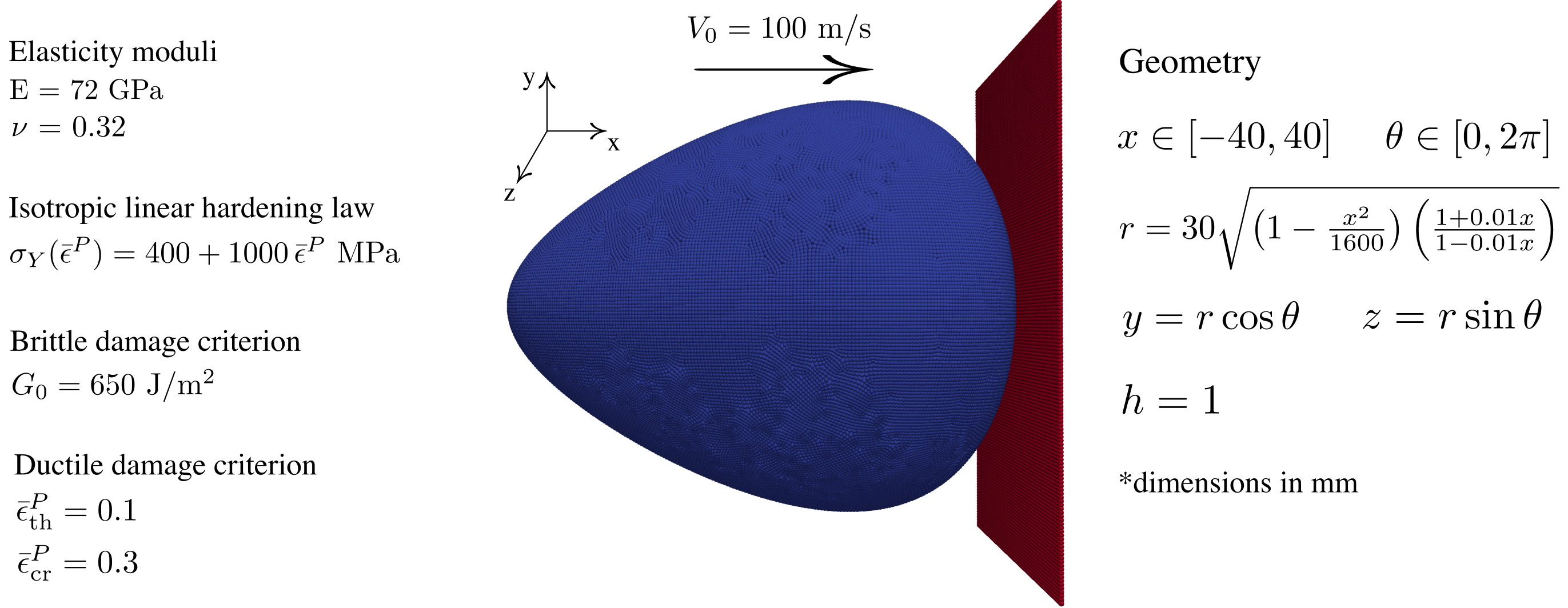}
    \caption{Impact of an egg-shaped shell structure on a rigid wall problem description. An elastic brittle material and a ductile elasto-plastic material are considered in the simulations. The egg is discretized using 83,000 PD nodes with the average node spacing of 0.5~mm.}
    \label{fig:impact_egg_setup}
\end{figure*}

In this problem, we consider high-speed impact of an egg-shaped shell structure on a rigid wall. The problem description including the geometry and material properties is provided in \cref{fig:impact_egg_setup}.  
Two scenarios are studied by modeling the shell structure as (a) an elastic brittle material (cf. \cref{sec:crack_branching}) and (b) a ductile elasto-plastic material (cf. \cref{sec:plastic_fracture}). Contact is modeled using a short-range repulsive-force approach~\cite{silling2005meshfree,kamensky2019peridynamic}. The simulations involved about 83,000 PD nodes and 50,000 time steps, which took approximately 1.5 hours to run on the Stampede 2 cluster at Texas Advanced Computing Center using 4 SKX compute nodes (192 processors in total).

\begin{figure*}[!hbpt]
    \centering
    \includegraphics[width=\columnwidth]{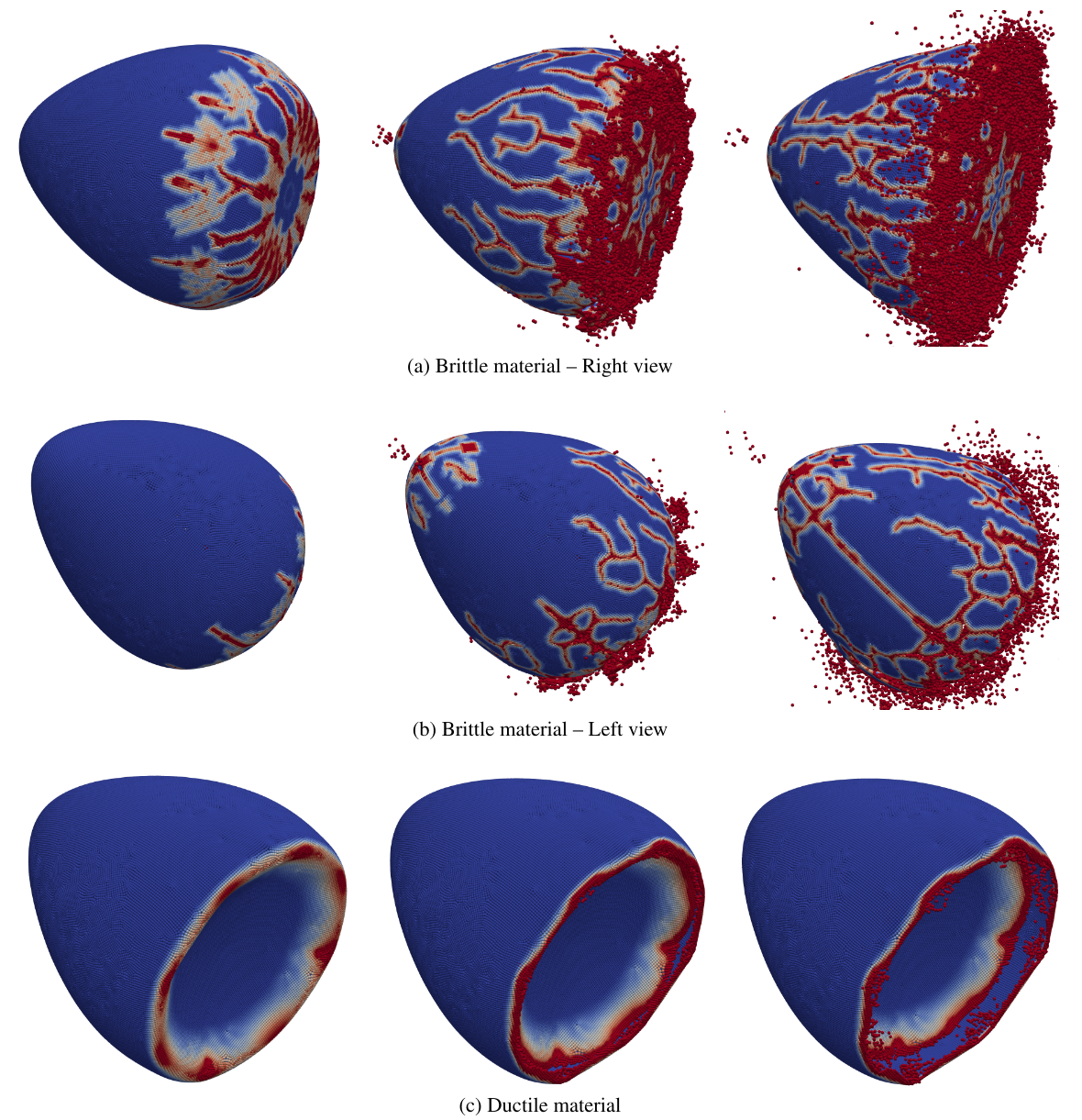}
    \caption{Impact of an egg-shaped shell structure on a rigid wall. Evolution of damage in the brittle (a--b) and ductile (c) materials.}
    \label{fig:impact_egg_damage}
\end{figure*}

\cref{fig:impact_egg_damage} shows the evolution of damage and failure in this problem. We note that, as expected~\cite{carmona2014fracture}, severe fragmentation and spalling on the side opposite the impact site occur in the brittle case, while a rupture-type failure occurs in the ductile case. Future studies are warranted to validate our formulation using experimental data that involve fragmentation. For this, a post-processing algorithm will need to be developed and utilized to obtain the histograms of the fragment size (see, e.g.,~\cite{littlewood2016identification,diehl2017extraction}) and compare with experiments (see, e.g.~\cite{vogler2003fragmentation,schraml2005simulating}).

\subsubsection{Impact and Penetration Problem}
\label{sec:bullet_penetration}

\begin{figure*}[!hbpt]
    \centering
    \includegraphics[width=\columnwidth]{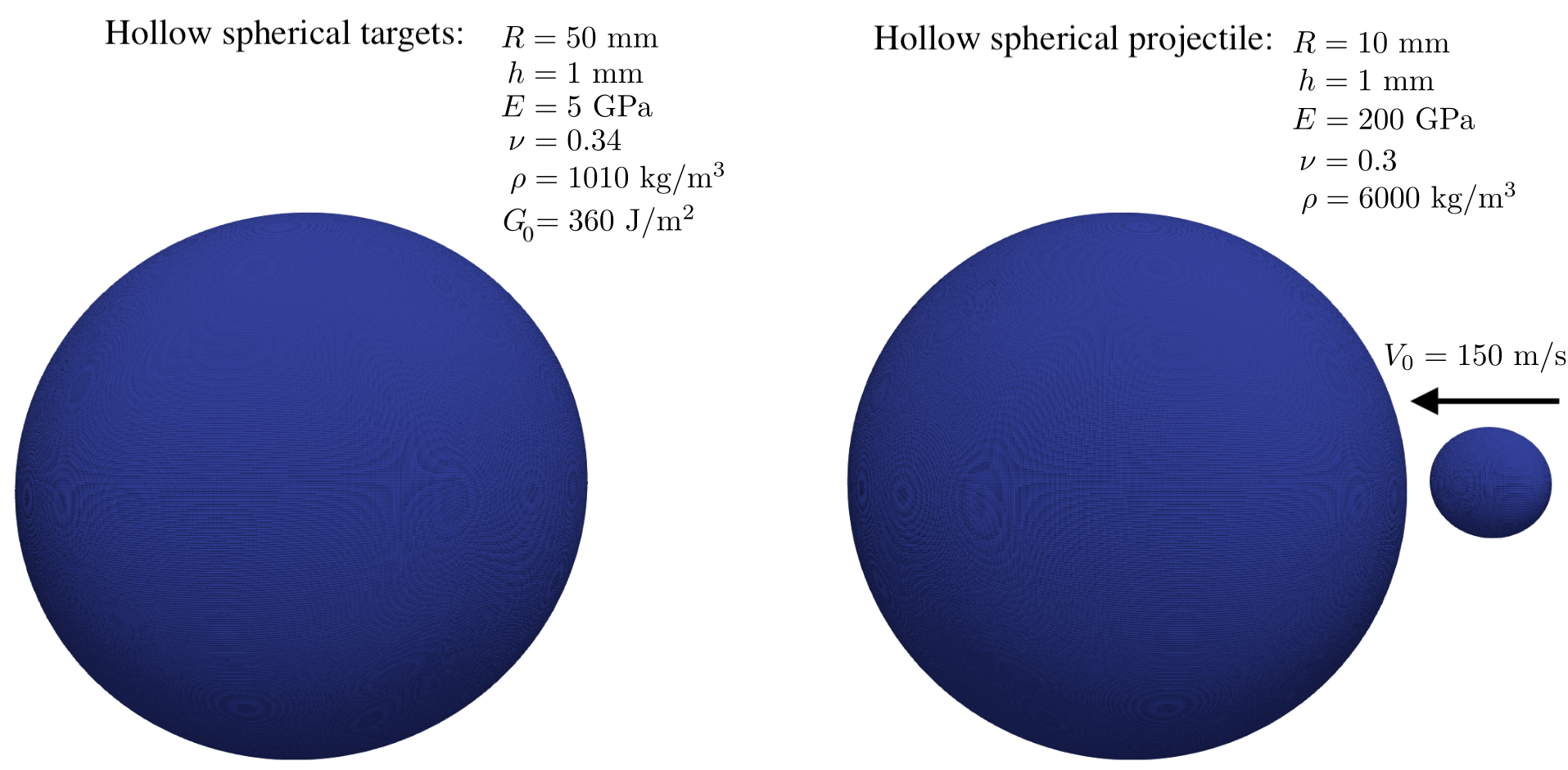}
    \caption{Pellet-like projectile impacting two back-to-back spherical targets. Problem description. The spherical target material is assumed to be elastic and brittle with the properties of PMMA.}
    \label{fig:bullet_penetration_setup}
\end{figure*}

We design a scenario where two identical back-to-back spherical shells are impacted by a hollow spherical pellet-like projectile. The problem description, including the geometric and material properties, is shown in \cref{fig:bullet_penetration_setup}. The material properties of the target structures are chosen as those of PMMA~\cite{mehrmashhadi2019uncovering}. An elastic brittle material is used to describe the spherical targets and the short-range repulsive-force-based approach~\cite{silling2005meshfree,kamensky2019peridynamic} is employed to model contact. The projectile is assumed to be elastic and intact. The simulation involved about 456,000 PD nodes and 16,500 time steps, which took approximately 1.5 hours to run on the same Stampede2 cluster using 10 SKX compute nodes or 480 processors in total.

\begin{figure*}[!hbpt]
    \centering
    \includegraphics[width=\columnwidth]{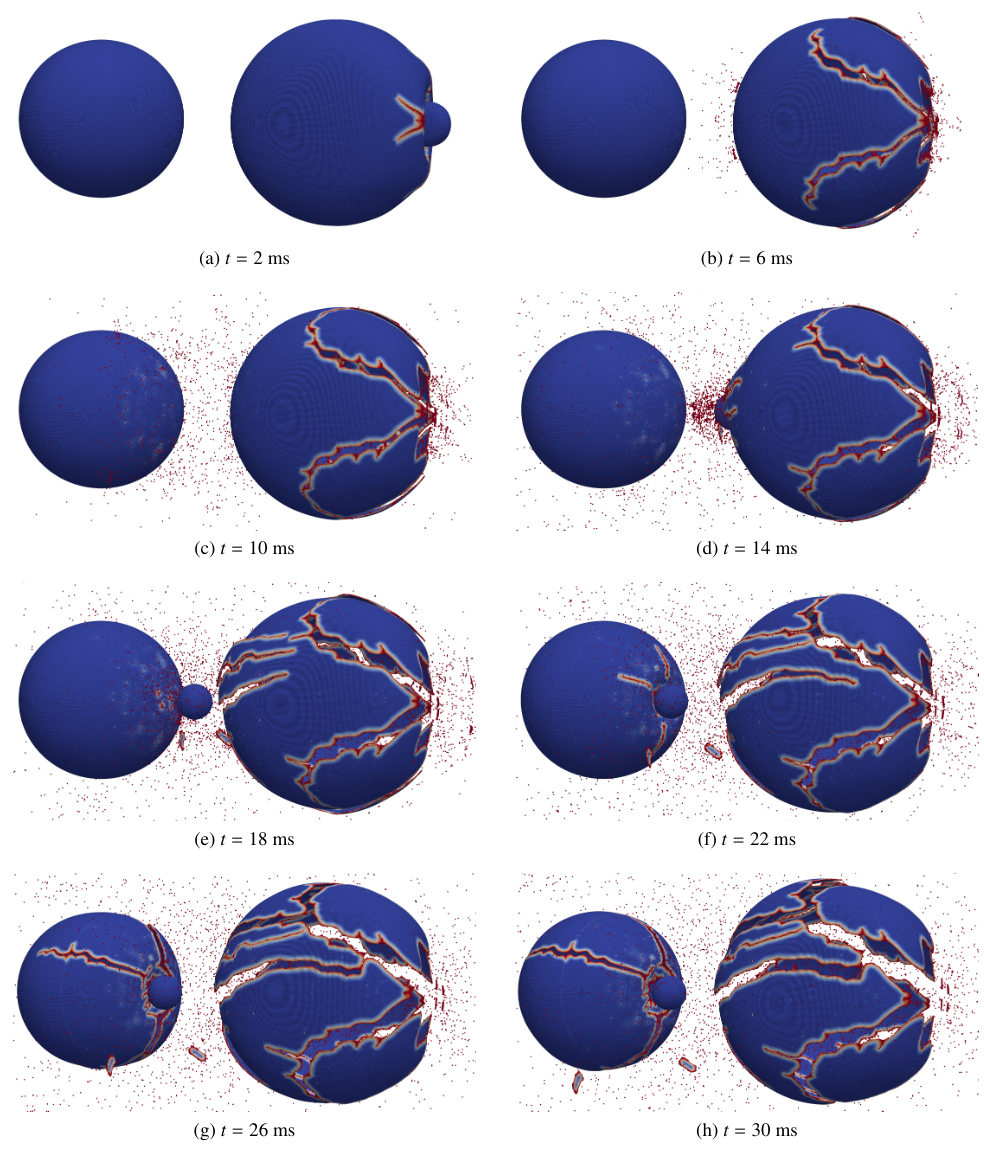}
    \caption{Damage plotted on the surfaces of the structure in the deformed and fragmented configuration.}
    \label{fig:bullet_penetration_damage}
\end{figure*}

\cref{fig:bullet_penetration_damage} shows the damage variable distribution on the structure surfaces in the deformed and fragmented configuration at different time instants. The simulation shows that the initial kinetic energy of the projectile is high enough to penetrate and exit the first target. The penetration event caused several dominant cracks, with micro-branches, to run along the sphere's surface. Additional fracturing is induced by the projectile's exit with fractures running in the opposite direction of the projectile's travel. Cracks due to the penetration and exit quickly merge causing the first sphere to fully fragment. The remaining projectile kinetic energy appears to be high enough to partially fracture, but not penetrate the second sphere. It is also worth noting that the debris from the first sphere caused minor damage to the second sphere and made it easier to fracture upon impact by the projectile.

\section{Conclusions}
\label{sec:conclusions}

In this work, a comprehensive PD-based framework is developed for the modeling of KL shells. There is no need for {\em a priori} global parameterization of the surface as we use a meshfree technique (PCA) to reconstruct manifolds locally at each nodal point using its neighborhood data. Classical KL kinematics are used in the correspondence PD framework. The rotation-free formulation is written in terms of the mid-surface velocity DOFs. Classical 3D rate-based constitutive models can be directly incorporated to handle different materials. A bond-stabilization technique (i.e., bond-associative modeling) is employed to attain numerical stability. Through studying several benchmark problems including elastostatics, dynamics, plasticity, and fracture growth, we demonstrate that the presented framework is of general applicability and is accurate, efficient, and robust in capturing large deformation and damage in thin-walled structures. It is also important to note that the present formulation enables the use of completely unstructured surface meshes to discretize a thin shell theory that uses higher-order derivatives of the kinematic fields.

We envision an extensive use of the proposed developments as a reliable method in a variety of applications for failure analysis of thin-shell structures as showcased in the last two numerical examples involving fracture and fragmentation. The present developments can be extended to other models. Semi-Lagrangian shell formulations based on this framework can be proposed with applications to extreme events. 

\section*{Acknowledgments}
\label{sec:acknowledge}
Y.~Bazilevs and M.~Alaydin were supported through the Sandia Contract No. 2111577. Y. Bazilevs and M. Behzadinasab were supported through the ONR Grant No. N00014-21-1-2670. N.~Trask acknowledges funding under the DOE ASCR PhILMS center (Grant number DE-SC001924) and the Laboratory Directed Research and Development program at Sandia National Laboratories. Sandia National Laboratories is a multi-mission laboratory managed and operated by National Technology and Engineering Solutions of Sandia, LLC., a wholly owned subsidiary of Honeywell International, Inc., for the U.S. Department of Energy’s National Nuclear Security Administration under contract DE-NA0003525. The authors acknowledge the Texas Advanced Computing Center (TACC) at The University of Texas at Austin for providing the HPC resources that have contributed to the research results reported within this paper.

\appendix

\section{Implementation Algorithm}
\label{sec:implementation}

A step-by-step guide for using the present formulation is provided in this section.\\ 

The following steps are only for the initialization phase:

\begin{steps}
\item Discretize the geometry and define the number of Gaussian quadrature points along the thickness direction $\xi_3$. (We use three Gauss points through the shell thickness in the examples presented.)

\item Obtain the local parametric coordinates $\xi_{1_{\rm \an{P-Q}}}$ and $\xi_{2_{\rm \an{P-Q}}}$ as described in \cref{sec:manifold}. 

\item Evaluate the first- and second-order gradient kernel functions $\left[\Phi_{\xi_1}, \, \Phi_{\xi_2}, \, \Phi_{\xi_1\xi_1}, \, \Phi_{\xi_1\xi_2}, \, \Phi_{\xi_2\xi_2}\right]_{\rm \an{P-Q}}$ as detailed in \cref{sec:PD_gradient}. Once damage grows, the weights should be recomputed.
\end{steps}

\noindent The following steps are performed at each time step:

\begin{steps}

\item Compute the parametric gradients of the mid-surface kinematics (i.e., position and velocity) at each PD node using \cref{eqn:dx2Ddxi,eqn:dv2Ddxi}.

\item Calculate the geometric entities at each PD node, i.e., the normal vector $\mathbf{n}_{({\rm P})}$ and the auxiliary tensors $\mathbf{A}_{({\rm P})}$, $\mathbf{B}^1_{({\rm P})}$ and $\mathbf{B}^2_{({\rm P})}$, using \cref{eqn:normal,eqn:A,eqn:B12}.

\item Obtain the nodal Jacobian matrices $\mathbf{F}^{\rm 3D}_{({\rm P})}$ at each through-thickness Gauss point $\xi_3$ using \cref{eqn:x3dd,eqn:defF}.

\item At each through-thickness Gauss point $\xi_3$, compute the auxiliary state $\boldsymbol{\beta}^{\rm 3D}_{\rm \an{P-Q}}$ for each PD bond using \cref{eqn:beta} and the nodal velocity gradient $v^{\rm 3D}_{{i,j}_{({\rm P})}}$ using \cref{eqn:dv3Ddx}. 

\item Calculate the auxiliary state $\boldsymbol{\gamma}^{\rm 3D}_{\rm \an{P-Q}}$, the position state $\mathbf{x}^{\rm 3D}_{\rm \an{P-Q}}$, and the velocity state $\mathbf{v}^{\rm 3D}_{\rm \an{P-Q}}$ using \cref{eqn:gamma}, \cref{eqn:x3Dpq}, and \cref{eqn:bondv3D}, respectively.

\item Compute the bond-associated velocity gradient $\mathbf{v}^{\rm 3D}_{\rm \an{P-Q}}$ using \cref{eqn:bondL}.

\item Use an objective stress rate and update the bond-associated Cauchy stress tensor $\boldsymbol{\sigma}^{\rm 3D}_{\rm \an{P-Q}}$ (cf. \cref{sec:stress_update,sec:plasticity}) and enforce the plane-stress condition as described in \cref{sec:zero_stress}.

\item Update the bond-associated Jacobian of the deformation gradient ${J}_{\rm \an{P-Q}}$ and the nodal thickness $h_{({\rm P})}$ as shown in \cref{sec:zero_stress}.

\item Calculate the bond-associated Kirchhoff stress tensor $\boldsymbol{\sigma}^{\rm 3D}_{\rm \an{P-Q}}$ using \cref{eqn:tau-sigma}.

\item Compute the auxiliary states $\mathbf{a}^{\rm 3D}_{\rm \an{P-Q}}$ and $\mathbf{b}^{\rm 3D}_{\rm \an{P-Q}}$ using \cref{eqn:ab}.

\item Obtain the force state $\mathbf{T}_{\rm \an{P-Q}}$ using \cref{eqn:force-discrete} and advance the PD equations of motion in time, i.e., \cref{eqn:PD_EOM}.

\item Evaluate the damage state at each PD bond (cf., \cref{sec:damage_modeling}).

\item If the damage grows, update the influence states using \cref{eqn:damagedOmega} and recalculate the gradient kernel functions for the next step.

\end{steps}

\section{Linearization}
\label{sec:linearization}

In this section, we linearize the PD equations of motion for an efficient implementation in (quasi-)static and implicit schemes. A left-hand-side (LHS) matrix $\mathbf{K}$ is calculated, which can be utilized in a predictor-multicorrector algorithm to reduce the residual vector $\mathbf{R}$ and update the solution field~\cite{chung1993time}. 

To linearize the formulation, we vary $\mathbf{v}^{\rm 2D}$ while keeping $\mathbf{x}^{\rm 2D}$ fixed, i.e., relaxing their dependence and treating the velocities as independent degrees of freedom in the variation phase. In the discrete setting, the internal force of a material point P at an iteration (k) is defined using \cref{eqn:PD_EOM}, i.e.,
\begin{equation}
f_{i_{\rm (P)}} = \sum_{\rm Q} \left( {T}_{i_{\rm \an{P-Q}}} - {T}_{i_{\rm \an{Q-P}}} \right) A_{\rm (Q)} ,
\label{eqn:f_def}
\end{equation}
where $A_{\rm (Q)}$ is the nodal area of Q.

In a static setting, the residual force for P at iteration (k) is defined as
\begin{equation}
R^{\rm (k)}_{i_{\rm (P)}} = \left( f^{\rm (k)}_{i_{\rm (P)}} + b_{i_{\rm (P)}} \right) A_{\rm (P)} .
\label{eqn:R_def}
\end{equation}
The variation of the residual vector with respect to the degrees of freedom defines the elements of the LHS matrix, i.e.,
\begin{equation}
K_{ij_{\rm \an{P-R}}} = \frac{\partial R^{\rm (k)}_{i_{\rm (P)}}}{\partial v^{\rm 2D}_{j_{\rm (R)}}} .
\label{eqn:K_def}
\end{equation}
The velocity increments $\Delta v^{\rm 2D}$ are calculated by solving the following linear system:
\begin{equation}
K_{ij_{\rm \an{P-R}}} \, \Delta v^{\rm 2D^{\rm (k+1)}}_{j_{\rm (R)}} = - R^{\rm (k)}_{i_{\rm (P)}} ,
\end{equation}
then,
\begin{equation}
\mathbf{v}^{\rm 2D^{\rm (k+1)}} = \mathbf{v}^{\rm 2D^{\rm (k)}} + \Delta \mathbf{v}^{\rm 2D^{\rm (k+1)}} .
\end{equation}
The positions are updated depending on the integration scheme. For example, for an Euler method, 
\begin{equation}
\mathbf{x}^{\rm 2D^{\rm (k+1)}} = \mathbf{x}^{\rm 2D^{\rm (k)}} + \Delta \mathbf{v}^{\rm 2D^{\rm (k+1)}} \, \Delta t .
\end{equation}

Using \cref{eqn:f_def,eqn:R_def,eqn:K_def}, $\mathbf{K}$ is obtained as
\begin{equation}
K_{ij_{\rm \an{P-R}}} = \sum_{\rm Q} \left( \frac{\partial {T}_{i_{\rm \an{P-Q}}}}{\partial v^{\rm 2D}_{j_{\rm (R)}}}  - \frac{\partial {T}_{i_{\rm \an{Q-P}}}}{\partial v^{\rm 2D}_{j_{\rm (R)}}}  \right) A_{\rm (Q)} \, A_{\rm (P)} ,
\label{eqn:K_lin}
\end{equation}
therefore, the variation of the force state with respect to the velocities should be considered. 

Using \cref{eqn:A,eqn:B12,eqn:defF} and observing that the matrices $\mathbf{A}$, $\mathbf{B}^1$, $\mathbf{B}^2$, and $\mathbf{F}^{\rm 3D}$ depend only on the positions, then using \cref{eqn:beta,eqn:gamma}, 
\begin{equation} 
\begin{aligned}
& \frac{\partial \beta^{{I}}_{{lik}_{\rm \an{P-Q}}}}{\partial v^{\rm 2D}_{j_{\rm (R)}}} = 0 , \\
& \frac{\partial \gamma^{{I}}_{{ki}_{\rm \an{P-Q}}}}{\partial v^{\rm 2D}_{j_{\rm (R)}}} = 0 , \\
& \frac{\partial \left(F^{{I}}_{{km}_{({\rm P})}}\right)^{-1}}{\partial v^{\rm 2D}_{j_{\rm (R)}}} = 0 .
\end{aligned}
\label{eqn:del-beta-gamma-F}
\end{equation}
Using \cref{eqn:force-discrete,eqn:del-beta-gamma-F} and the chain rule,
\begin{equation}
\begin{aligned}
& \frac{\partial {T}_{i_{\rm \an{P-Q}}}}{\partial v^{\rm 2D}_{j_{\rm (R)}}} = \sum_{I=1}^{N_L} w^I \, \frac{\partial {T}_{i_{\rm \an{P-Q}}}^I}{\partial v^{\rm 2D}_{j_{\rm (R)}}} , \\
& \frac{\partial {T}_{i_{\rm \an{P-Q}}}^I}{\partial v^{\rm 2D}_{j_{\rm (R)}}} = 
\frac{\partial {a}^{{I}}_{i_{\rm \an{P-Q}}}}{\partial v^{\rm 2D}_{j_{\rm (R)}}}
+ \frac{\partial \bar{a}^{{I}}_{k_{\rm (P)}}}{\partial v^{\rm 2D}_{j_{\rm (R)}}} \, \gamma^{{I}}_{{ki}_{\rm \an{P-Q}}} 
+ \frac{\partial \bar{b}^{{I}}_{lk_{\rm (P)}}}{\partial v^{\rm 2D}_{j_{\rm (R)}}} \, \beta^{{I}}_{{lik}_{\rm \an{P-Q}}} , \\
& \frac{\partial \bar{a}^{{I}}_{k_{\rm (P)}}}{\partial v^{\rm 2D}_{j_{\rm (R)}}} = 
\sum_{\rm S} \Bigg( \frac{\partial {a}^{{I}}_{k_{\rm \an{S-P}}}}{\partial v^{\rm 2D}_{j_{\rm (R)}}} 
- \frac{\partial {a}^{{I}}_{k_{\rm \an{P-S}}}}{\partial v^{\rm 2D}_{j_{\rm (R)}}} \Bigg) A_{\rm (S)} ,\\
& \frac{\partial \bar{b}^{{I}}_{lk_{\rm (P)}}}{\partial v^{\rm 2D}_{j_{\rm (R)}}} = 
\left(F^{{I}}_{{km}_{({\rm P})}}\right)^{-1} \sum_{\rm S} \Bigg( \frac{\partial {b}^{{I}}_{lm_{\rm \an{P-S}}}}{\partial v^{\rm 2D}_{j_{\rm (R)}}} + \frac{\partial {b}^{{I}}_{lm_{\rm \an{S-P}}}}{\partial v^{\rm 2D}_{j_{\rm (R)}}} \Bigg) A_{\rm (S)} .
\end{aligned}
\label{eqn:del_T}
\end{equation}
Using \cref{eqn:ab}, 
\begin{equation} 
\begin{aligned}  
&\frac{\partial {a}^{\rm 3D}_{i_{\rm \an{P-Q}}}}{\partial v^{\rm 2D}_{j_{\rm (R)}}} = \frac{h_{\rm (P)}}{2} \, \alpha_{\rm \an{P-Q}} \, \frac{\partial {\sigma}^{\rm 3D}_{ik_{\rm \an{P-Q}}}}{\partial v^{\rm 2D}_{j_{\rm (R)}}} \frac{{x}^{\rm 3D}_{k_{\rm \an{P-Q}}}}{|{x}^{\rm 3D}_{\rm \an{P-Q}}|^2} , \\
&\frac{\partial {b}^{\rm 3D}_{il_{\rm \an{P-Q}}}}{\partial v^{\rm 2D}_{j_{\rm (R)}}} = \frac{h_{\rm (P)}}{4} \, \alpha_{\rm \an{P-Q}} \, \frac{\partial {\sigma}^{\rm 3D}_{ik_{\rm \an{P-Q}}}}{\partial v^{\rm 2D}_{j_{\rm (R)}}} \Bigg( \delta_{kl} - \frac{{x}^{\rm 3D}_{k_{\rm \an{P-Q}}} \, {x}^{\rm 3D}_{l_{\rm \an{P-Q}}}}{|{x}^{\rm 3D}_{\rm \an{P-Q}}|^2} \Bigg) .
\end{aligned}  
\label{eqn:del-ab}
\end{equation}
Using the chain rule,
\begin{equation} 
\begin{aligned}  
\frac{\partial {\sigma}^{\rm 3D}_{ik_{\rm \an{P-Q}}}}{\partial v^{\rm 2D}_{j_{\rm (R)}}} = \frac{\partial {\sigma}^{\rm 3D}_{ik_{\rm \an{P-Q}}}}{\partial {v}^{\rm 3D}_{t,u_{\rm \an{P-Q}}}} \frac{\partial {v}^{\rm 3D}_{t,u_{\rm \an{P-Q}}}}{\partial v^{\rm 2D}_{j_{\rm (R)}}} .
\end{aligned}  
\label{eqn:del-sigma-v2D}
\end{equation}
For elasticity, according to the classical law, the stress-strain relation reads: 
\begin{equation} 
\begin{aligned}  
\frac{\partial {\sigma}^{\rm 3D}_{ik_{\rm \an{P-Q}}}}{\partial {v}^{\rm 3D}_{t,u_{\rm \an{P-Q}}}} = \mathbb{C}_{iktu} \, \Delta t,
\end{aligned}  
\label{eqn:del-sigma}
\end{equation}
where $\mathbb{C}$ is the fourth-order material tangent stiffness tensor. For elastoplastic analysis, a consistent stiffness tensor should be used. Using \cref{eqn:bondL},
\begin{equation} 
\frac{\partial {v}^{\rm 3D}_{t,u_{\rm \an{P-Q}}}}{\partial v^{\rm 2D}_{j_{\rm (R)}}} = \frac{{x}^{\rm 3D}_{u_{\rm \an{P-Q}}}}{|{x}^{\rm 3D}_{\rm \an{P-Q}}|^2} \frac{\partial {v}^{\rm 3D}_{t_{\rm \an{P-Q}}}}{\partial v^{\rm 2D}_{j_{\rm (R)}}} + \frac{1}{2} \left(\delta_{vu} - \frac{{x}^{\rm 3D}_{v_{\rm \an{P-Q}}} \, {x}^{\rm 3D}_{u_{\rm \an{P-Q}}}}{|{x}^{\rm 3D}_{\rm \an{P-Q}}|^2} \right) \left( \frac{\partial v^{\rm 3D}_{{t,v}_{({\rm P})}}}{\partial v^{\rm 2D}_{j_{\rm (R)}}} + \frac{\partial v^{\rm 3D}_{{t,v}_{({\rm Q})}}}{\partial v^{\rm 2D}_{j_{\rm (R)}}} \right) .
\label{eqn:del-bondL}
\end{equation}
Using \cref{eqn:bondv3D,eqn:del-beta-gamma-F},
\begin{equation} 
\begin{aligned}
\frac{\partial {v}^{\rm 3D}_{t_{\rm \an{P-Q}}}}{\partial v^{\rm 2D}_{j_{\rm (R)}}}
& = \delta_{tj} (\delta_{\rm QR} - \delta_{\rm PR}) 
+ \left( \sum_{\rm S} \gamma^{\rm 3D}_{{tv}_{\rm \an{Q-S}}} \, \delta_{vj} (\delta_{\rm SR} - \delta_{\rm QR}) \, A_{\rm (S)} - \sum_{\rm S} \gamma^{\rm 3D}_{{tv}_{\rm \an{P-S}}} \, \delta_{vj} (\delta_{\rm SR} - \delta_{\rm PR}) \, A_{\rm (S)} \right) \\
& = \delta_{tj} (\delta_{\rm QR} - \delta_{\rm PR}) 
+ \left( \gamma^{\rm 3D}_{{tj}_{\rm \an{Q-R}}} - \gamma^{\rm 3D}_{{tj}_{\rm \an{P-R}}} \right) A_{\rm (R)} - \delta_{\rm QR} \sum_{\rm S} \gamma^{\rm 3D}_{{tj}_{\rm \an{Q-S}}} \, A_{\rm (S)} + \delta_{\rm PR} \sum_{\rm S} \gamma^{\rm 3D}_{{tj}_{\rm \an{P-S}}} \, A_{\rm (S)} .
\end{aligned}
\label{eqn:del-bondv3D}
\end{equation}
Using \cref{eqn:dv3Ddx,eqn:del-beta-gamma-F},
\begin{equation}
\begin{aligned}
\frac{\partial v^{\rm 3D}_{{t,v}_{({\rm P})}}}{\partial v^{\rm 2D}_{j_{\rm (R)}}} 
& = \left( \sum_{\rm S}  \beta^{\rm 3D}_{{trw}_{\rm \an{P-S}}} \, \delta_{rj} (\delta_{\rm SR} - \delta_{\rm PR}) \, A_{\rm (S)} \right) F^{\rm 3D^{-1}}_{{wv}_{({\rm P})}} \\
& = F^{\rm 3D^{-1}}_{{wv}_{({\rm P})}} \left( \beta^{\rm 3D}_{{tjw}_{\rm \an{P-R}}} \, A_{\rm (R)} - \delta_{\rm PR} \sum_{\rm S} \beta^{\rm 3D}_{{tjw}_{\rm \an{P-S}}} \, A_{\rm (S)} \right) .
\end{aligned}
\label{eqn:del-dv3Ddx}
\end{equation}
Introducing \cref{eqn:del-dv3Ddx,eqn:del-bondv3D,eqn:del-bondL,eqn:del-sigma,eqn:del-sigma-v2D,eqn:del-ab,eqn:del_T} in \cref{eqn:K_lin}, the tangent matrix $\mathbf{K}$ can be obtained.

\bibliographystyle{plainnat}
\bibliography{main}

\begin{thebibliography}{123}
\providecommand{\natexlab}[1]{#1}
\providecommand{\url}[1]{\texttt{#1}}
\expandafter\ifx\csname urlstyle\endcsname\relax
  \providecommand{\doi}[1]{doi: #1}\else
  \providecommand{\doi}{doi: \begingroup \urlstyle{rm}\Url}\fi

\bibitem[Alaydin et~al.(2021)Alaydin, Benson, and Bazilevs]{alaydin2021updated}
Mert~D Alaydin, David~J Benson, and Yuri Bazilevs.
\newblock An updated {Lagrangian} framework for {Isogeometric Kirchhoff--Love}
  thin-shell analysis.
\newblock \emph{Computer Methods in Applied Mechanics and Engineering},
  384:\penalty0 113977, 2021.

\bibitem[Ambati and De~Lorenzis(2016)]{ambati2016phase}
Marreddy Ambati and Laura De~Lorenzis.
\newblock Phase-field modeling of brittle and ductile fracture in shells with
  isogeometric {NURBS}-based solid-shell elements.
\newblock \emph{Computer Methods in Applied Mechanics and Engineering},
  312:\penalty0 351--373, 2016.

\bibitem[Ambati et~al.(2018)Ambati, Kiendl, and
  De~Lorenzis]{ambati2018isogeometric}
Marreddy Ambati, Josef Kiendl, and Laura De~Lorenzis.
\newblock Isogeometric {Kirchhoff--Love} shell formulation for
  elasto-plasticity.
\newblock \emph{Computer Methods in Applied Mechanics and Engineering},
  340:\penalty0 320--339, 2018.

\bibitem[Amenta and Kil(2004)]{amenta2004defining}
Nina Amenta and Yong~Joo Kil.
\newblock Defining point-set surfaces.
\newblock \emph{ACM Transactions on Graphics (TOG)}, 23\penalty0 (3):\penalty0
  264--270, 2004.

\bibitem[Bazilevs et~al.(2009)Bazilevs, Hsu, Benson, Sankaran, and
  Marsden]{bazilevs2009computational}
Yuri Bazilevs, Ming-Chen Hsu, David~J Benson, Sethu Sankaran, and Alison~L
  Marsden.
\newblock Computational fluid--structure interaction: methods and application
  to a total cavopulmonary connection.
\newblock \emph{Computational Mechanics}, 45\penalty0 (1):\penalty0 77--89,
  2009.

\bibitem[Behzadinasab(2020)]{behzadinasab2020peridynamic}
Masoud Behzadinasab.
\newblock \emph{Peridynamic modeling of large deformation and ductile
  fracture}.
\newblock PhD thesis, The University of Texas at Austin, 2020.

\bibitem[Behzadinasab and Foster(2019)]{behzadinasab2019third}
Masoud Behzadinasab and John~T Foster.
\newblock {The third Sandia Fracture Challenge}: peridynamic blind prediction
  of ductile fracture characterization in additively manufactured metal.
\newblock \emph{International Journal of Fracture}, 218\penalty0 (1):\penalty0
  97--109, 2019.

\bibitem[Behzadinasab and
  Foster(2020{\natexlab{a}})]{behzadinasab2020revisiting}
Masoud Behzadinasab and John~T Foster.
\newblock Revisiting the third {Sandia Fracture Challenge}: a bond-associated,
  semi-{Lagrangian} peridynamic approach to modeling large deformation and
  ductile fracture.
\newblock \emph{International Journal of Fracture}, 224:\penalty0 261--267,
  2020{\natexlab{a}}.

\bibitem[Behzadinasab and Foster(2020{\natexlab{b}})]{behzadinasab2020semi}
Masoud Behzadinasab and John~T Foster.
\newblock A semi-{Lagrangian} constitutive correspondence framework for
  peridynamics.
\newblock \emph{Journal of the Mechanics and Physics of Solids}, 137:\penalty0
  103862, 2020{\natexlab{b}}.

\bibitem[Behzadinasab and
  Foster(2020{\natexlab{c}})]{behzadinasab2020stability}
Masoud Behzadinasab and John~T Foster.
\newblock On the stability of the generalized, finite deformation
  correspondence model of peridynamics.
\newblock \emph{International Journal of Solids and Structures}, 182:\penalty0
  64--76, 2020{\natexlab{c}}.

\bibitem[Behzadinasab et~al.(2018)Behzadinasab, Vogler, Peterson, Rahman, and
  Foster]{behzadinasab2018peridynamics}
Masoud Behzadinasab, Tracy~J Vogler, Amanda~M Peterson, Rezwanur Rahman, and
  John~T Foster.
\newblock Peridynamics modeling of a shock wave perturbation decay experiment
  in granular materials with intra-granular fracture.
\newblock \emph{Journal of Dynamic Behavior of Materials}, 4\penalty0
  (4):\penalty0 529--542, 2018.

\bibitem[Behzadinasab et~al.(2021{\natexlab{a}})Behzadinasab, Foster, and
  Bazilevs]{behzadinasab2021unifiedII}
Masoud Behzadinasab, John~T Foster, and Yuri Bazilevs.
\newblock A unified, stable and accurate meshfree framework for peridynamic
  correspondence modeling—{Part II: Wave} propagation and enforcement of
  stress boundary conditions.
\newblock \emph{Journal of Peridynamics and Nonlocal Modeling}, 3:\penalty0
  46–--66, 2021{\natexlab{a}}.

\bibitem[Behzadinasab et~al.(2021{\natexlab{b}})Behzadinasab, Trask, and
  Bazilevs]{behzadinasab2021unifiedI}
Masoud Behzadinasab, Nathaniel~A Trask, and Yuri Bazilevs.
\newblock A unified, stable and accurate meshfree framework for peridynamic
  correspondence modeling—{Part I: Core} methods.
\newblock \emph{Journal of Peridynamics and Nonlocal Modeling}, 3:\penalty0
  24--45, 2021{\natexlab{b}}.

\bibitem[Belytschko et~al.(1985)Belytschko, Stolarski, Liu, Carpenter, and
  Ong]{belytschko1985stress}
Ted Belytschko, Henryk Stolarski, Wing~Kam Liu, Nicholas Carpenter, and Jame~SJ
  Ong.
\newblock Stress projection for membrane and shear locking in shell finite
  elements.
\newblock \emph{Computer Methods in Applied Mechanics and Engineering},
  51\penalty0 (1-3):\penalty0 221--258, 1985.

\bibitem[Belytschko et~al.(2013)Belytschko, Liu, Moran, and
  Elkhodary]{belytschko2013nonlinear}
Ted Belytschko, Wing~Kam Liu, Brian Moran, and Khalil Elkhodary.
\newblock \emph{Nonlinear finite elements for continua and structures}.
\newblock John wiley \& sons, 2013.

\bibitem[Benson et~al.(2011)Benson, Bazilevs, Hsu, and Hughes]{benson2011large}
David~J Benson, Yuri Bazilevs, Ming-Chen Hsu, and Thomas~JR Hughes.
\newblock A large deformation, rotation-free, isogeometric shell.
\newblock \emph{Computer Methods in Applied Mechanics and Engineering},
  200\penalty0 (13-16):\penalty0 1367--1378, 2011.

\bibitem[Bischoff et~al.(2018)Bischoff, Ramm, and
  Irslinger]{bischoff2018models}
Manfred Bischoff, E~Ramm, and J~Irslinger.
\newblock Models and finite elements for thin-walled structures.
\newblock \emph{Encyclopedia of Computational Mechanics Second Edition}, pages
  1--86, 2018.

\bibitem[Bobaru and Zhang(2015)]{bobaru2015cracks}
Florin Bobaru and Guanfeng Zhang.
\newblock Why do cracks branch{?} a peridynamic investigation of dynamic
  brittle fracture.
\newblock \emph{International Journal of Fracture}, 196\penalty0
  (1-2):\penalty0 59--98, 2015.

\bibitem[Bobaru et~al.(2009)Bobaru, Yang, Alves, Silling, Askari, and
  Xu]{bobaru2009convergence}
Florin Bobaru, Mijia Yang, Leonardo~Frota Alves, Stewart~A Silling, Ebrahim
  Askari, and Jifeng Xu.
\newblock Convergence, adaptive refinement, and scaling in 1d peridynamics.
\newblock \emph{International Journal for Numerical Methods in Engineering},
  77\penalty0 (6):\penalty0 852--877, 2009.

\bibitem[Bobaru et~al.(2016)Bobaru, Foster, Geubelle, and
  Silling]{bobaru2016handbook}
Florin Bobaru, John~T Foster, Philippe~H Geubelle, and Stewart~A Silling.
\newblock \emph{Handbook of peridynamic modeling}.
\newblock CRC press, 2016.

\bibitem[Bobaru et~al.(2018)Bobaru, Mehrmashhadi, Chen, and
  Niazi]{bobaru2018intraply}
Florin Bobaru, Javad Mehrmashhadi, Ziguang Chen, and Sina Niazi.
\newblock Intraply fracture in fiber-reinforced composites: {A} peridynamic
  analysis.
\newblock In \emph{ASC 33rd Annual Technical Conference \& 18th US-Japan
  Conference on Composite Materials, Seattle}, page~9, 2018.

\bibitem[Bowden et~al.(1967)Bowden, Brunton, Field, and
  Heyes]{bowden1967controlled}
F.P. Bowden, J.H. Brunton, J.E. Field, and A.D. Heyes.
\newblock Controlled fracture of brittle solids and interruption of electrical
  current.
\newblock \emph{Nature}, 216\penalty0 (5110):\penalty0 38--42, 1967.

\bibitem[Boyce et~al.(2014)Boyce, Kramer, Fang, Cordova, Neilsen, Dion,
  Kaczmarowski, Karasz, Xue, Gross, et~al.]{boyce2014sandia}
Brad~L Boyce, Sharlotte~LB Kramer, H~Eliot Fang, Theresa~E Cordova, Michael~K
  Neilsen, K~Dion, Amy~K Kaczmarowski, Erin Karasz, Liang Xue, Andrew~J Gross,
  et~al.
\newblock {The Sandia Fracture Challenge}: blind round robin predictions of
  ductile tearing.
\newblock \emph{International Journal of Fracture}, 186\penalty0
  (1-2):\penalty0 5--68, 2014.

\bibitem[Breitenfeld et~al.(2014)Breitenfeld, Geubelle, Weckner, and
  Silling]{breitenfeld2014non}
Michael~S Breitenfeld, Philippe~H Geubelle, Olaf Weckner, and Stewart~A
  Silling.
\newblock Non-ordinary state-based peridynamic analysis of stationary crack
  problems.
\newblock \emph{Computer Methods in Applied Mechanics and Engineering},
  272:\penalty0 233--250, 2014.

\bibitem[Breitzman and Dayal(2018)]{breitzman2018bond}
Timothy Breitzman and Kaushik Dayal.
\newblock Bond-level deformation gradients and energy averaging in
  peridynamics.
\newblock \emph{Journal of the Mechanics and Physics of Solids}, 110:\penalty0
  192--204, 2018.

\bibitem[Carmona et~al.(2014)Carmona, Wittel, and Kun]{carmona2014fracture}
Humberto~A Carmona, Falk~K Wittel, and Ferenc Kun.
\newblock From fracture to fragmentation: {Discrete} element modeling.
\newblock \emph{The European Physical Journal Special Topics}, 223\penalty0
  (11):\penalty0 2369--2382, 2014.

\bibitem[Chen(2018)]{chen2018bond}
Hailong Chen.
\newblock Bond-associated deformation gradients for peridynamic correspondence
  model.
\newblock \emph{Mechanics Research Communications}, 90:\penalty0 34--41, 2018.

\bibitem[Chen and Bobaru(2015)]{chen2015peridynamic}
Ziguang Chen and Florin Bobaru.
\newblock Peridynamic modeling of pitting corrosion damage.
\newblock \emph{Journal of the Mechanics and Physics of Solids}, 78:\penalty0
  352--381, 2015.

\bibitem[Chen et~al.(2019)Chen, Niazi, and Bobaru]{chen2019peridynamic}
Ziguang Chen, Sina Niazi, and Florin Bobaru.
\newblock A peridynamic model for brittle damage and fracture in porous
  materials.
\newblock \emph{International Journal of Rock Mechanics and Mining Sciences},
  122:\penalty0 104059, 2019.

\bibitem[Chi et~al.(2013)Chi, Chen, Hu, and Yang]{chi2013gradient}
Sheng-Wei Chi, Jiun-Shyan Chen, Hsin-Yun Hu, and Judy~P Yang.
\newblock A gradient reproducing kernel collocation method for boundary value
  problems.
\newblock \emph{International Journal for Numerical Methods in Engineering},
  93\penalty0 (13):\penalty0 1381--1402, 2013.

\bibitem[Chowdhury et~al.(2016)Chowdhury, Roy, Roy, and
  Reddy]{chowdhury2016peridynamic}
Shubhankar~R Chowdhury, Pranesh Roy, Debasish Roy, and J.N. Reddy.
\newblock A peridynamic theory for linear elastic shells.
\newblock \emph{International Journal of Solids and Structures}, 84:\penalty0
  110--132, 2016.

\bibitem[Chowdhury et~al.(2019)Chowdhury, Roy, Roy, and
  Reddy]{chowdhury2019modified}
Shubhankar~R Chowdhury, Pranesh Roy, Debasish Roy, and J.N. Reddy.
\newblock A modified peridynamics correspondence principle: {Removal} of
  zero-energy deformation and other implications.
\newblock \emph{Computer Methods in Applied Mechanics and Engineering},
  346:\penalty0 530--549, 2019.

\bibitem[Chung and Hulbert(1993)]{chung1993time}
Jintai Chung and Gregory~M Hulbert.
\newblock A time integration algorithm for structural dynamics with improved
  numerical dissipation: the generalized-$\alpha$ method.
\newblock \emph{Journal of Applied Mechanics}, pages 371--375, 1993.

\bibitem[Coox et~al.(2017)Coox, Maurin, Greco, Deckers, Vandepitte, and
  Desmet]{coox2017flexible}
Laurens Coox, Florian Maurin, Francesco Greco, Elke Deckers, Dirk Vandepitte,
  and Wim Desmet.
\newblock A flexible approach for coupling {NURBS} patches in rotationless
  isogeometric analysis of {Kirchhoff--Love} shells.
\newblock \emph{Computer Methods in Applied Mechanics and Engineering},
  325:\penalty0 505--531, 2017.

\bibitem[Cottrell et~al.(2009)Cottrell, Hughes, and
  Bazilevs]{cottrell2009isogeometric}
John~A Cottrell, Thomas~JR Hughes, and Yuri Bazilevs.
\newblock \emph{Isogeometric analysis: toward integration of CAD and FEA}.
\newblock John Wiley \& Sons, 2009.

\bibitem[de~Souza~Neto et~al.(2011)de~Souza~Neto, Peric, and
  Owen]{de2011computational}
Eduardo~A de~Souza~Neto, Djordje Peric, and David~RJ Owen.
\newblock \emph{Computational methods for plasticity: {Theory} and
  applications}.
\newblock John Wiley \& Sons, 2011.

\bibitem[Diehl et~al.(2017)Diehl, Bu{\ss}ler, Pfl{\"u}ger, Frey, Ertl, Sadlo,
  and Schweitzer]{diehl2017extraction}
Patrick Diehl, Michael Bu{\ss}ler, Dirk Pfl{\"u}ger, Steffen Frey, Thomas Ertl,
  Filip Sadlo, and Marc~Alexander Schweitzer.
\newblock Extraction of fragments and waves after impact damage in
  particle-based simulations.
\newblock In \emph{Meshfree Methods for Partial Differential Equations VIII},
  pages 17--34. Springer, 2017.

\bibitem[Dienes(1979)]{dienes1979analysis}
John~K Dienes.
\newblock On the analysis of rotation and stress rate in deforming bodies.
\newblock \emph{Acta mechanica}, 32\penalty0 (4):\penalty0 217--232, 1979.

\bibitem[Dipasquale et~al.(2014)Dipasquale, Zaccariotto, and
  Galvanetto]{dipasquale2014crack}
Daniele Dipasquale, Mirco Zaccariotto, and Ugo Galvanetto.
\newblock Crack propagation with adaptive grid refinement in 2d peridynamics.
\newblock \emph{International Journal of Fracture}, 190\penalty0
  (1-2):\penalty0 1--22, 2014.

\bibitem[Diyaroglu et~al.(2015)Diyaroglu, Oterkus, Oterkus, and
  Madenci]{diyaroglu2015peridynamics}
C~Diyaroglu, E~Oterkus, S~Oterkus, and Erdogan Madenci.
\newblock Peridynamics for bending of beams and plates with transverse shear
  deformation.
\newblock \emph{International Journal of Solids and Structures}, 69:\penalty0
  152--168, 2015.

\bibitem[Duong et~al.(2017)Duong, Roohbakhshan, and Sauer]{duong2017new}
Thang~X Duong, Farshad Roohbakhshan, and Roger~A Sauer.
\newblock A new rotation-free isogeometric thin shell formulation and a
  corresponding continuity constraint for patch boundaries.
\newblock \emph{Computer Methods in applied Mechanics and engineering},
  316:\penalty0 43--83, 2017.

\bibitem[Flanagan and Taylor(1987)]{flanagan1987accurate}
D.P. Flanagan and L.M. Taylor.
\newblock An accurate numerical algorithm for stress integration with finite
  rotations.
\newblock \emph{Computer methods in applied mechanics and engineering},
  62\penalty0 (3):\penalty0 305--320, 1987.

\bibitem[Fuselier and Wright(2013)]{fuselier2013high}
Edward~J Fuselier and Grady~B Wright.
\newblock A high-order kernel method for diffusion and reaction-diffusion
  equations on surfaces.
\newblock \emph{Journal of Scientific Computing}, 56\penalty0 (3):\penalty0
  535--565, 2013.

\bibitem[Gerstle et~al.(2007)Gerstle, Sau, and Silling]{gerstle2007peridynamic}
Walter Gerstle, Nicolas Sau, and Stewart Silling.
\newblock Peridynamic modeling of concrete structures.
\newblock \emph{Nuclear engineering and design}, 237\penalty0 (12-13):\penalty0
  1250--1258, 2007.

\bibitem[Gross et~al.(2020)Gross, Trask, Kuberry, and
  Atzberger]{gross2020meshfree}
Ben~J Gross, Nathaniel~A Trask, Paul Kuberry, and Paul~J Atzberger.
\newblock Meshfree methods on manifolds for hydrodynamic flows on curved
  surfaces: {A Generalized Moving Least-Squares (GMLS)} approach.
\newblock \emph{Journal of Computational Physics}, 409:\penalty0 109340, 2020.

\bibitem[Guan et~al.(2011)Guan, Chi, Chen, Slawson, and Roth]{guan2011semi}
Pai-Chen Guan, Sheng-Wei Chi, Jiun-Shyan Chen, Thomas~R Slawson, and Michael~J
  Roth.
\newblock Semi-{Lagrangian} reproducing kernel particle method for
  fragment-impact problems.
\newblock \emph{International Journal of Impact Engineering}, 38\penalty0
  (12):\penalty0 1033--1047, 2011.

\bibitem[Ha and Bobaru(2010)]{ha2010studies}
Youn~Doh Ha and Florin Bobaru.
\newblock Studies of dynamic crack propagation and crack branching with
  peridynamics.
\newblock \emph{International Journal of Fracture}, 162\penalty0 (1):\penalty0
  229--244, 2010.

\bibitem[Hallquist et~al.(2006)]{hallquist2006ls}
John~O Hallquist et~al.
\newblock Ls-dyna theory manual.
\newblock \emph{Livermore software Technology corporation}, 3:\penalty0 25--31,
  2006.

\bibitem[Hashin(1980)]{hashin1980failure}
Zvi Hashin.
\newblock Failure criteria for unidirectional fiber composites.
\newblock \emph{Journal of Applied Mechanics}, 47\penalty0 (2):\penalty0
  329--334, 1980.

\bibitem[Hillman et~al.(2020)Hillman, Pasetto, and
  Zhou]{hillman2020generalized}
Michael Hillman, Marco Pasetto, and Guohua Zhou.
\newblock Generalized reproducing kernel peridynamics: unification of local and
  non-local meshfree methods, non-local derivative operations, and an
  arbitrary-order state-based peridynamic formulation.
\newblock \emph{Computational Particle Mechanics}, 7\penalty0 (2):\penalty0
  435--469, 2020.

\bibitem[Hoppe et~al.(1992)Hoppe, DeRose, Duchamp, McDonald, and
  Stuetzle]{hoppe1992surface}
Hugues Hoppe, Tony DeRose, Tom Duchamp, John McDonald, and Werner Stuetzle.
\newblock Surface reconstruction from unorganized points.
\newblock In \emph{Proceedings of the 19th annual conference on computer
  graphics and interactive techniques}, pages 71--78, 1992.

\bibitem[Hu et~al.(2020)Hu, Feng, Li, Sheng, and Zhang]{hu2020numerical}
Yumeng Hu, Guoqing Feng, Shaofan Li, Weijia Sheng, and Chaoyi Zhang.
\newblock Numerical modelling of ductile fracture in steel plates with
  non-ordinary state-based peridynamics.
\newblock \emph{Engineering Fracture Mechanics}, 225:\penalty0 106446, 2020.

\bibitem[Hughes et~al.(2005)Hughes, Cottrell, and
  Bazilevs]{hughes2005isogeometric}
Thomas~JR Hughes, John~A Cottrell, and Yuri Bazilevs.
\newblock Isogeometric analysis: {CAD}, finite elements, {NURBS}, exact
  geometry and mesh refinement.
\newblock \emph{Computer methods in applied mechanics and engineering},
  194\penalty0 (39-41):\penalty0 4135--4195, 2005.

\bibitem[Jafarzadeh et~al.(2018)Jafarzadeh, Chen, and
  Bobaru]{jafarzadeh2018peridynamic}
Siavash Jafarzadeh, Ziguang Chen, and Florin Bobaru.
\newblock Peridynamic modeling of intergranular corrosion damage.
\newblock \emph{Journal of The Electrochemical Society}, 165\penalty0
  (7):\penalty0 C362, 2018.

\bibitem[Javili et~al.(2019)Javili, Morasata, Oterkus, and
  Oterkus]{javili2019peridynamics}
Ali Javili, Rico Morasata, Erkan Oterkus, and Selda Oterkus.
\newblock Peridynamics review.
\newblock \emph{Mathematics and Mechanics of Solids}, 24\penalty0
  (11):\penalty0 3714--3739, 2019.

\bibitem[Johnson and Cook(1985)]{johnson1985fracture}
Gordon~R Johnson and William~H Cook.
\newblock Fracture characteristics of three metals subjected to various
  strains, strain rates, temperatures and pressures.
\newblock \emph{Engineering fracture mechanics}, 21\penalty0 (1):\penalty0
  31--48, 1985.

\bibitem[Kamensky et~al.(2019)Kamensky, Behzadinasab, Foster, and
  Bazilevs]{kamensky2019peridynamic}
David Kamensky, Masoud Behzadinasab, John~T Foster, and Yuri Bazilevs.
\newblock Peridynamic modeling of frictional contact.
\newblock \emph{Journal of Peridynamics and Nonlocal Modeling}, 1\penalty0
  (2):\penalty0 107--121, 2019.

\bibitem[Kelly(2021)]{kelly2021solid}
Piaras Kelly.
\newblock Mechanics lecture notes: {Engineering} solid mechanics – small
  strain.
\newblock
  \url{http://homepages.engineering.auckland.ac.nz/~pkel015/SolidMechanicsBooks/Part_II/06_PlateTheory/06_PlateTheory_Complete.pdf},
  February 2021.

\bibitem[Kiendl et~al.(2009)Kiendl, Bletzinger, Linhard, and
  W{\"u}chner]{kiendl2009isogeometric}
Josef Kiendl, Kai-Uwe Bletzinger, Johannes Linhard, and Roland W{\"u}chner.
\newblock Isogeometric shell analysis with {Kirchhoff--Love} elements.
\newblock \emph{Computer Methods in Applied Mechanics and Engineering},
  198\penalty0 (49-52):\penalty0 3902--3914, 2009.

\bibitem[Kiendl et~al.(2010)Kiendl, Bazilevs, Hsu, W{\"u}chner, and
  Bletzinger]{kiendl2010bending}
Josef Kiendl, Yuri Bazilevs, Ming-Chen Hsu, Roland W{\"u}chner, and Kai-Uwe
  Bletzinger.
\newblock The bending strip method for isogeometric analysis of
  {Kirchhoff--Love} shell structures comprised of multiple patches.
\newblock \emph{Computer Methods in Applied Mechanics and Engineering},
  199\penalty0 (37-40):\penalty0 2403--2416, 2010.

\bibitem[Kramer et~al.(2019)Kramer, Jones, Mostafa, Ravaji,
  et~al.]{kramer2019third}
Sharlotte~LB Kramer, Amanda Jones, Ahmed Mostafa, Babak Ravaji, et~al.
\newblock {The third Sandia Fracture Challenge}: predictions of ductile
  fracture in additively manufactured metal.
\newblock \emph{International Journal of Fracture}, 218\penalty0 (1):\penalty0
  5--61, 2019.

\bibitem[Lai et~al.(2013)Lai, Liang, and Zhao]{lai2013local}
Rongjie Lai, Jiang Liang, and Hong-Kai Zhao.
\newblock A local mesh method for solving pdes on point clouds.
\newblock \emph{Inverse Problems \& Imaging}, 7\penalty0 (3):\penalty0 737,
  2013.

\bibitem[Lancaster and Salkauskas(1981)]{lancaster1981surfaces}
Peter Lancaster and Kes Salkauskas.
\newblock Surfaces generated by moving least squares methods.
\newblock \emph{Mathematics of computation}, 37\penalty0 (155):\penalty0
  141--158, 1981.

\bibitem[Lapczyk and Hurtado(2007)]{lapczyk2007progressive}
Ireneusz Lapczyk and Juan~A Hurtado.
\newblock Progressive damage modeling in fiber-reinforced materials.
\newblock \emph{Composites Part A: Applied Science and Manufacturing},
  38\penalty0 (11):\penalty0 2333--2341, 2007.

\bibitem[Leng et~al.(2019)Leng, Tian, and Foster]{leng2019super}
Yu~Leng, Xiaochuan Tian, and John~T Foster.
\newblock Super-convergence of reproducing kernel approximation.
\newblock \emph{Computer Methods in Applied Mechanics and Engineering},
  352:\penalty0 488--507, 2019.

\bibitem[Leng et~al.(2020)Leng, Tian, Trask, and
  Foster]{leng2020asymptotically}
Yu~Leng, Xiaochuan Tian, Nathaniel~A Trask, and John~T Foster.
\newblock Asymptotically compatible reproducing kernel collocation and meshfree
  integration for the peridynamic {Navier} equation.
\newblock \emph{Computer Methods in Applied Mechanics and Engineering},
  370:\penalty0 113264, 2020.

\bibitem[Leng et~al.(2021)Leng, Tian, Trask, and
  Foster]{leng2021asymptotically}
Yu~Leng, Xiaochuan Tian, Nathaniel~A Trask, and John~T Foster.
\newblock Asymptotically compatible reproducing kernel collocation and meshfree
  integration for nonlocal diffusion.
\newblock \emph{SIAM Journal on Numerical Analysis}, 59\penalty0 (1):\penalty0
  88--118, 2021.

\bibitem[Leung et~al.(2011)Leung, Lowengrub, and Zhao]{leung2011grid}
Shingyu Leung, John Lowengrub, and Hongkai Zhao.
\newblock A grid based particle method for solving partial differential
  equations on evolving surfaces and modeling high order geometrical motion.
\newblock \emph{Journal of Computational Physics}, 230\penalty0 (7):\penalty0
  2540--2561, 2011.

\bibitem[Li et~al.(2011)Li, Wierzbicki, Sutton, Yan, and Deng]{li2011mixed}
Yaning Li, Tomasz Wierzbicki, Michael~A Sutton, Junhui Yan, and Xiaomin Deng.
\newblock Mixed mode stable tearing of thin sheet {AI 6061-T6} specimens:
  {Experimental} measurements and finite element simulations using a modified
  {Mohr-Coulomb} fracture criterion.
\newblock \emph{International Journal of Fracture}, 168\penalty0 (1):\penalty0
  53--71, 2011.

\bibitem[Littlewood(2010)]{littlewood2010simulation}
David~J Littlewood.
\newblock Simulation of dynamic fracture using peridynamics, finite element
  modeling, and contact.
\newblock In \emph{ASME 2010 International Mechanical Engineering Congress and
  Exposition}, pages 209--217. American Society of Mechanical Engineers, 2010.

\bibitem[Littlewood et~al.(2016)Littlewood, Silling, and
  Demmie]{littlewood2016identification}
David~J Littlewood, Stewart~A Silling, and Paul~N Demmie.
\newblock Identification of fragments in a meshfree peridynamic simulation.
\newblock In \emph{ASME International Mechanical Engineering Congress and
  Exposition}, volume 50633, page V009T12A071. American Society of Mechanical
  Engineers, 2016.

\bibitem[Love(1888)]{love1888xvi}
Augustus~EH Love.
\newblock The small free vibrations and deformation of a thin elastic shell.
\newblock \emph{Philosophical Transactions of the Royal Society of
  London.(A.)}, 179:\penalty0 491--546, 1888.

\bibitem[Macdonald et~al.(2013)Macdonald, Merriman, and
  Ruuth]{macdonald2013simple}
Colin~B Macdonald, Barry Merriman, and Steven~J Ruuth.
\newblock Simple computation of reaction--diffusion processes on point clouds.
\newblock \emph{Proceedings of the National Academy of Sciences}, 110\penalty0
  (23):\penalty0 9209--9214, 2013.

\bibitem[Macneal and Harder(1985)]{macneal1985proposed}
Richard~H Macneal and Robert~L Harder.
\newblock A proposed standard set of problems to test finite element accuracy.
\newblock \emph{Finite elements in analysis and design}, 1\penalty0
  (1):\penalty0 3--20, 1985.

\bibitem[Madenci and Oterkus(2014)]{madenci2014peridynamic}
Erdogan Madenci and Erkan Oterkus.
\newblock Peridynamic theory.
\newblock In \emph{Peridynamic Theory and Its Applications}, pages 19--43.
  Springer, 2014.

\bibitem[Madenci et~al.(2016)Madenci, Barut, and Futch]{madenci2016peridynamic}
Erdogan Madenci, Atila Barut, and Michael Futch.
\newblock Peridynamic differential operator and its applications.
\newblock \emph{Computer Methods in Applied Mechanics and Engineering},
  304:\penalty0 408--451, 2016.

\bibitem[Matzenmiller et~al.(1995)Matzenmiller, Lubliner, and
  Taylor]{matzenmiller1995constitutive}
Anton Matzenmiller, Jacob Lubliner, and Robert~L Taylor.
\newblock A constitutive model for anisotropic damage in fiber-composites.
\newblock \emph{Mechanics of materials}, 20\penalty0 (2):\penalty0 125--152,
  1995.

\bibitem[Maurin et~al.(2018)Maurin, Greco, Coox, Vandepitte, and
  Desmet]{maurin2018isogeometric}
Florian Maurin, Francesco Greco, Laurens Coox, Dirk Vandepitte, and Wim Desmet.
\newblock Isogeometric collocation for {Kirchhoff--Love} plates and shells.
\newblock \emph{Computer Methods in Applied Mechanics and Engineering},
  329:\penalty0 396--420, 2018.

\bibitem[Mehrmashhadi et~al.(2019)Mehrmashhadi, Wang, and
  Bobaru]{mehrmashhadi2019uncovering}
Javad Mehrmashhadi, Longzhen Wang, and Florin Bobaru.
\newblock Uncovering the dynamic fracture behavior of {PMMA} with peridynamics:
  {The} importance of softening at the crack tip.
\newblock \emph{Engineering Fracture Mechanics}, 219:\penalty0 106617, 2019.

\bibitem[Mindlin(1951)]{mindlin1951influence}
Raymond~D Mindlin.
\newblock Influence of rotatory inertia and shear on flexural motions of
  isotropic, elastic plates.
\newblock \emph{Journal of Applied Mechanics}, 18:\penalty0 31--38, 1951.

\bibitem[Mohammadi et~al.(2019)Mohammadi, Dehghan, Khodadadian, and
  Wick]{mohammadi2019generalized}
Vahid Mohammadi, Mehdi Dehghan, Amirreza Khodadadian, and Thomas Wick.
\newblock Generalized moving least squares and moving kriging least squares
  approximations for solving the transport equation on the sphere.
\newblock \emph{arXiv preprint arXiv:1904.05831}, 2019.

\bibitem[Nguyen and Oterkus(2019)]{nguyen2019peridynamics}
Cong~Tien Nguyen and Selda Oterkus.
\newblock Peridynamics for the thermomechanical behavior of shell structures.
\newblock \emph{Engineering Fracture Mechanics}, 219:\penalty0 106623, 2019.

\bibitem[Nikravesh and Gerstle(2018)]{nikravesh2018improved}
Siavash Nikravesh and Walter Gerstle.
\newblock Improved state-based peridynamic lattice model including elasticity,
  plasticity and damage.
\newblock \emph{Computer Modeling in Engineering \& Sciences}, 116\penalty0
  (3):\penalty0 323--347, 2018.

\bibitem[Oterkus et~al.(2020)Oterkus, Madenci, and
  Oterkus]{oterkus2020peridynamic}
Erkan Oterkus, Erdogan Madenci, and Selda Oterkus.
\newblock Peridynamic shell membrane formulation.
\newblock \emph{Procedia Structural Integrity}, 28:\penalty0 411--417, 2020.

\bibitem[Ouchi et~al.(2015)Ouchi, Katiyar, York, Foster, and
  Sharma]{ouchi2015fully}
Hisanao Ouchi, Amit Katiyar, Jason York, John~T Foster, and Mukul~M Sharma.
\newblock A fully coupled porous flow and geomechanics model for fluid driven
  cracks: a peridynamics approach.
\newblock \emph{Computational Mechanics}, 55\penalty0 (3):\penalty0 561--576,
  2015.

\bibitem[O’Grady and Foster(2014)]{o2014peridynamic}
James O’Grady and John~T Foster.
\newblock Peridynamic plates and flat shells: A non-ordinary, state-based
  model.
\newblock \emph{International Journal of Solids and Structures}, 51\penalty0
  (25-26):\penalty0 4572--4579, 2014.

\bibitem[Parks et~al.(2012)Parks, Littlewood, Mitchell, and
  Silling]{parks2012peridigm}
Michael~L Parks, David~J Littlewood, John~A Mitchell, and Stewart~A Silling.
\newblock Peridigm users’ guide v1. 0.0.
\newblock \emph{SAND Report}, 7800, 2012.

\bibitem[Piret(2012)]{piret2012orthogonal}
C{\'e}cile Piret.
\newblock The orthogonal gradients method: A radial basis functions method for
  solving partial differential equations on arbitrary surfaces.
\newblock \emph{Journal of Computational Physics}, 231\penalty0 (14):\penalty0
  4662--4675, 2012.

\bibitem[Pressley(2010)]{pressley2010elementary}
Andrew~N Pressley.
\newblock \emph{Elementary differential geometry}.
\newblock Springer Science \& Business Media, 2010.

\bibitem[Ravi-Chandar and Knauss(1984)]{ravi1984experimental}
K.~Ravi-Chandar and Wolfgang~G Knauss.
\newblock An experimental investigation into dynamic fracture: {IV. On} the
  interaction of stress waves with propagating cracks.
\newblock \emph{International Journal of Fracture}, 26\penalty0 (3):\penalty0
  189--200, 1984.

\bibitem[Reissner(1945)]{reissner1945effect}
Eric Reissner.
\newblock The effect of transverse shear deformation on the bending of elastic
  plates.
\newblock \emph{Journal of Applied Mechanics}, 12:\penalty0 69--77, 1945.

\bibitem[Rokkam et~al.(2019)Rokkam, Gunzburger, Brothers, Phan, and
  Goel]{rokkam2019nonlocal}
Srujan Rokkam, Max Gunzburger, Michael Brothers, Nam Phan, and Kishan Goel.
\newblock A nonlocal peridynamics modeling approach for corrosion damage and
  crack propagation.
\newblock \emph{Theoretical and Applied Fracture Mechanics}, 101:\penalty0
  373--387, 2019.

\bibitem[Sarego et~al.(2016)Sarego, Le, Bobaru, Zaccariotto, and
  Galvanetto]{sarego2016linearized}
Giulia Sarego, Quang~V Le, Florin Bobaru, Mirco Zaccariotto, and Ugo
  Galvanetto.
\newblock Linearized state-based peridynamics for 2-d problems.
\newblock \emph{International Journal for Numerical Methods in Engineering},
  108\penalty0 (10):\penalty0 1174--1197, 2016.

\bibitem[Schraml et~al.(2005)Schraml, Meyer, Kleponis, and
  Kimsey]{schraml2005simulating}
Stephen~J Schraml, Hubert~W Meyer, David~S Kleponis, and Kent~D Kimsey.
\newblock Simulating the formation and evolution of behind armor debris fields.
\newblock In \emph{2005 Users Group Conference (DOD-UGC'05)}, pages 215--221.
  IEEE, 2005.

\bibitem[Shankar and Wright(2018)]{shankar2018mesh}
Varun Shankar and Grady~B Wright.
\newblock Mesh-free semi-{Lagrangian} methods for transport on a sphere using
  radial basis functions.
\newblock \emph{Journal of Computational Physics}, 366:\penalty0 170--190,
  2018.

\bibitem[Shankar et~al.(2018)Shankar, Narayan, and Kirby]{shankar2018rbf}
Varun Shankar, Akil Narayan, and Robert~M Kirby.
\newblock Rbf-loi: Augmenting radial basis functions (rbfs) with least
  orthogonal interpolation (loi) for solving pdes on surfaces.
\newblock \emph{Journal of Computational Physics}, 373:\penalty0 722--735,
  2018.

\bibitem[Silling(2000)]{silling2000reformulation}
Stewart~A Silling.
\newblock Reformulation of elasticity theory for discontinuities and long-range
  forces.
\newblock \emph{Journal of the Mechanics and Physics of Solids}, 48\penalty0
  (1):\penalty0 175--209, 2000.

\bibitem[Silling(2017)]{silling2017stability}
Stewart~A Silling.
\newblock Stability of peridynamic correspondence material models and their
  particle discretizations.
\newblock \emph{Computer Methods in Applied Mechanics and Engineering},
  322:\penalty0 42--57, 2017.

\bibitem[Silling and Askari(2005)]{silling2005meshfree}
Stewart~A Silling and Ebrahim Askari.
\newblock A meshfree method based on the peridynamic model of solid mechanics.
\newblock \emph{Computers \& structures}, 83\penalty0 (17-18):\penalty0
  1526--1535, 2005.

\bibitem[Silling and Bobaru(2005)]{silling2005peridynamic}
Stewart~A Silling and Florin Bobaru.
\newblock Peridynamic modeling of membranes and fibers.
\newblock \emph{International Journal of Non-Linear Mechanics}, 40\penalty0
  (2-3):\penalty0 395--409, 2005.

\bibitem[Silling and Lehoucq(2010)]{silling2010peridynamic}
Stewart~A Silling and Richard~B Lehoucq.
\newblock Peridynamic theory of solid mechanics.
\newblock \emph{Advances in applied mechanics}, 44:\penalty0 73--168, 2010.

\bibitem[Silling et~al.(2007)Silling, Epton, Weckner, Xu, and
  Askari]{silling2007peridynamic}
Stewart~A Silling, M~Epton, O~Weckner, Ji~Xu, and Ebrahim Askari.
\newblock Peridynamic states and constitutive modeling.
\newblock \emph{Journal of Elasticity}, 88\penalty0 (2):\penalty0 151--184,
  2007.

\bibitem[Simo and Taylor(1985)]{simo1985consistent}
Juan~C Simo and Robert~L Taylor.
\newblock Consistent tangent operators for rate-independent elastoplasticity.
\newblock \emph{Computer methods in applied mechanics and engineering},
  48\penalty0 (1):\penalty0 101--118, 1985.

\bibitem[Suchde(2021)]{suchde2021meshfree}
Pratik Suchde.
\newblock A meshfree lagrangian method for flow on manifolds.
\newblock \emph{International Journal for Numerical Methods in Fluids},
  93\penalty0 (6):\penalty0 1871--1894, 2021.

\bibitem[Suchde and Kuhnert(2019)]{suchde2019fully}
Pratik Suchde and J{\"o}rg Kuhnert.
\newblock A fully lagrangian meshfree framework for pdes on evolving surfaces.
\newblock \emph{Journal of Computational Physics}, 395:\penalty0 38--59, 2019.

\bibitem[Talamini and Radovitzky(2017)]{talamini2017parallel}
Brandon~L Talamini and Ra{\'u}l Radovitzky.
\newblock A parallel discontinuous {Galerkin}/cohesive-zone computational
  framework for the simulation of fracture in shear-flexible shells.
\newblock \emph{Computer Methods in Applied Mechanics and Engineering},
  317:\penalty0 480--506, 2017.

\bibitem[Taylor and Steigmann(2015)]{taylor2015two}
Michael Taylor and David~J Steigmann.
\newblock A two-dimensional peridynamic model for thin plates.
\newblock \emph{Mathematics and Mechanics of Solids}, 20\penalty0 (8):\penalty0
  998--1010, 2015.

\bibitem[Tian and Du(2014)]{tian2014asymptotically}
Xiaochuan Tian and Qiang Du.
\newblock Asymptotically compatible schemes and applications to robust
  discretization of nonlocal models.
\newblock \emph{SIAM Journal on Numerical Analysis}, 52\penalty0 (4):\penalty0
  1641--1665, 2014.

\bibitem[Timoshenko and Woinowsky-Krieger(1959)]{timoshenko1959theory}
Stephen~P Timoshenko and Sergius Woinowsky-Krieger.
\newblock \emph{Theory of plates and shells}.
\newblock McGraw-hill, 1959.

\bibitem[Torres-S{\'a}nchez et~al.(2020)Torres-S{\'a}nchez, Santos-Oliv{\'a}n,
  and Arroyo]{torres2020approximation}
Alejandro Torres-S{\'a}nchez, Daniel Santos-Oliv{\'a}n, and Marino Arroyo.
\newblock Approximation of tensor fields on surfaces of arbitrary topology
  based on local monge parametrizations.
\newblock \emph{Journal of Computational Physics}, 405:\penalty0 109168, 2020.

\bibitem[Trask and Kuberry(2020)]{trask2020compatible}
Nathaniel~A Trask and Paul Kuberry.
\newblock Compatible meshfree discretization of surface {PDEs}.
\newblock \emph{Computational Particle Mechanics}, 7\penalty0 (2):\penalty0
  271--277, 2020.

\bibitem[Trask et~al.(2019)Trask, You, Yu, and Parks]{trask2019asymptotically}
Nathaniel~A Trask, Huaiqian You, Yue Yu, and Michael~L Parks.
\newblock An asymptotically compatible meshfree quadrature rule for nonlocal
  problems with applications to peridynamics.
\newblock \emph{Computer Methods in Applied Mechanics and Engineering},
  343:\penalty0 151--165, 2019.

\bibitem[Tupek and Radovitzky(2014)]{tupek2014extended}
Michael~R Tupek and Raul Radovitzky.
\newblock An extended constitutive correspondence formulation of peridynamics
  based on nonlinear bond-strain measures.
\newblock \emph{Journal of the Mechanics and Physics of Solids}, 65:\penalty0
  82--92, 2014.

\bibitem[Tupek et~al.(2013)Tupek, Rimoli, and Radovitzky]{tupek2013approach}
Michael~R Tupek, Julian~J Rimoli, and Raul Radovitzky.
\newblock An approach for incorporating classical continuum damage models in
  state-based peridynamics.
\newblock \emph{Computer methods in applied mechanics and engineering},
  263:\penalty0 20--26, 2013.

\bibitem[Ugural(2009)]{ugural2009stresses}
Ansel~C Ugural.
\newblock \emph{Stresses in beams, plates, and shells}.
\newblock CRC press, 2009.

\bibitem[Vazic et~al.(2020)Vazic, Oterkus, and Oterkus]{vazic2020peridynamic}
Bozo Vazic, Erkan Oterkus, and Selda Oterkus.
\newblock Peridynamic model for a {Mindlin} plate resting on a {Winkler}
  elastic foundation.
\newblock \emph{Journal of Peridynamics and Nonlocal Modeling}, pages 1--10,
  2020.

\bibitem[Vogler et~al.(2003)Vogler, Thornhill, Reinhart, Chhabildas, Grady,
  Wilson, Hurricane, and Sunwoo]{vogler2003fragmentation}
Tracy~J Vogler, Tom~F Thornhill, William~D Reinhart, Lalit~C Chhabildas,
  Dennis~E Grady, Leonard~T Wilson, Omar~A Hurricane, and Anne Sunwoo.
\newblock Fragmentation of materials in expanding tube experiments.
\newblock \emph{International journal of impact engineering}, 29\penalty0
  (1-10):\penalty0 735--746, 2003.

\bibitem[Wold et~al.(1987)Wold, Esbensen, and Geladi]{wold1987principal}
Svante Wold, Kim Esbensen, and Paul Geladi.
\newblock Principal component analysis.
\newblock \emph{Chemometrics and intelligent laboratory systems}, 2\penalty0
  (1-3):\penalty0 37--52, 1987.

\bibitem[Wu and Ren(2015)]{wu2015stabilized}
C.T. Wu and Bo~Ren.
\newblock A stabilized non-ordinary state-based peridynamics for the nonlocal
  ductile material failure analysis in metal machining process.
\newblock \emph{Computer Methods in Applied Mechanics and Engineering},
  291:\penalty0 197--215, 2015.

\bibitem[Yaghoobi and Chorzepa(2017)]{yaghoobi2017fracture}
Amin Yaghoobi and Mi~G Chorzepa.
\newblock Fracture analysis of fiber reinforced concrete structures in the
  micropolar peridynamic analysis framework.
\newblock \emph{Engineering Fracture Mechanics}, 169:\penalty0 238--250, 2017.

\bibitem[Yang et~al.(2020)Yang, Vazic, Diyaroglu, Oterkus, and
  Oterkus]{yang2020kirchhoff}
Zhenghao Yang, Bozo Vazic, Cagan Diyaroglu, Erkan Oterkus, and Selda Oterkus.
\newblock A {Kirchhoff} plate formulation in a state-based peridynamic
  framework.
\newblock \emph{Mathematics and Mechanics of Solids}, 25\penalty0 (3):\penalty0
  727--738, 2020.

\bibitem[Zhang et~al.(2021)Zhang, Li, Zhang, Peng, and
  Yan]{zhang2021peridynamic}
Qi~Zhang, Shaofan Li, A-Man Zhang, Yuxiang Peng, and Jiale Yan.
\newblock A peridynamic {Reissner-Mindlin} shell theory.
\newblock \emph{International Journal for Numerical Methods in Engineering},
  122\penalty0 (1):\penalty0 122--147, 2021.

\bibitem[Zhu et~al.(2011)Zhu, Mobasher, Rajan, and
  Peralta]{zhu2011characterization}
Deju Zhu, Barzin Mobasher, Subramaniam~D Rajan, and Pedro Peralta.
\newblock Characterization of dynamic tensile testing using aluminum alloy
  {6061-T6} at intermediate strain rates.
\newblock \emph{Journal of Engineering Mechanics}, 137\penalty0 (10):\penalty0
  669--679, 2011.

\end{thebibliography}

\end{document}